\documentclass[12pt,a4paper]{article}
\usepackage{t1enc}
\usepackage[latin1]{inputenc}
\usepackage[english]{babel}
\usepackage{amssymb}
\usepackage{amsmath}
\usepackage{graphics}
\usepackage{amsthm, amssymb}
\usepackage{amsfonts}
\usepackage{epsfig,multicol}

\newtheorem{theorem}[subsubsection]{Theorem}
\newtheorem{proposition}[subsubsection]{Proposition}
\newtheorem{lemma}[subsubsection]{Lemma}
\newtheorem{corollary}[subsubsection]{Corollary}
\newtheorem{definition}[subsubsection]{Definition}
\newtheorem{example}[subsubsection]{Example}
\newtheorem{remark}[subsubsection]{Remark}
  
\newtheorem{diagram}[subsubsection]{Fig.}

\newfont{\gothic} { ygoth scaled \magstep{1.5}}

\newcommand{\5}{\vskip 5pt}

\def\<{\langle}
\def\>{\rangle}

\begin{document}

\def\hpic #1 #2 {\mbox{$\begin{array}[c]{l} \epsfig{file=#1,height=#2}
\end{array}$}}
 
\def\vpic #1 #2 {\mbox{$\begin{array}[c]{l} \epsfig{file=#1,width=#2}
\end{array}$}}

\title{Intermediate subfactors with no extra structure.}
\author{Pinhas Grossman and Vaughan F.R. Jones
\thanks{Supported in part by NSF Grant DMS93--22675, the Marsden fund UOA520,
 and the Swiss National Science Foundation.}
} \maketitle 
\begin{abstract}
 If $N\subseteq P,Q\subseteq M$ are type II$_1$ factors with $N'\cap M =\mathbb C id$ 
 and $[M:N]<\infty$ we show that restrictions on the standard invariants of
the elementary inclusions $N\subseteq P$, $N\subseteq Q$, $P\subseteq M$ and $Q\subseteq M$
imply drastic restrictions on the indices and angles between the subfactors. In particular
we show that if these standard invariants are trivial and the conditional expectations
onto $P$ and $Q$ do not commute, then $[M:N]$ is $6$ or $6+4\sqrt 2$. In the former case
$N$ is the fixed point algebra for an outer action of $S_3$ on $M$ and the angle is
$\pi/3$, and in the latter case the angle is $\cos^{-1}(\sqrt 2-1)$ and
an example may be found in the GHJ subfactor family. The techniques of proof rely heavily on
planar algebras.

  \end{abstract}

\section{Introduction}

Let $N\subseteq M$ be II$_1$ factors with $[M:N]<\infty$.
There is a "standard invariant" for $N\subseteq M$ which we
shall describe using the planar algebra formalism of \cite{J18}. The
vector spaces $P_k$ of $N-N$ invariant vectors in the $N-N$ bimodule $\otimes^k M$  
admit an action of the operad of planar tangles
as in \cite{J18} and \cite{J21}. In more usual notation the vector space $P_k$ is
the relative commutant $N'\cap M_{k-1}$ in the tower $M_k$ of \cite{J3}.
The conditional expectation $E_N$ from  $M$ to $N$ is in
$P_2$ and generates a sub-planar algebra called the Temperley-Lieb
algebra. In \cite{BJ}, Bisch and the second author studied the planar subalgebra
of the $P_k$ generated by the conditional expectation onto a single 
intermediate subfactor $N\subseteq P\subseteq M$. The resulting planar
algebra is called the Fuss-Catalan algebra and was generalised by Bisch and the second author to a chain
of intermediate subfactors -see also \cite{La}. These planar algebras are
universal in that they are always planar subalgebras of the standard invariant
for any subfactor possessing a chain of intermediate subfactors.  
If $P_i \subseteq P_{i+1}$ is the chain, there are no restricitions on
the individual inclusions of $P_i$ in $P_{i+1}$. Moreover the existence
of the Fuss Catalan planar algebra together with a theorem of Popa in \cite{P20}
allows one to construct a "free" increasing chain where the individual
inclusions  $P_i \subseteq P_{i+1}$ have "no extra structure", i.e. their
own standard invariants are just the Temperley-Lieb algebra. Thus the 
standard invariants for the $P_i \subseteq P_{i+1}$ are "decoupled"
from the algebraic symmetries coming from the existence of a chain
of intermediate subfactors.

In \cite{SaW}, Sano and Watatani considered the angle between two subfactors
$P\subseteq M$ and $Q\subseteq M$ which we shall here define via the square
of its cosine, namely the spectrum of the positive self-adjoint operator
$E_PE_QE_P$ (on $L^2(M)$). In \cite{JX}, Feng Xu and the second author proved
that finiteness of the angle (as a substet of $[0,1]$) is equivalent to 
finiteness of the index of $P\cap Q$ in $M$. If we suppose that $P\cap Q$
is an irreducible finite index subfactor of $M$ we might expect that the
angle is "quantized", i.e. only a certain discrete countable family 
of numbers occurs-at least in a range close to $0$ and $\pi/2$. Determining
these allowed angle values is becoming a significant question in the abstract
theory of subfactors. This paper can be considered a first step in answering
that question.

In \cite{Wat3}, Watatani considered the lattice of intermediate subfactors for 
a finite index inclusion and showed that if the inclusion is irreducible
the lattice is finite. He gave some constructions which allowed him to 
realise many simple finite lattices, but even for two lattices with
only six elements,
the question of their realisation as intermediate subfactor lattices
remains entirely open.

The present paper grew out of an attempt by Dietmar Bisch and the second
author to extend the methods of \cite{BJ} to attack both the angle quantization and the 
intermediate lattice problems. The hope was to construct universal planar
algebras depending only on the lattice of intermediate subfactors, and possibly
the angles between them, and use Popa's theorem to construct subfactors 
realising the lattice and angle values. This project is probably sound but
it is hugely more difficult in the case where the lattice is not a chain
or the angles are not all $0$ or $\pi /2$. The reason is very simple-the
planar algebra generated by the conditional expectations can no longer be decoupled
from the standard invariants of the elementary subfactor inclusions in
the lattice. This is surprisingly easy to see. The spectral subspaces of
$E_PE_QE_P$ are $N-N$-bimodules contained in $P$ so that as soon as the
angle operator has a significant spectrum the subfactor $N\subseteq M$ 
must have elements in its planar algebra that are not in the Temperley-Lieb
subalgebra-a situation we shall refer to as having "extra structure" and which
we will quantify using the notion of supertransitivity introduced in \cite{J31}. 
In particular if there is no extra structure the spectrum of $E_PE_QE_P$
can consist of at most one number besides $0$ and $1$. We will call the
angle whose cosine is the square root of this number "\underline{the} angle between
$P$ and $Q$.
Or "dually" if $PQP$ is not equal to all of $M$, it is a
non-trivial $P-P$ bimodule between $P$ and $M$ so that the inclusion
$P\subseteq M$ must have extra structure. 

Thus we are led to the question-what are the possible pairs of subfactors
$P$ and $Q$ in $M$ with $P\cap Q$ a finite index irreducible subfactor
of $M$, for which the four elementary subfactors $N\subseteq P$, $N\subseteq Q$,
$P\subseteq M$ and $Q\subseteq M$ all have no extra structure? 
More properly, since we are not trying to control the isomorphism type
of the individual factors, one should ask what are the standard invariants
that arise. One situation
is rather easy to take care of: if the subfactors form a commuting cocommuting
square in the sense of \cite{SaW}, there is no obstruction-it was essentially observed
by Sano and Watatani that in this case $E_P$ and $E_Q$ generate a tensor
product of their individual Temperley-Lieb algebras. And to realise any 
$N\subseteq P$ and $N\subseteq Q$ just take the tensor product II$_1$ factors.
However if we assume that the subfactors either do not commute or do not cocommute,
we will show in this paper the following unexpected result.

\begin{theorem} Suppose $\begin{array}{ccc}
P&\subset &M \cr
\cup& &\cup \cr
N&\subset &Q
\end{array}$
 is a quadrilateral of
subfactors with $N'\cap M =\mathbb C$ , $[M:N]<\infty$ and no extra structure. 
Then either the quadrilateral is commuting or  one
of the following two cases occurs:

a) $[M:N]=6$ and $N$ is the fixed point algebra for an outer action of $S_3$
on $M$ with $P$ and $Q$ being the fixed point algebras for two transpositions
in $S_3$. In this case the angle between $P$ and $Q$ is $\pi/3$ and the
full intermediate subfactor lattice is 
\vpic{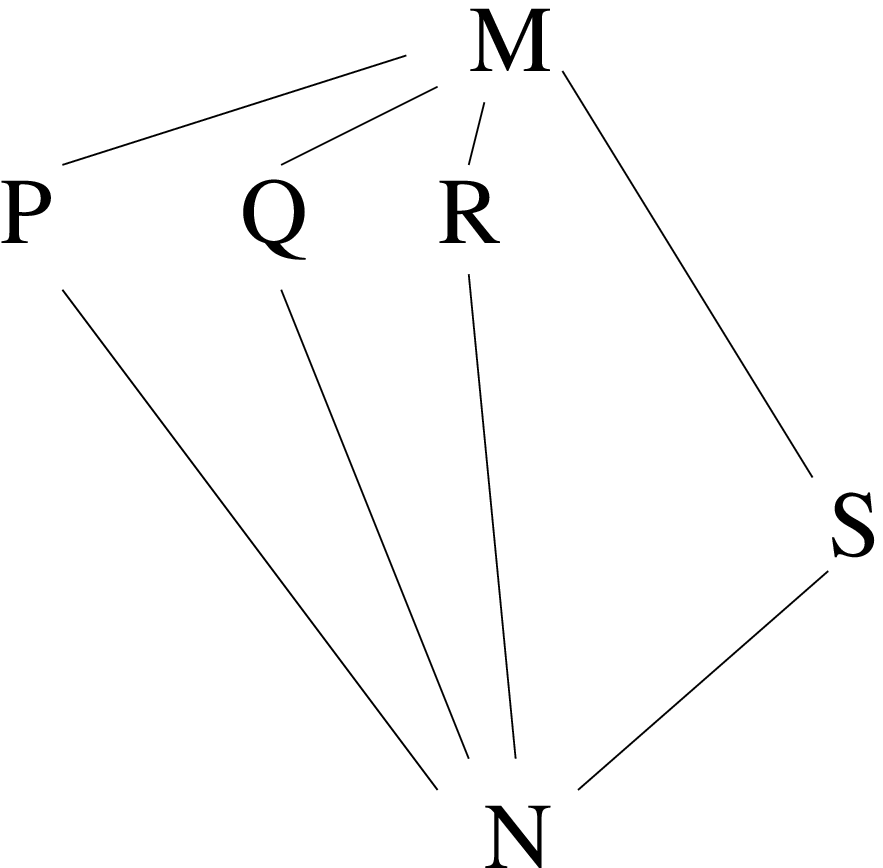} {1.0 in}
.

(Note that the dual of this quadrilateral is a commuting square.)

b) The subfactor $N$ is of depth $3$, 
$[M:N]=(2+\sqrt 2)^2$ and the planar algebra of $N\subseteq M$ is the same as
that coming from the GHJ subfactor (see \cite{GHJ}) constructed from the Coxeter graph 
$D_5$ with the distinguished vertex being the trivalent one. 
Each of the intermediate inclusions has index $2+\sqrt 2$ and the angle between
$P$ and $Q$ is $\theta=\cos^{-1} (\sqrt 2 - 1)$. The principal graph of $N\subseteq M$
is \vpic{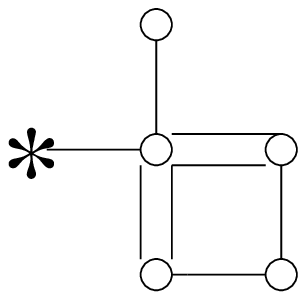} {0.7 in}  and the full intermediate subfactor lattice is
\vpic{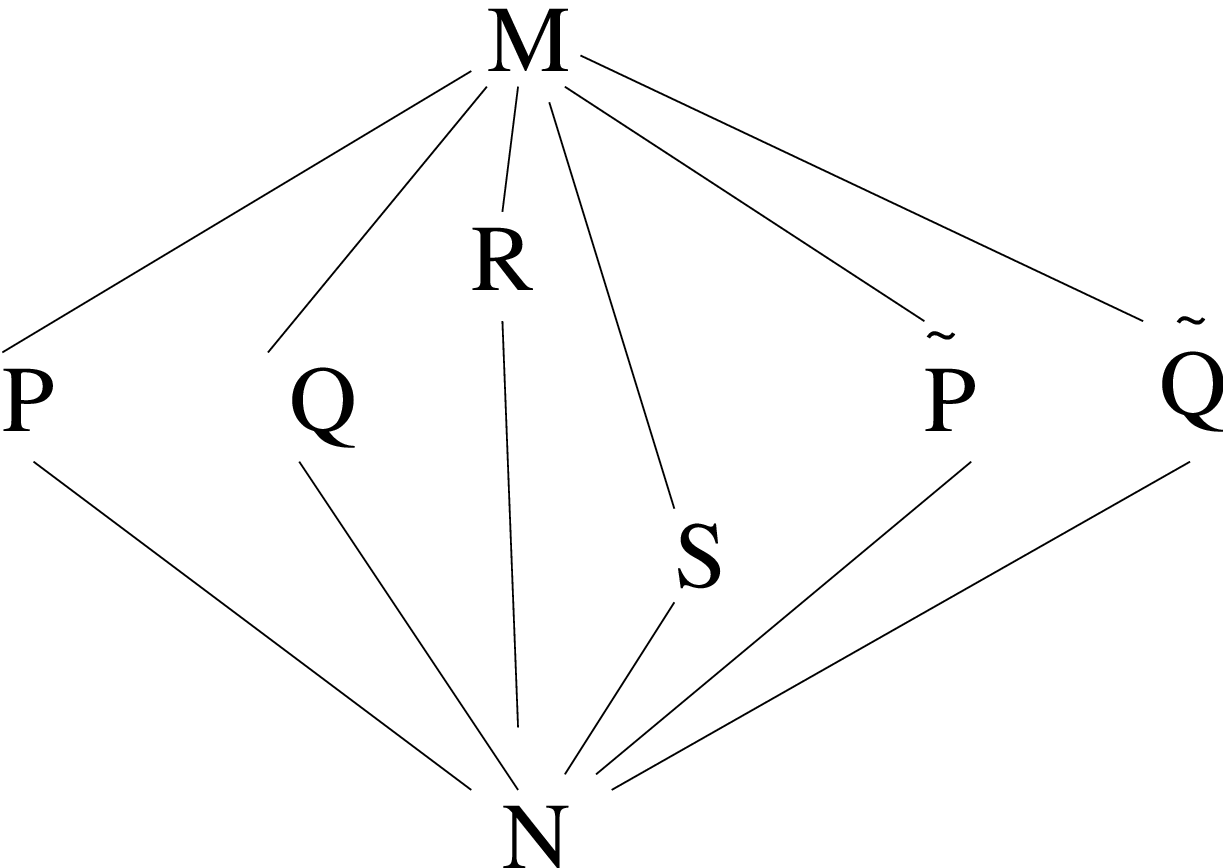} {2.0 in} \quad where the angle between $\tilde P$ and $\tilde Q$ 
is also $\theta$ but $P$ and $Q$ both commute with $\tilde P$ and $\tilde Q$. Moreover
$[M:R]=[S:N]=2$ and $M,N,R$ and $S$ form a commuting cocommuting square.
The planar algebra of $N\subseteq M$ is isomorphic to its dual-the planar algebra
of $M\subseteq M_1$.
\end{theorem} 

Note that from Ocneanu's paragroup point of view $N$ is the fixed point algebra
of an action of the paragroup given by the planar algebra on $M$. Thus if
the ambient factor $M$ is hyperfinite, Popa's theorem in \cite{P6} guarantees that
the subactors are unique up to an automorphism of $M$. Also note that
it is a consequence of the theorem that any intermediate subfactor lattice with
four elements and no extra structure is a commuting square.

Our methods rely heavily on planar algebras. Of crucial importance is the diagram
discovered by Landau for the projection onto the product $PQ$. We give a proof
of Landau's result and some general consequences.
The uniqueness of the subfactor of index $6+4\sqrt 2$ mentioned in the theorem is
proved using the "exchange relation" of \cite{La2} - the planar algebras have a very simple
skein theory in the sense of \cite{J21}. The no extra structure hypothesis necessary for
the theorem is in fact weaker than the one we have stated above. For a precise
statement of the required supertransitivity see \ref{cocommthm} and \ref{mainchance}. 

The authors would like to thank Dietmar Bisch and Zeph Landau for several
fruitful discussions concerning this paper.

\section{Background.}
\subsection{Bimodules.}

We recall some basic facts about bimodules over II$_{1}$ factors.
 The treatment follows \cite{Bs7}. For more on this, look there and in \cite{JS}.

\begin{definition}
Let $M$ be a II$_{1}$ factor. A left $M$-module is a pair $(H,\pi)$ where $H$ 
is a Hilbert space and $\pi$ is a unital normal homomorphism from $M$ into the 
algebra of bounded operators on $H$. The dimension of $H$ over $M$, denoted $dim_{M}H$, is the
 extended positive number given by the Murray-von Neumann coupling constant of $\pi(M)$.
Let $M^{OP}$ be the opposite algebra of $M$ (i.e. the algebra with the same underlying vector
 space 
but with multiplication reversed).
Then a right $M$-module is defined as a left $M^{OP}$-module. An $M-N$-bimodule is a triple 
$(H,\pi,\phi)$, where $H$ is a Hilbert space and $\pi$ and $\phi$ are normal unital homomorphisms
 from, respectively, $M$ and $N^{OP}$ into the algebra of bounded operators on $H$, such that
$\pi(M)$ and $\phi(N^{OP})$ commute. Such a bimodule is denoted by $_{M} H _{N}$
, or sometimes simply by $H$, if the action is understood. We write $m \xi n$ 
for $\pi(m) \phi(n) \xi$, where $m \in M$, $n \in N$, and $\xi \in H$.
\end{definition}

There are obvious notions of submodules and direct sums. An $M-N$ bimodule is in particular both 
a left $M$-module and a left $N^{OP}$-module.

\begin{definition}
An $M-N$-bimodule is bifinite if $dim_{M}H$ and $dim_{N^{OP}}H$ are both finite.
\end{definition}

All bimodules will be assumed to be bifinite.

\begin{definition}
Let $_{M} H^{1} _{N}$ and $_{M} H^{2} _{N}$ be bimodules.  The intertwiner space, 
denoted $Hom_{M-N}(H^{1},H^{2})$, is the subspace of bounded operators from $H^{1}$ to $H^{2}$ 
consisting of those operators which commute with the bimodule action: $T \in Hom_{M-N}(H^{1},H^{2})$ 
iff $T(m \xi n) = m (T\xi) n$ for all $m \in M$, $n \in N$, $\xi \in H^{1}$.
\end{definition}

\begin{example}
Let $M$ be a II$_{1}$ factor. $L^{2}(M)$ is the Hilbert space completion of $M$ with respect to the
 inner product induced
by the unique normalized trace on $M$. Then $L^{2}(M)$ is an $M-M$ bimodule, and the left and 
right actions are simply the continuous extensions of ordinary left and right multiplication in $M$.
If $P$ and $Q$ are subfactors of $M$, then $L^{2}(M)$ is a $P-Q$-bimodule by restriction, 
and it is bifinite iff the indices $[M:P]$ and $[M:Q]$ are finite.
\end{example}

\begin{definition}
Let $_{M} H _{N}$ be a bimodule. There is a dense subspace $H^{0}$ of $H$, called the space 
of bounded vectors, defined by the rule that 
$\xi \in H^{0}$ iff the map $m \mapsto m \xi$ extends to a bounded operator from $L^{2}(M)$ to $H$.
 To each pair of bounded vectors $(\xi, \eta)$ there is associated an element of $M$, 
denoted \\$\langle \xi,\eta \rangle_{M}$, determined by the relation $\langle m \xi, \eta \rangle =tr(m \langle \xi,\eta \rangle_{M})$.
\end{definition}

\begin{remark}
It is in fact also true that $\xi \in H^{0}$ iff the map $n \mapsto n \xi$ extends to a 
bounded operator from $L^{2}(N)$ to $H$.
\end{remark}

\begin{remark}
Let $M$ be a II$_{1}$ factor, and consider $L^{2}(M)$ as a bimodule over a pair of finite 
index subfactors as in Example 2.0.5. Then $L^{2}(M)^{0}$ is simply the image of $M$ in $L^{2}(M)$.
\end{remark}

\begin{definition}
Let $_{M} H _{N}$ and $_{N} K _{P}$ be bimodules. There is an $M-P$ bimodule, 
denoted $(_{M} H _{N})  \otimes_{N} (_{N} K _{P})$, called the relative tensor 
product, or fusion, of 
$_{M} H _{N}$ and $_{N} K _{P}$, which is characterized by the following 
property: there is a surjective linear map
from the algebraic tensor product $H^{0} \odot K^{0}$ to 
$((_{M} H _{N}) \otimes_{N} ( _{N} K _{P}))^{0}$, $\xi \otimes \eta \mapsto \xi \otimes_{N} \eta$ 
satisfying the following three conditions:
\vskip 2pt
\noindent (i) $\xi n \otimes_{N} \eta=\xi \otimes_{N} n \eta$ \\
(ii) $m(\xi \otimes_{N} \eta)p=(m \xi)\otimes_{N}(\eta p)$\\(iii)
 $\langle \xi \otimes_{N} \eta,\xi ' \otimes_{N} \eta ' \rangle_{M}=\langle \xi \langle \eta,\eta '\rangle_{M},\xi '\rangle_{M} $\\
 (for all $m \in M$, $n \in N$, and $p \in P$).
\end{definition}

\begin{remark}
Among the properties enjoyed by fusion are: it is distributive over direct sums, it is associative,
 and it is multiplicative in dimension:
 $\dim_{M}( {}_{M}H _{N}  \otimes_{N} \ {}_{N}K_{P})=(\dim_{M} H )(\dim_{N} K$).
\end{remark}

Let $N \subset M$ be an inclusion of II$_1$ factors with finite index. $L^{2}(N)$ can be
 identified with a subspace of $L^{2}(M)$. Let $e_{1}$ denote the corresponding projection 
on $L^{2}(M)$, and let $M_{1}$ be the von Neumann algebra generated by $M$ and $e_{1}$. Then
 $M_{1}$ is a II$_{1}$ factor and $[M_{1}:M]=[M:N]$. This procedure is called the basic 
construction \cite{J3}. Recall that the space of bounded vectors in $L^{2}(M)$ can be identified with $M$.
 $e_{1}$ leaves this space invariant, inducing a trace-preserving expectation of $M$ onto $N$.

Iterating the basic construction we get a sequence of projections $e_{1},e_{2}...$, and 
a tower of algebras $M_{-1} \subset M_{0} \subset M_{1} \subset M_{2} \subset ...$,
 where $M_{-1}=N$, $M_{0}=M$, $e_{k}$ is the projection onto
 $L^{2}(M_{k-2})$ in $B(L^{2}(M_{k-1}))$, and $M_{k}$ is the von Neumann algebra generated
 by $M_{k-1}$ and $e_{k}$, for $k \geq 1$. Restricting the tower to those elements which 
commute with $N$, we get a tower of finite dimensional algebras, called the tower of relative
 commutants $N'\cap M_k$. 

Each $L^{2}(M_{k})$, $k \geq 0$ is an $N-N$ bimodule, and   

\begin{proposition}
$L^{2}(M_{k}) \cong L^{2}(M) \otimes_{N} ... \otimes_{N} L^{2}(M)$, ($k+1$ factors), as an
$N-N$ bimodule. Moreover, $Hom_{N-N}L^{2}(M_{k}) \cong N' \cap M_{2k+1}$. So an $N-N$ bimodule decomposition of 
$L^{2}(M) \otimes_{N} ... \otimes_{N} L^{2}(M)$,($k+1$ factors), corresponds to
 a decomposition of the identity in $N' \cap M_{2k+1}$. Under this correspondence
projections in $N' \cap M_{2k+1}$ correspond to submodules of 
$L^{2}(M) \otimes_{N} ... \otimes_{N} L^{2}(M)$,($k+1$ factors), minimal projections 
correspond to irreducible submodules (those which have no proper nonzero closed submodules),
 and simple summands of $N'\cap M_{2k+1}$ to equivalence classes of irreducible submodules.
\end{proposition}

\subsection{Planar algebras.}
In \cite{J18} a diagrammatic calculus was introduced as an axiomatisation and 
calculational tool for the standard invariant of a finite index subfactor. 
We will use it heavily in this paper so we recall some of the essentials 
The specific uses of the calculus in this paper make possible 
a couple of simplifying conventions for the pictures.

In its most recent formulation in \cite{J31} a planar algebra {\gothic {P}}  consists of 
vector spaces $P_k^\pm$ indexed by a non-negative integer $n$
and a sign $+$ or $-$. For the planar algebra of a
subfactor $N\subseteq M$, $P_k^+=N'\cap M_{k-1}$ and $P_k^-=M'\cap M_k$.
The vector spaces $P_k^\pm$ form an algebra over the planar operad which means
that there are multilinear maps between the $P_k^\pm$ indexed by planar
tangles. A planar $k$-tangle $\mathfrak T$ consists of \\
(\romannumeral 1) The unit disc $D_0$ with $2k$ distinguished boundary points,
a finite number of disjoint interior discs $D_j\subset D_0$ for $k\geq 0$,
each with an even number of distinguished boundary points, and smooth disjoint curves
called \underline{strings}, in $D_0$
meeting the $D_j$ exactly (transversally) in the distinguished boundary points.\\
(\romannumeral 2) A black and white \underline{shading} of the regions of  $\mathfrak T$
whose boundaries consist of the strings and the boundaries of the discs between
the distinguished points. Regions of the tangle whose closures intersect are
shaded different colours.\\
(\romannumeral 3) For each disc $D_j$ there is a choice of distinguished boundary
interval between two adjacent distinguished points.\\

An example of a $k$-tangle is shown below (where we have used a $*$ near a
boundary interval to indicate the chosen one).\\

\vpic{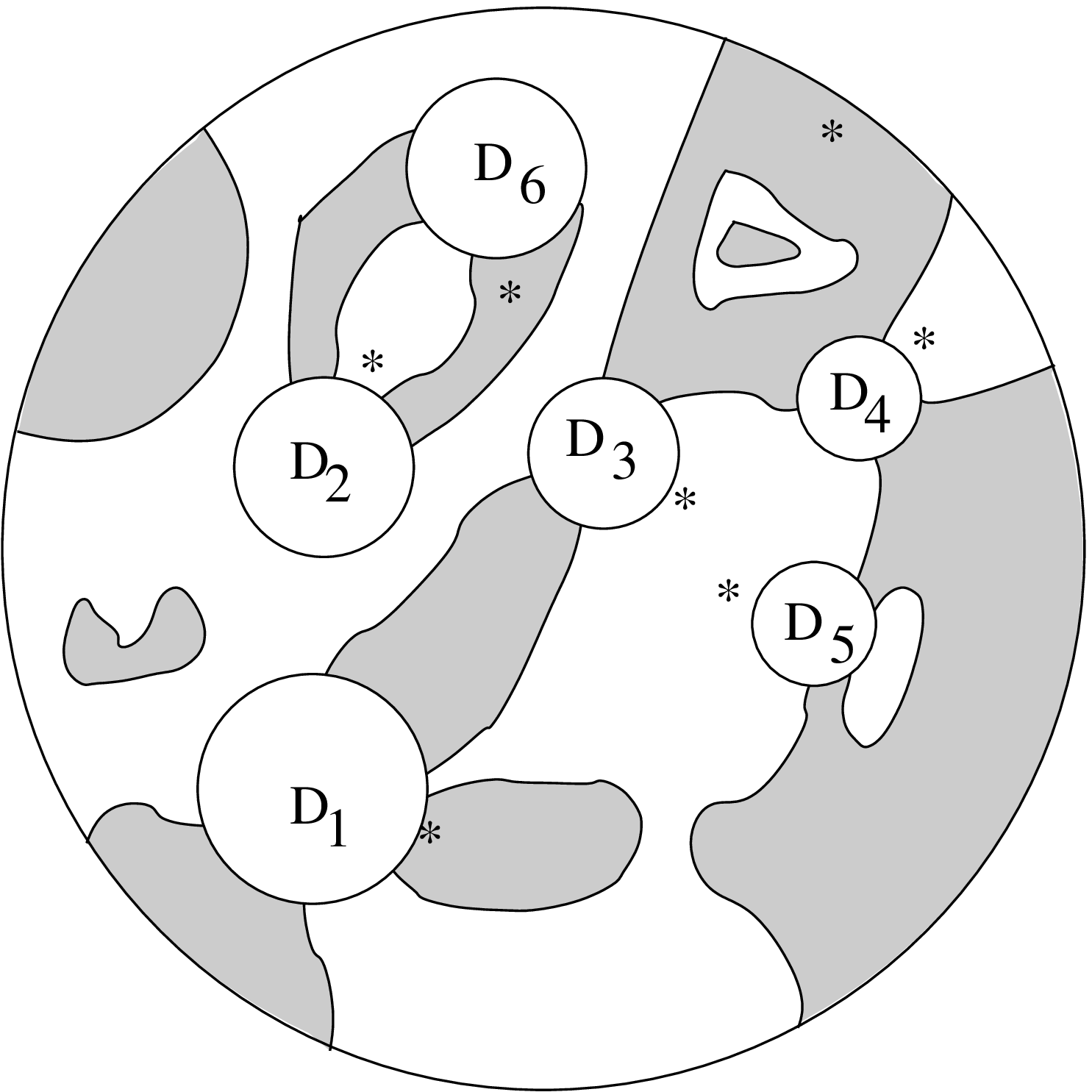} {2.5in}

The multilinear map associated to the $k$-tangle $\mathfrak T$ goes from the product of 
the $P_{k_j}^\pm$ for each internal disc where $k$ is half the number of boundary
points for $D_j$ to $P_k^\pm$, the signs being chosen $+$ if the distinguised boundary
region is shaded and $-$ if it is unshaded. The axioms of a planar algebra are
that the multilinear maps be independent of isotopies globably fixing the boundary of $D_0$
and compatible with
the gluing of tangles in a sense made clear in \cite{J18}.  To indicate the value of a tangle
on its arguments one simply inserts the arguments in the internal discs. This notation
for an element of $P_k$ is called a labelled tangle. For instance for 
$x\in P_3^+$, $y \in P_2^+$ and $a,b,c,d \in P_2^-$, the labelled
tangle below is the element of $P_4^+$ obtained by applying
the multilinear map defined by the tangle above to the elements
$x,y,a,b,c,d$ according to the discs in which they are placed.\\

\vpic{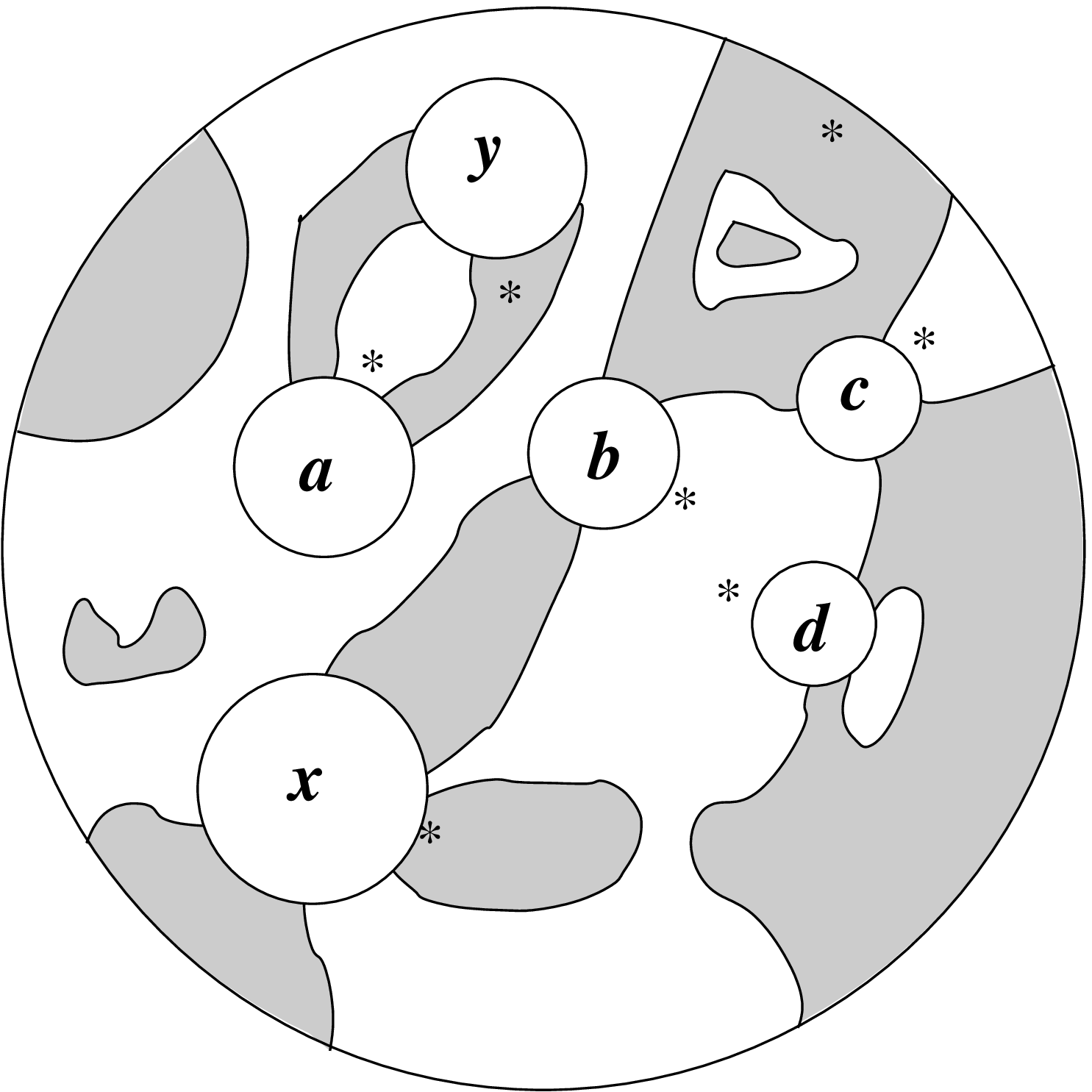} {2.5in}

We refer to \cite{J18} for details on the meaning of various tangles and the
fact that the standard invariant of a subfactor is a planar algebra. Recall
that closed strings in a tangle can always be removed, each one counting
for a multiplicative factor of the parameter $\delta$ which is the square
root of the index for a subfactor planar algebra.
%A subfactor planar algebra is one with several properties, first
%$\dim P_0^\pm =1$ so that one may identify these vector spaces with $\mathbb C$,
%the tangle with no internal discs and no strings being the identity. 

 To avoid both the shading and the marking of the distinguished
boundary interval we will adopt the following convention:

All discs will be replaced by rectangles called ``boxes". The distinguished boundary
points will be on a pair of opposite edges of each box, called the top and bottom.
Labels will be well chose letters which have
a top and bottom which will allow us to say which edge is top and which is bottom.
 The distinguished interval will be supposed shaded and always
be between the first and second strings on the top of a box. This allows us to 
put elements of $P_k^+$ in the boxes so we further adopt the convention that
if $t$ is in $P_k^-$ it will be inserted at right angles to the top-bottom axis
of its rectangle, which is to be interpreted as an internal disc whose distinguished
(unshaded) interval is the edge of the rectangle to which the letter points upwards.

Thus the two diagrams below, with $a\in P_2^-$ and $b,e\in P_2^+$
represent the same thing according to our convention.\\

\vpic{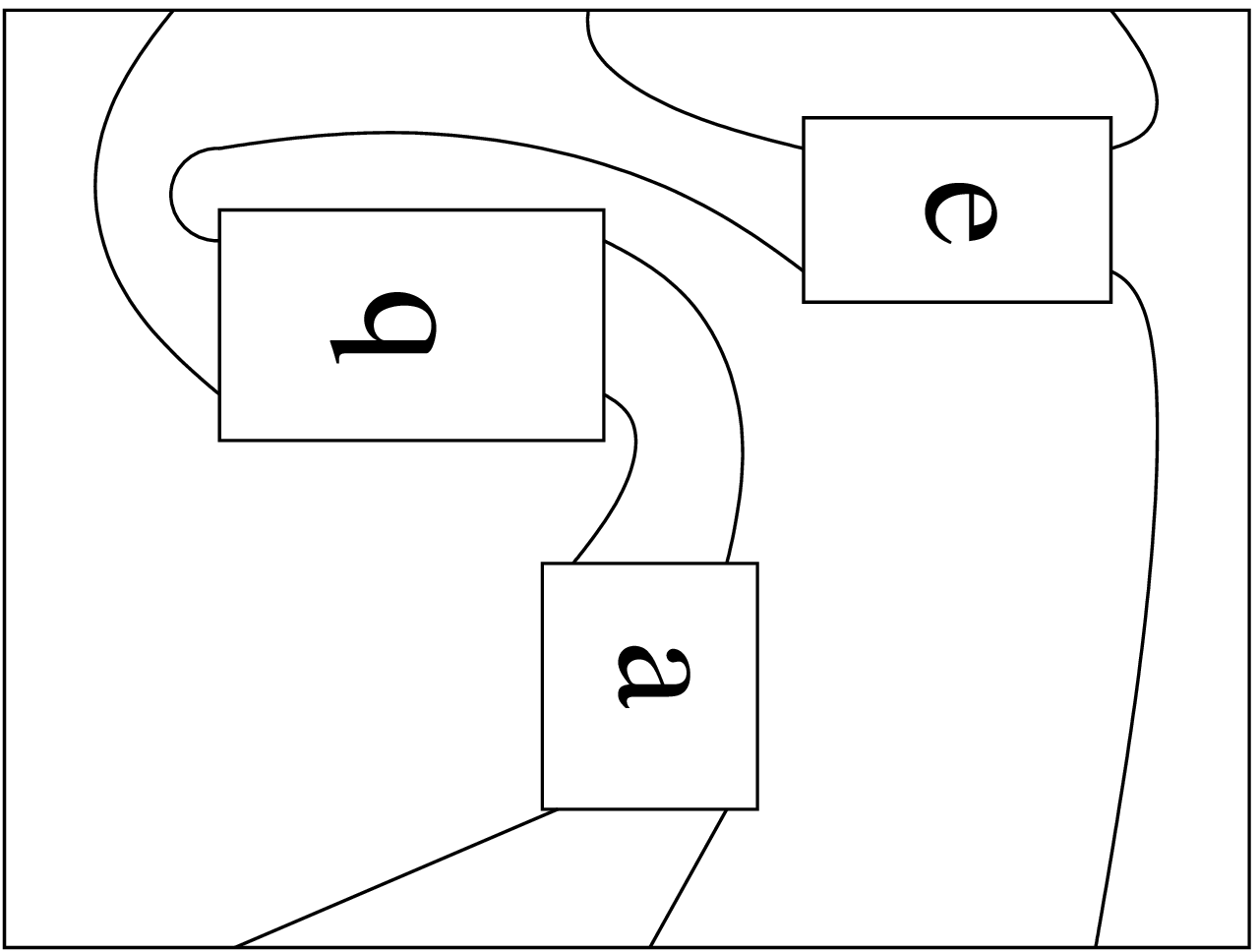} {1.6in}   \hbox{\quad}  \vpic{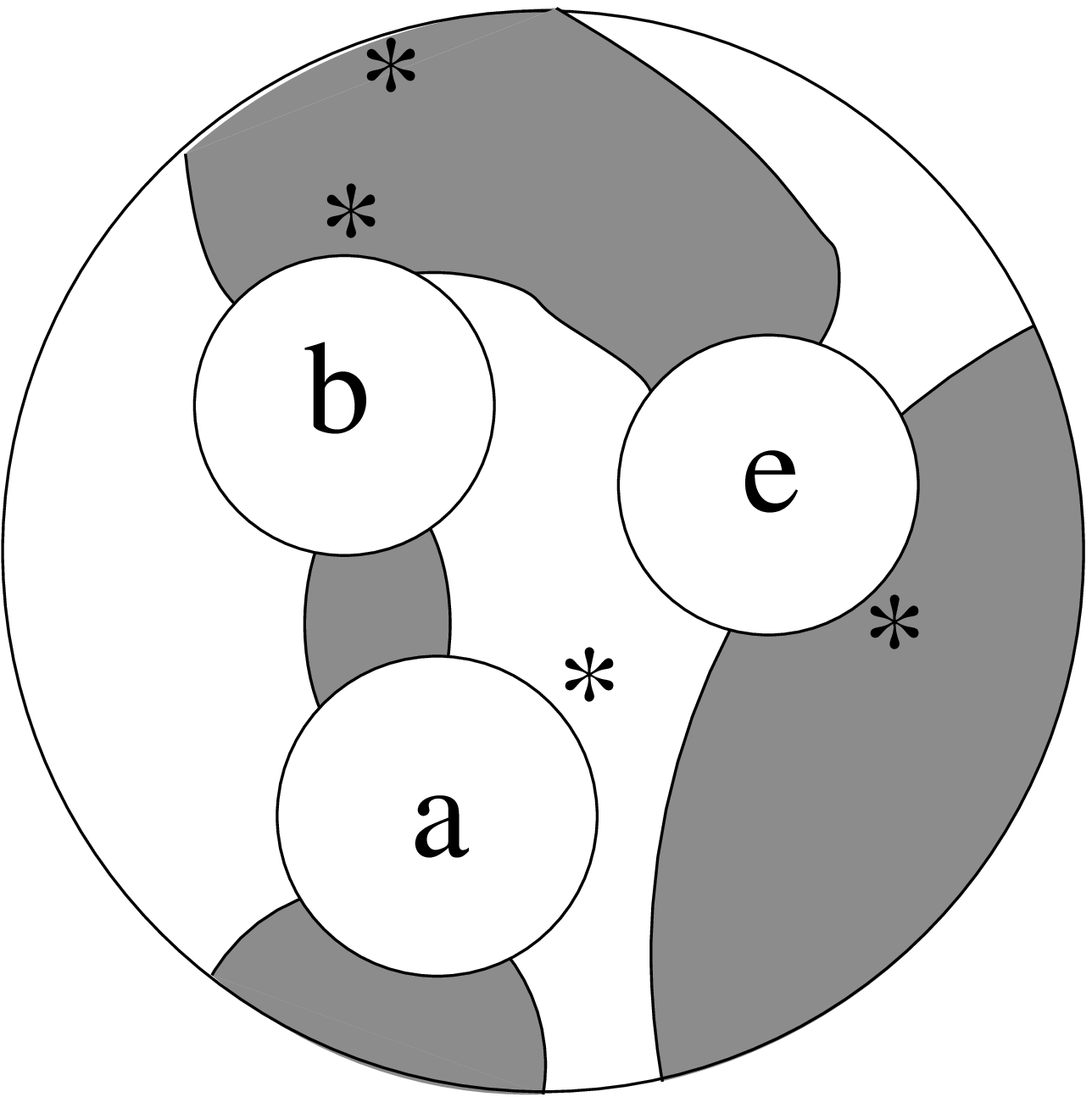} {1.6in}

 We will also from time to time simplify the diagrams further by suppressing
the outside rectangle. Thus both the above pictures are the same as
the one below:\\
\hbox{\qquad \qquad \qquad \qquad} \vpic{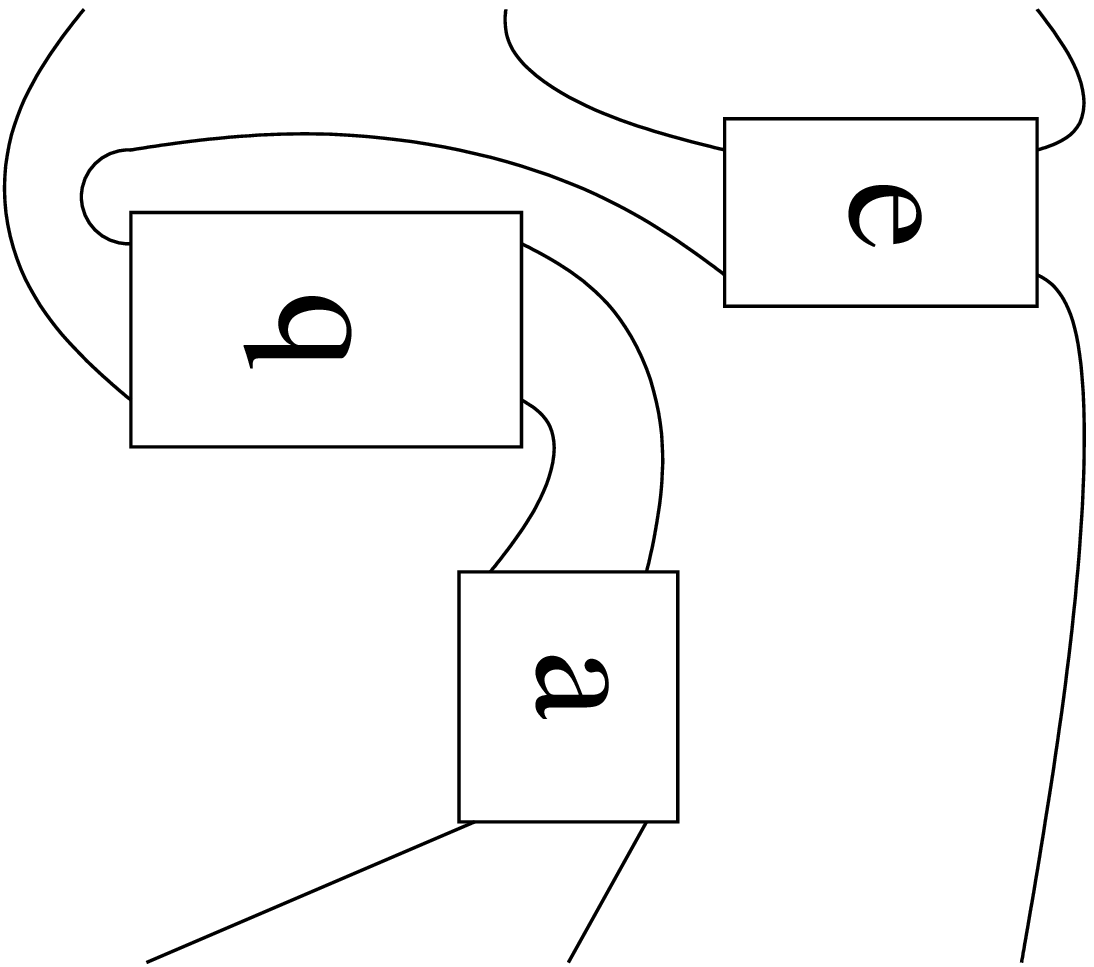} {1.6in}

\section{Generalities.}

\subsection{Multiplication.}
Let $N \subset M$ be an irreducible inclusion of $II_{1}$ factors with finite index, and
 suppose that $P$ and $Q$ are intermediate subfactors of this inclusion. Following Sano and Watatani \cite{SaW} , we say that $N \subset P,Q \subset M$ is a quadrilateral if $(P \cup Q)''=M$ and $P \cap Q=N$. (There is no real loss of generality here since in any case we can restrict our attention to $(P \cup Q)''$ and $P \cap Q$.) If $N \subset P,Q \subset M$ is a quadrilateral, there is also a dual quadrilateral $M \subset \bar{P},\bar{Q} \subset M_{1}$, where $M_{1}$ as usual is the extension of $M$ by $e_{N}$ and $\bar{P}$ and $\bar{Q}$ are the extensions of $M$ by $e_{P}$ and $e_{Q}$ respectively.

\begin{proposition} \label{mulmap}
The multiplication map from $P \otimes_{N} Q$ to $M$ extends to a surjective $N-N$ bimodule intertwiner from $L^{2}(P) \otimes_{N} L^{2}(Q)$ to $L^{2}(PQ)$.
\end{proposition}

\begin{proof}
The extension is simply (a scalar multiple of) the composition $$L^{2}(P) \otimes_{N} L^{2}(Q) \rightarrow L^{2}(M) \otimes_{N} L^{2}(M) \cong L^{2}(M_{1}) \rightarrow L^{2}(M)$$ where the first map is the tensor product of the inclusions and the last map is the conditional expectation $e_{M}$. 
\end{proof}

\begin{corollary} \label{isosub}
$L^{2}(PQ)$ is isomorphic as an $N-N$ bimodule to a submodule of $L^{2}(P) \otimes_{N} L^{2}(Q)$.
\end{corollary}

\begin{remark} \label{pqpmulmap}
In a similar way we can define a multiplication map from $\otimes_{N}^{k} (L^{2}(P) \otimes_{N} L^{2}(Q))$ to $L^{2}((PQ)^{k})$ for any $k$.
\end{remark}

\subsection{Comultiplication.}
    
Let $N \subset M$ be an irreducible inclusion of $II_{1}$ factors with finite index.
 (Irreducible here means that $N' \cap M \cong \mathbb{C}$). Consider also the dual
 inclusion $M \subset M_{1}$. 
\begin{proposition}
The first relative commutants $N' \cap M_{1}$ and $M' \cap M_{2}$ have the same vector space dimension, and the map $\phi:N' \cap M_{1} \rightarrow M' \cap M_{2}$, $a \mapsto \delta^{3} E_{M'}(ae_{2}e_{1})$, is a linear isomorphism with inverse $a \mapsto \delta^{3} E_{M_{1}}(ae_{1}e_{2})$, where $\delta=[M:N]^{\frac{1}{2}}$ and $E_{M'}$,$E_{M_{1}}$ are the conditional expectations of $N' \cap M_{2}$ onto $M' \cap M_{2}$ and $N' \cap M_{1}$ respectively.  
\end{proposition}

\begin{remark}
In the planar picture, $\phi$ is simply \vpic{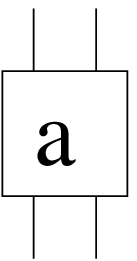} {0.3 in}  $\mapsto$ \vpic{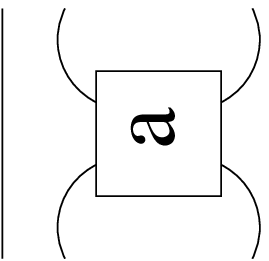} {0.6 in}  .
\end{remark}

Pulling back the multiplication in $M' \cap M_{2}$ via $\phi$ gives a second 
multiplication 
on $N' \cap M_{1}$. Using the inner product given by the trace one may identify
the vector space $N' \cap M_{1}$ with its dual and the second multiplication
may thus be pulled back to the dual. If the depth of the subfactor is 2 this
multiplication on the dual induces a Hopf
algebra structure on $N' \cap M_{1}$, but in general this does not work. We will
abuse terminology by calling the second multiplication on $N' \cap M_{1}$
``comultiplication" and use the symbol $\circ$ for it.

\begin{definition}
Let $a$ and $b$ be elements of $N' \cap M_{1}$. Then 
$a \circ b =\phi^{-1}(\phi(b)\phi(a))=\delta^{9} 
E_{M_{1}}(E_{M'}(be_{2}e_{1})E_{M'}(ae_{2}e_{1})e_{1}e_{2})$.
 Diagrammatically, $a \circ b$ is given by the picture \vpic{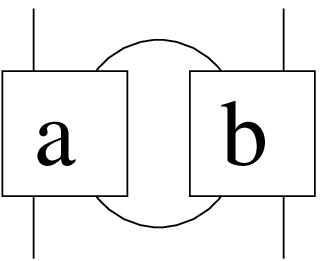} {0.5 in}  .
\end{definition}

\begin{remark}
Dually, there is a comultiplication on $M' \cap M_{2}$, also denoted by $\circ$, defined by pulling back the multipication via $\phi^{-1}$. Consequently, all of the formulas involving comultiplication have dual versions. 
\end{remark}

If $V$ is a vector subspace of $M$ which is closed under left and right multiplication 
by elements of $N$, then the closure of the image of $V$ in $L^{2}(M)$ is an $N-N$
 submodule of $L^{2}(M)$, denoted by $L^{2}(V)$, and the corresponding projection 
(necessarily in $N' \cap M_{1}$) by $e_{V}$. Conversely, any projection
 in $N' \cap M_{1}$ is of the form $e_{V}$ for a strongly closed $N-N$
 submodule $V$ of $M$, which is self-adjoint and multiplicatively closed iff $V$ is an 
intermediate subfactor. Bisch has shown that if $e$ is an arbitrary projection 
in $N' \cap M_{1}$, then $e$ is of the form $e_P$ for an intermediate subfactor $P$ iff $e$ commutes with the modular conjugation on $L^{2}(M)$ and $e \circ e$ is a scalar multiple of
 $e$(\cite{Bs3}). In that case we call $e$ a biprojection. 

NOTATION: In the planar algebra pictures, discs will be labelled simply by $V$ instead of $e_V$. 

Note that the set of biprojections inherits a partial order fron the intermediate subfactor lattice.
 In particualr, $e_{1}=e_{N}=e_{N}e_{P}$ for any intermediate subfactor $P$. We will need the following 
technical result:

\begin{lemma} \label{Ebarp}
Suppose that $e_{P} \in N' \cap M_{1}$ is an intermediate subfactor projection. 
Let $\bar{P}=<M,e_{P}> \subseteq M_1$, with corresponding projection $e_{\bar{P}}$ in $L^{2}(M_{1})$. 
Then $E_{\bar{P}}(e_{1})=\delta^{-2}tr(e_{P})^{-1}e_{P}$.
\end{lemma}

\begin{proof}
Let $x$ and $y$ be elements of $M$. Then \\
$tr(e_{1}xe_{P}y)=tr(e_{1}e_{P}xe_{P}y)=tr(e_{1}E_{P}(x)y)=$ $\delta^{-2}tr(E_{P}(x)y)=$ \\
$\delta^{-2}tr(e_{P})^{-1}tr(e_{P}E_{P}(x)y)=
\delta^{-2}tr(e_{P})^{-1}tr(e_{P}xe_{P}y)$. 
\end{proof}

\begin{lemma} \label{epbarpic}
With notation as above, $\phi(e_{P})=$\vpic{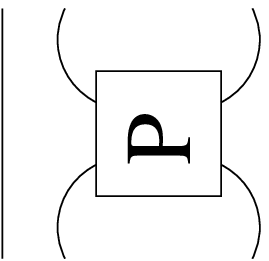} {0.5 in}   $=\delta tr(e_{P}) e_{\bar{P}}$. 
\end{lemma}

\begin{proof}
We have $\phi^{-1}(e_{\bar{P}})=\delta^{3}E_{M_{1}}(e_{\bar{P}}e_{1}e_{2})=
\delta^{3}E_{M_{1}}(e_{\bar{P}}e_{1}e_{\bar{P}}e_{2})=\\
\delta^{3}E_{M_{1}}(E_{\bar{P}}(e_{1})e_{2})=\delta^{3}E_{\bar{P}}(e_{1})E_{M_{1}}(e_{2})=
\delta^{3}\delta^{-2}tr(e_{P})^{-1}e_{P}\delta^{-2}=\\ \delta^{-1}tr(e_{P})^{-1}e_{P}$. 
Applying $\phi$ to both sides of the equation gives the result.
\end{proof}

Let $P$ and $Q$ be intermediate subfactors of the inclusion $N \subset M$ with corresponding 
projections $e_{P}$ and $e_{Q}$. Then $\bar{P}=\langle M,e_{P}\rangle$ and 
$\bar{Q}=\langle M,e_{Q}\rangle$ are intermediate
 subfactors of the dual inclusion $M \subset M_{1}$, with corresponding projections
 $e_{\bar{P}}$ and $e_{\bar{Q}}$ in $M' \cap M_{2}$.

The following result is due to Zeph Landau:

\begin{theorem} \label{landauproj}
(Landau) $e_{P} \circ e_{Q}=\delta tr (e_{P}e_{Q}) e_{PQ}$.
\end{theorem}

\begin{proof}
We have $e_{P} \circ e_{Q} = \phi^{-1}(\phi(e_{Q})\phi(e_{P}))=
\delta^{2}tr(e_{Q})tr(e_{P})\phi^{-1}(e_{\bar{Q}}e_{\bar{P}})=
\delta^{5}tr(e_{Q})tr(e_{P})E_{M_{1}}(e_{\bar{Q}}e_{\bar{P}}e_{1}e_{2})$.

By a small abuse of notation, we shall identify $M$ with its image 
in $L^{2}(M)$. Let $x \in M$. For any $a \in N' \cap M_{1}$, we have $a(x)=\delta^{2}E_{M}(axe_{1})$.
 In particular, $e_{P} \circ e_{Q}(x)= 
\delta^{5}tr(e_{P})tr(e_{Q})E_{M_{1}}(e_{\bar{Q}}e_{\bar{P}}e_{1}e_{2})(x)=\\
\delta^{7}tr(e_{P})tr(e_{Q})E_{M}(e_{\bar{Q}}e_{\bar{P}}e_{1}e_{2}xe_{1})$. 
Let $y$ be another element of $M$. Then $tr(e_{\bar{Q}}e_{\bar{P}}e_{1}e_{2}xe_{1}y)=
tr(e_{\bar{P}}e_{1}e_{\bar{P}}e_{2}xe_{2}e_{\bar{Q}}e_{1}e_{\bar{Q}}y)\\=
\delta^{-4}tr(e_{P})^{-1}tr(e_{Q})^{-1}tr(e_{P}e_{2}xe_{Q}y)$ \qquad \qquad (by \ref{Ebarp}  )\\ 
$=\delta^{-6}tr(e_{P})^{-1}tr(e_{Q})^{-1}tr(e_{P}xe_{Q}y)$. \\
Thus $E_{M}(e_{\bar{Q}}e_{\bar{P}}e_{1}e_{2}xe_{1})=\delta^{-6}tr(e_{P})^{-1}tr(e_{Q})^{-1}$ and
 $e_{P} \circ e_{Q}(x)=\delta E_{M}(e_{P}xe_{Q})$.

So if $x=pq$, with $p \in P$ and $q \in Q$, then 
$e_{P} \circ e_{Q} (x) = \delta E_{M}(e_{P} x e_{Q})=
\delta E_{M}(e_{P}e_{Q})x=\delta tr (e_{P}e_{Q}) (x)$.

To finish the proof, it suffices to show that $e_{P} \circ e_{Q}$ vanishes on the orthogonal 
complement of $L^{2}(PQ)$, or equivalently, that if $tr(xqp)=0$ for all $p \in P, q \in Q$ 
then $E_{M}(e_{P} x e_{Q})=0$. So suppose $tr(xqp)=0$ for all $p \in P, q \in Q$. 
Let $\{p_{i}\}$,$\{q_{j}\}$ be Pimsner-Popa bases over $N$ for $P$ and $Q$, respectively. 
Then $e_{P}=\sum{p_{i} e_{1} p_{i}^*}$ and $e_{Q}=\sum{q_{i} e_{N} q_{i}^*}$. Suppose $y \in M$. 
For any $i$, $j$, we have: $tr(p_{i} e_{1} p_{i}^* x q_{j} e_{1} q_{j}^* y)=
\delta^{-2} tr(p_{i}^* x q_{i} E_{N}(q_{j}^* y p_{i}))=
\delta^{-2} tr (x q_{i} E_{N}(q_{j}^* y p_{i}) p_{i}^*)=0$.
This implies that $E_{M}(e_{P}xe_{Q})=0$.    
\end{proof}

\begin{corollary} \label{epqM=PQ}
 $e_{PQ}(M)=PQ$. 
\end{corollary}
\begin{proof}
From the proof we have $e_{PQ}(x)=tr(e_P e_Q)^{-1} E_{M}(e_Pxe_Q)$. Moreover 
$e_{P}=\sum{p_{i} e_{1} p_{i}^*}$ and $e_{Q}=\sum{q_{i} e_{N} q_{i}^*}$ with the same
notation as in \ref{landauproj}. We 
 see that $e_{PQ}(M)\subset PQ$. 
\end{proof}

\begin{corollary}\label{pqclosed}
$PQ$ is strongly closed in $M$.
\end{corollary}
\begin{proof}
Since $e_{PQ}$ is strongly continuous and the identity
on $PQ$, $e_{PQ}$ is the identity on the strong closure of $PQ$.
\end{proof}

\begin{lemma}
Let $a \in N' \cap M_{1}$. Then \vpic{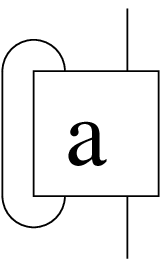} {0.3 in}  $=$ \vpic{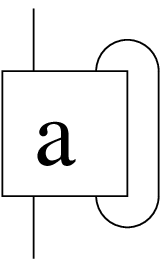} {0.3 in}  $=\delta tr(a)$.
\end{lemma}

\begin{proof}
Labelled tangles with two boundary points are elements of $N' \cap M$, which by irreducibility must be 
scalars. \\So \vpic{a_right} {0.3 in}  $=tr($ \vpic{a_right} {0.3 in}  $)= \delta ^{-1} tr($ \vpic{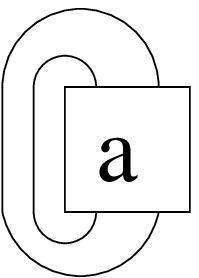} {0.3 in}  $ )=\delta tr(a)$.
\end{proof}

One corollary of \ref{landauproj} is the following multiplication formula:

\begin{proposition} \label{mulfor}
$tr(e_{PQ}) tr(e_{P}e_{Q}) = tr(e_{P})tr(e_{Q})$.
\end{proposition}

\begin{proof}
$\delta tr(e_{PQ})tr(e_{P}e_{Q})=tr($ \vpic{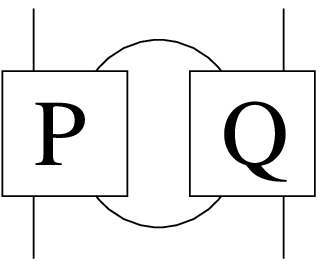} {0.6 in} $)
=\delta^{-2}$ \vpic{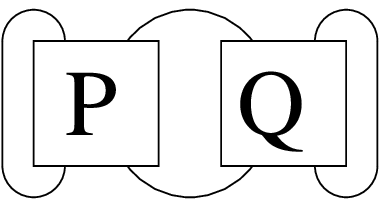} {0.8 in} \\  $=\delta^{-2}($\vpic{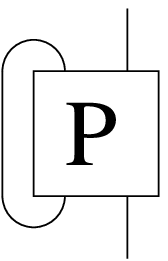} {0.3 in}  $)
($ \vpic{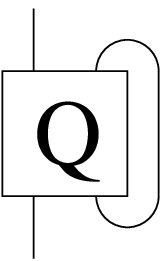} {0.3 in}  $)($ \vpic{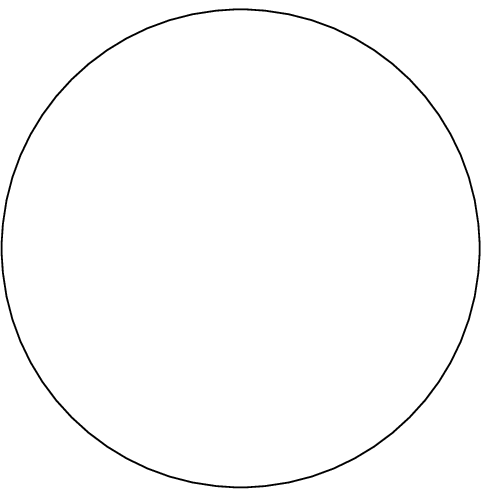} {0.2 in} $)=\delta tr(e_{P}) tr(e_{Q})$.
\end{proof}

\begin{corollary} \label{treqpeqltrepq}
$tr(e_{PQ})=tr(e_{QP})$.
\end{corollary}

And another trace formula:

\begin{lemma} \label{trfor2}
$tr(e_{P}e_{Q}) = \displaystyle \frac{1}{dim_{M}L^{2}(\bar{P}\bar{Q})}$.
\end{lemma}

\begin{proof}
${\displaystyle \frac{1}{dim_{M}L^{2}(\bar{P}\bar{Q})}=
\displaystyle \frac{1}{\delta^{2}tr(e_{\bar{P}\bar{Q}})}=
\displaystyle \frac{tr(e_{\bar{P}}e_{\bar{Q}})}{\delta^{2}tr(e_{\bar{P}})tr(e_{\bar{Q}})}}$,
 by the (dual version of) the multiplication formula. 
By \ref{epbarpic}  ${\displaystyle e_{\bar{P}}=\frac{1}{\delta tr(e_{P})} \phi(e_{P})}$,
 so that becomes ${\displaystyle \frac{1}{dim_{M}L^{2}(\bar{P}\bar{Q})}=
\frac{1}{\delta^{4}tr(e_{\bar{P}})tr(e_{\bar{Q}}) tr(e_{P}) tr(e_{Q})} tr(\phi(e_{P})\phi(e_{Q}))= 
tr(\phi(e_{P})\phi(e_{Q}))}$\\ $ =\delta^{-2} \vpic{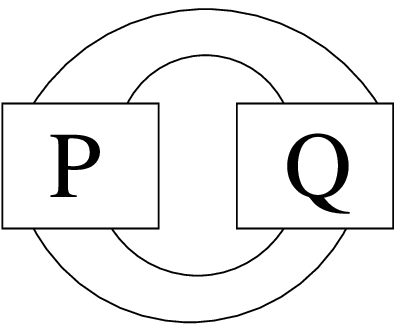} {0.6 in}   $.
On the other hand, $\vpic{cotrpcoq} {0.6 in}  =\delta^{3} tr(e_{1}(e_{P} 
\circ e_{Q})e_{1})=\delta^{4}tr(e_{P}e_{Q})tr(e_{1})=\delta^{2}tr(e_{P}e_{Q})$. 
Combining these two equations gives the result.
\end{proof}

We mention one more formula which we will need later.

\begin{lemma} \label{trepqeqp}
$tr(e_{PQ}e_{QP})=(\delta tr(e_{PQ}))^{2}tr((e_{\bar{P}}e_{\bar{Q}}e_{\bar{P}})^{2})$.
\end{lemma}

\begin{proof}
By \ref{landauproj}  , $tr(e_{PQ}e_{QP})=\displaystyle \frac{tr((e_{P} \circ e_{Q})(e_{Q} \circ e_{P}) )} 
  {(\delta tr(e_{P}e_{Q}))(\delta tr(e_{Q}e_{P}))}=
\displaystyle \frac{1}{\delta^{4} (tr(e_{P}e_{Q}))^{2}}$ \vpic{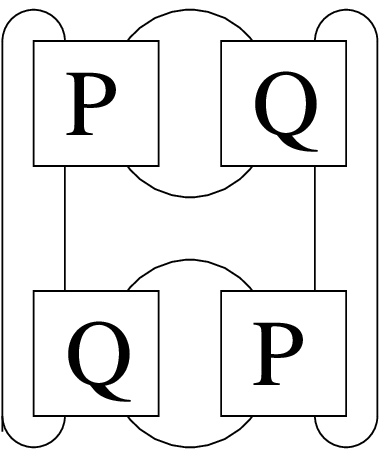} {0.5 in}  . 
On the other hand, by \ref{epbarpic}  , $tr((e_{\bar{P}}e_{\bar{Q}}e_{\bar{P}})^{2})=
tr(e_{\bar{P}}e_{\bar{Q}}e_{\bar{P}}e_{\bar{Q}})=\displaystyle 
\frac{1}{\delta^{4}tr(e_{P})^{2}tr(e_{Q})^{2}}$ \vpic{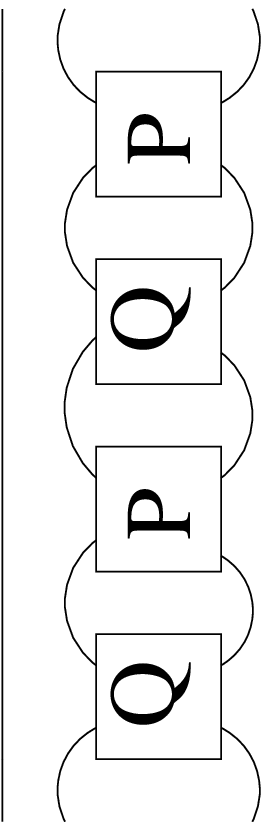} {0.4 in}   $=
\displaystyle \frac{1}{\delta^{6}tr(e_{P})^{2}tr(e_{Q})^{2}}$ \vpic {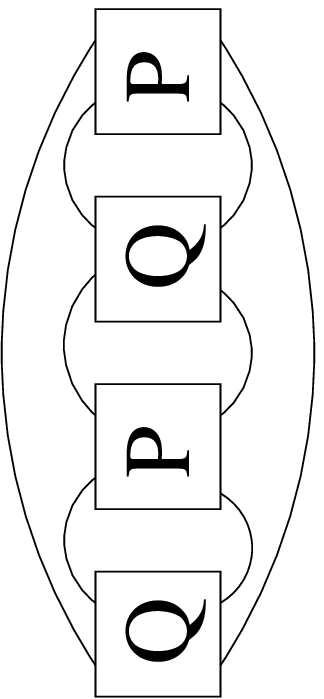} {0.5 in}   . 
By \cite{Bs3} the 2-box for a biprojection is invariant under rotation by $\pi$, so the two trace
 pictures are the same. Combining these two equations then gives 
${\displaystyle tr(e_{PQ}e_{QP})}$\\${\displaystyle =\delta^{2} \frac{tr(e_{P})^{2}tr(e_{Q})^{2}} 
{tr(e_{P}e_{Q})^{2}} 
tr((e_{\bar{P}}e_{\bar{Q}}e_{\bar{P}})^{2})}$, which by \ref{mulfor}  equals\\
 $(\delta tr(e_{PQ}))^{2} tr((e_{\bar{P}}e_{\bar{Q}}e_{\bar{P}})^{2})$.
\end{proof}

\subsection{Commuting and cocommuting quadrilaterals.}
Following Sano and Watatani \cite{SaW} we consider the condition that a quadrilateral forms a commuting square,
 which means that $e_{P}e_{Q}=e_{Q}e_{P}$. A quadrilateral is called a cocommuting 
square if the dual quadrilateral is a commuting square.

\begin{lemma} \label{commcocomm}
Let $N \subset P, Q \subset M$ be a quadrilateral of $II_{1}$ factors, where $N \subset M$ is an 
irreducible finite-index inclusion. Consider the multiplication map of \ref{mulmap}
from $L^{2}(P) \otimes_{N} L^{2}(Q)$ to $L^{2}(PQ)$. The quadrilateral commutes iff 
this map is injective and cocommutes iff the map is surjective.
\end{lemma}

\begin{proof}
The quadrilateral commutes iff $e_{P}e_{Q}=e_{Q}e_{P}$ iff $e_{P}e_{Q}=e_{N}$. By \ref{mulfor} 
 this is equivalent to $\displaystyle \frac{1}{[M:N]}=tr(e_{N})=tr(e_{P}e_{Q})=\displaystyle \frac{tr(e_{P})tr(e_{Q})}{tr(e_{PQ})}$, or $dim_{N}L^{2}(PQ)=[M:N]tr(e_{PQ})=[M:N]^{2}tr(e_{P})tr(e_{Q})=
$\\$dim_{N}L^{2}(P)dim_{N}L^{2}(Q)$. But by \ref{isosub} $L^{2}(PQ)$ is isomorphic to a submodule of $L^{2}(P) \otimes_{N} L^{2}(Q)$, so the two have the same $N$-dimension iff they are in fact isomorphic, which is equivalent to the injectivity of the multiplication map.
 
The quadrilateral cocommutes iff $e_{\bar{P} \bar{Q}}=e_{\bar{Q}\bar{P}}=e_{M}$. By \ref{trfor2} , 
this is equivalent to $dim_{N}L^{2}(PQ)=\frac{1}{tr(e_{\bar{P} 
\bar{Q}})}=\frac{1}{tr(e_{M})}=dim_{N}L^{2}(M)$, which is clearly equivalent to $L^{2}(PQ)=L^{2}(M)$, or 
$e_{PQ}=1$.
\end{proof}

\begin{corollary} \label{dimprod}
The quadrilateral commutes iff $$dim_{N}L^{2}(PQ)=dim_{N}(L^{2}(P) \otimes_{N} L^{2}(Q))=[P:N][Q:N]$$.
\end{corollary}

\begin{corollary}  \label{pqeqqp}
The quadrilateral cocommutes iff $L^{2}(PQ)=L^{2}(QP)$.
\end{corollary}

\begin{proof}
If the quadrilateral cocommutes, then $L^{2}(PQ)=L^{2}(M)=L^{2}(QP)$. 
Conversely, if $L^{2}(PQ)=L^{2}(QP)$, then $e_{PQ}=e_{QP}$. By \ref{landauproj} 
 $e_{PQ}$ is a scalar multiple of $e_{P} \circ e_{Q}$, so $e_{PQ} \circ e_{PQ}$ 
is a scalar multiple of $(e_{P} \circ e_{Q}) \circ (e_{P} \circ e_{Q})=
e_{P} \circ (e_{Q} \circ e_{P}) \circ e_{Q}=e_{P} \circ (e_{P} \circ e_{Q}) 
\circ e_{Q}=(e_{P} \circ e_{P}) \circ (e_{Q} \circ e_{Q})$, which is a scalar multiple 
of $e_{P} \circ e_{Q}$. This implies that $e_{PQ}$ is a biprojection. The corresponding
subfactor has to contain both $P$ and $Q$ so is all of $M$.
So $L^{2}(PQ)=L^{2}(M)$ and the quadrilateral cocommutes.
\end{proof}

In fact one doesn't need the Hilbert space completion for this:

\begin{theorem}
Let $N \subset P, Q \subset M$ be a quadrilateral of $II_{1}$ factors,
 where $N \subset M$ is an irreducible finite-index inclusion. Consider 
the multiplication map from the (algebraic) bimodule tensor product $P \otimes_{N} Q$ to $M$. 
The quadrilateral commutes iff this map is injective and cocommutes iff the map is surjective.
\end{theorem}

\begin{proof}
(a) Injectivity. If the algebraic map from $P \otimes_{N} Q$ to $M$ has a kernel then it is 
obvious that the $L^2$ map  does. On the other hand the kernel $\mathfrak K$
of the $L^2$ map $\mu$ is a closed
$N$-$N$ sub-bimodule of $L^2(M_1)$ (under the isomorphism of  $L^{2}(M) \otimes_{N} L^{2}(M)$
with $L^2(M_1)$, and by the form of elements in the first relative commutant
 the orthogonal projection onto $\mathfrak K$ sends $M_1$ to
itself so there are elements of $M_1$ in $\ker \mu$. Moreover since $M_1\cong M\otimes_N M$
the map $E_P \otimes E_Q$ produces an element of $\ker \mu$ in $Pe_NQ$.

%Consider the inclusion $L^{2}(P) \otimes_{N} L^{2}(Q) \subseteq L^{2}(M) \otimes_{N} L^{2}(M) \cong L^{2}(M_{1})$. The bounded vectors of $L^{2}(M_{1})$ may be identified with $M_{1}=Me_{1}M$, and the bounded vectors of $L^{2}(P) \otimes_{N} L^{2}(Q)$ with $Pe_{1}Q \subseteq M_{1}$. $Pe_{1}Q$ is isomorphic to the algebraic tensor product $P \otimes_{N} Q$ via $pe_{1}q \mapsto p \otimes_{N} q$, so the algebraic multiplication map $P \otimes_{N} Q \rightarrow M $ is just the restriction of the Hilbert space multiplication map of \ref{mulmap}  to the subspace of bounded vectors. Clearly if the algebraic map has a nonzero kernel then so does the Hilbert space map. Conversely, if the Hilbert space map has a nonzero kernel, that kernel is an $N-N$-submodule of $L^{2}(M_{1})$, so it contains bounded vectors and the algebraic map has a nonzero kernel as well. So by \ref{commcocomm}  the quadrilateral commutes iff the algebraic map is injective.   

(b)Surjectivity. The algebraic map is surjective iff $PQ=M$. Clearly $PQ=M$ implies 
$L^{2}(PQ)=L^{2}(M)$. 
Conversely if $L^{2}(PQ)=L^{2}(M)$, $e_{PQ}$ is the identity so $M=PQ$ by
\ref{epqM=PQ}.
%Conversely, by \ref{pqclosed}  , $PQ$ is strongly closed,
% so $PQ$ may be identified with the subspace of bounded vectors of $L^{2}(PQ)$. 
%Therefore $L^{2}(PQ)=L^{2}(M)$ implies that $PQ=M$. So again by \ref{commcocomm} 
% the quadrilateral cocommutes iff the algebraic map is surjective.
\end{proof}

\begin{remark}
Sano and Watatani have already shown that the quadrilateral is a 
cocommuting square iff $PQ=M$ under the additional hypothesis that the quadrilateral is a commuting 
square\cite{SaW}.
\end{remark}

\section{No extra structure}
\subsection{Definition}
Let $N \subset M$ be an inclusion of II$_{1}$ factors with associated tower 
$M_{-1} \subset M_{0} \subset M_{1} \subset ...$, where $M_{-1}=N$, $M_{0}=M$,
 and $M_{k+1}$, $k \geq 0$ is the von Neumann algebra on $L^{2}(M_{k})$ generated by
 $M_{k}$ and $e_{k+1}$, the projection onto $L^{2}(M_{k-1})$.
 Each $e_{k}$ commutes with $N$, so $\{ 1,e_{1},..,e_{k} \}$ generates a *-subalgebra, 
which we will call $TL_{k+1}$, of the $k^{th}$ relative commutant $N' \cap M_{k}$. 

To motivate the following definition (which first occurs in \cite{J31})
consider the case where $N=R^G, M=R^H$ where
$G$ is a finite group of outer automorphisms of the II$_1$ factor $R$.
It is well known that, as a vector space, $N'\cap M_k$ is the set of $G-$invariant
functions on $X^{k+1}$ where $X=G/H$. Thus the {\it transivity} of the action of $G$ on
$X$ is measured by the dimension of $N'\cap M_k$ - an action is $k+1$-transitive
if its dimension is the same as that for the full symmetric group $S_X$. Moreover any
function invariant under $S_X$ is necessarily invariant under $G$ so the relative
commutants for $R^G\subseteq R^H$ always contain a copy of those coming from $S_X$.
The invariants under $S_X$ in this context are sometimes called the partition algebra
so transivity (or rather lack of it) is measured by how much bigger $N'\cap M_k$
is than the partition algebra. Now for a general subfactor $N\subseteq M$ a similar
situation occurs: $N'\cap M_k$ aways contains $TL_{k+1}$. Since this is, for $k>3$, 
strictly smaller in dimension than the partition algebra we see that if we think of
subfactors as "quantum"  spaces $G/H$ they might be "more transitive" than finite
group actions.   

\begin{definition}
Call a finite-index subfactor $N \subseteq M$ 
 $k$-{\rm supertransitive} (for $k>1$) if $N' \cap M_{k-1} = TL_k$. 
We will say $N\subseteq M$ is {\rm supertransitive} if it is $k$-supertransitive
for all $k$.   
\end{definition}

Since $\dim TL_k$ is the same as the partition algebra for $k=1,2,3$ it is
natural to call a $1,2$ or $3$-supertransitive subfactor transitive, $2$-transitive
or $3$-transitive respectively.

\begin{remark}  \label{transi}
$N\subseteq M$ is transitive iff it is irreducible, i.e. $N'\cap M \cong \mathbb C$,
it is $2$-transitive iff the $N-N$ bimodule $L^2(M)$ has two irreducible components
and $3$-transitive iff $\dim N'\cap M_2 \leq 5$.
Supertransitivity of $N\subseteq M$ is the same as saying its principal graph 
is $A_n$ for some $n=2,3,4,..., \infty$.
\end{remark}

\begin{lemma} \label{fusionrules}
Suppose  $N \subset M$ is supertransitive. If $[M:N] \geq 4$ then there is a sequence 
of irreducible $N-N$ bimodules
 $ V_{0}, V_{1}, V_{2}... $ such that $L^{2}(N) \cong V_{0}$,
 $L^{2}(M) \cong V_{0} \oplus V_{1}$, and $V_{i} \otimes V_{j} \cong \oplus_{k=|i-j|}^{i+j}V_{k}$.
 If $[M:N]=4cos^{2}(\frac{\pi}{n})$ then the sequence terminates at $V_{l}$, 
where $l=[\frac{n-2}{2}]$, and the fusion rule
 is: $V_{i} \otimes V_{j} \cong \oplus_{k=|i-j|}^{(\frac{n-2}{2})-|(\frac{n-2}{2})-(i+j)|}V_{k}$. (see \cite{BJ2} ) \\
In either case, we have $\displaystyle \dim_{N}V_{k}=[M:N]^{k}T_{2k+1}(\frac{1}{[M:N]})$, 
where $\{ T_{k}(x) \}$ is the sequence of polynomials defined recursively by
 $T_{0}(x)=0$, $p_{1}(x)=1$, and $T_{k+2}(x)=T_{k+1}(x)-xT_{k}(x)$.
\end{lemma}
\begin{corollary} \label{dimfor}
$dim_{N}V_{1}=[M:N]-1$ and \\$dim_{N}V_{2}=[M:N]^{2}-3[M:N]+1$.
\end{corollary}

\begin{remark} \label{modfusionrules}
If $N \subset M$ is $2k$-supertransitive, then there is a sequence 
of irreducible bimodules $V_{0}, ..., V_{k}$ for which the above fusion rules and 
dimension formula hold as long as $i + j \leq k$. 
\end{remark}

Now let $N\subseteq P,Q \subseteq M$ be a quadrilateral of finite index subfactors. We
will call the four subfactors $N\subseteq P$,$N \subseteq Q$,$P \subseteq M$, and $Q\subseteq M$ 
the {\it elementary subfactors}.

\begin{definition} A quadrilateral as above will be said to have {\rm no extra structure}
if all the elementary subfactors are supertransitive.
\end{definition}

\begin{example}
Let $G=S_{3}$ and let $H$ and $K$ be distinct two-element subgroups of $G$. Given an outer action 
of $G$ on a $II_{1}$ factor $M$, let $N=M^{G}$, and let $P=M^{H}$ and $Q=M^{K}$. 
Then $N \subset P,Q \subset M$ is a quadrilateral which cocommutes
 (since $M' \cap M_{2} \cong l^{\infty}(G)$ is Abelian) but does not commute (since $HK \neq KH$). 

%Ocneanu (?) has shown that the principal graph for any subfactor of index less than four must be one of the $A-D-E$ Dynkin diagrams. In particular, index two subfactors have principal graph $A_{3}$ and index three subfactors have principal graph $A_{5}$. Therefore this quadrilateral has no extra structure, since $[M:P]=[M:Q]=2$ and $[P:N]=[Q:N]=3$. 
This quadrilateral has no extra structure since the permutations actions
 of $S_2$ and $S_3$ are as transitive as
possible. The dual quadrilateral also has no extra structure.
\end{example}

\subsection{Consequences of  supertransitivity.}
Let $N \subset P,Q \subset M$ be a quadrilateral of $II_1$ factors, where $N \subset M$ is 
an irreducible inclusion with finite index. We also have the dual quadrilateral 
$M \subset \bar{P},\bar{Q} \subset M_1$. Let $N \subset P \subset P_{1} ...$ be the tower 
for $N \subset P$ , and similarly for $Q$.

\begin{lemma} \label{pisoq}
If $N\subseteq P$ and $N\subseteq Q$ are 2-transitive  and the quadrilateral does not 
commute then $L^{2}(P) \cong L^{2}(Q)$ as $N-N$-bimodules, and therefore $[P:N]=[Q:N]$.
\end{lemma}

\begin{proof}
By \ref{transi}  write
$L^{2}(P)=L^{2}(N) \oplus V$, where $V$ is an irreducible $N-N$ bimodule. 
Similarly $L^{2}(Q)=L^{2}(N) \oplus W$, for some irreducible $N-N$ bimodule $W$.
Since $e_{P}e_{Q}$ is an $N-N$ intertwiner of $L^{2}(M)$ which fixes $L^{2}(N)$, 
leaves $L^{2}(N)^{\perp}$ invariant and whose range is contained in $L^{2}(P)$, it
 maps $W$ into $V$. Since $W$ is irreducible, $ker(e_{P}e_{Q}|_{W})$ must either be zero 
or all of $W$. The former is impossible since that would imply $e_{P}e_{Q}=e_{N}$, which 
is contrary to our assumption that the quadrilateral does not commute. Thus $V \cong W$, 
and $dim_{N}V=dim_{N}W$.
\end{proof}

\begin{corollary} 
$L^{2}(P) \otimes_{N} L^{2}(Q) \cong L^{2}(P) \otimes_{N} L^{2}(P) \cong L^{2}(P_{1})$.
\end{corollary} 

\begin{lemma}  \label{meqpqp}
If $P \subseteq M$ is $2$-transitive then $L^{2}(PQP)=L^{2}(M)$.
\end{lemma}

\begin{proof}
By \ref{transi} write $L^{2}(M) \cong L^{2}(P) \oplus W$ for some irreducible $P-P$ bimodule $W$. 
Since $L^{2}(PQP)$ is a $P-P$ submodule of $L^{2}(M)$ which is strictly larger than $L^{2}(P)$, it
 must in fact be equal to $L^{2}(M)$. 
\end{proof}

\begin{remark}
Suppose all of the elementary inclusions of the quadrilateral are $2k$-supertransitive 
for some $k \geq 1$. Then the elementary inclusions of the dual quadrilateral are also 
$2k$-supertransitive.  Putting together \ref{pqpmulmap}, \ref{meqpqp}, and \ref{pisoq}, 
we find that as an $N-N$ bimodule, $L^{2}(M)$ is a quotient of $\otimes_{N}^{3} L^{2}(P)$. 
If $k \geq 3$ then the irreducible submodules of $L^{2}(M)$ belong to 
$\{ V_{0},V_{1},V_{2},V_{3} \}$, where the $\{ V_{i} \}$ are as in \ref{modfusionrules} 
 for the $6$-supertransitive inclusion $N \subset P$. Similarly, as an $M-M$-bimodule,
 $L^{2}(M_{1})$ is a quotient of $\otimes_{M}^{3} L^{2}(\bar{P})$. We will write
 $U_{0}$, $U_{1}$ etc. for the irreducible $M-M$ bimodules ocurring in the decomposition of 
the first $k$ tensor powers of $L^{2}(\bar{P})$.  
\end{remark}

For convenience we state the following rewording of a lemma in \cite{PP2}
 which we will be using repeatedly:

\begin{lemma} \label{pplemma}
If the $N-N$-bimodule decomposition of $L^{2}(M)$ contains $k$ copies of the
 $N-N$-bimodule $R$, then $k \leq dim_{N} R$. In particular, $L^{2}(M)$ contains only 
one copy of $L^{2}(N)$.
\end{lemma}

\begin{proof}
$_{N} L^{2}(M)_{N} \cong  (_{N}L^{2}(M)_{M}) \otimes_{M} ({}_{M} L^{2}(M)_{N})$,
 so if $_{N} L^{2}(M)_{N}$ contains $k$ copies of $R$, then by Frobenius reciprocity
 $R \otimes_{N} ({}_{N} L^{2}(M)_{M})$ contains $k$ copies of the $N-M$ bimodule
 $_{N} L^{2}(M)_{M}$, which implies that \\
$dim_{N} (R \otimes_{N} {}_{N} L^{2}(M)_{M})=dim_{N}(R) [M:N] \geq k dim_{N} (_{N} L^{2}(M)_{M})
\\=k[M:N]$. 
\end{proof}

\begin{lemma} \label{pqdecomp}
If $N\subseteq P$ and $N\subseteq Q$ are 4-supertransitive and the quadrilateral does not commute
 then the $N-N$-bimodule $L^{2}(PQ)$ isomorphic to one of the 
following: $V_{0} \oplus 2V_{1} \oplus V_{2}$, $V_{0} \oplus 3V_{1} \oplus V_{2}$, 
or $V_{0} \oplus 3V_{1}$, where the $V_{i}$ are as in \ref{fusionrules}   
(for the 4-supertransitive inclusion $N \subset P$).
\end{lemma}

\begin{proof}
By \ref{isosub}, $L^{2}(PQ)$ is isomorphic to a submodule of $L^{2}(P_{1})$. 
A decomposition of $L^{2}(P_{1})$ into $N-N$-submodules corresponds to a decomposition 
of the identity in $N' \cap P_{3}$. 

If $dim(N' \cap P_{3})=14$ then
 $N' \cap P_{3} \cong M_{2}(\mathbb{C}) \oplus M_{3}(\mathbb{C}) \oplus \mathbb{C}$,
 where the first summand corresponds to $V_{0}$, the second to $V_{1}$, and the third to $V_{2}$.
 So $L^{2}(P_{1}) \cong 2V_{0} \oplus 3V_{1} \oplus V_{2}$. By \ref{pplemma}  , 
$L^{2}(PQ)$ contains only one copy of $L^{2}(N)$. Also, by \ref{pisoq}  , 
$L^{2}(Q) \cong L^{2}(P)$, but $L^{2}(P) \neq L^{2}(Q) $ so $L^{2}(PQ)$ contains
 at least two copies of $V_{1}$. It is impossible that $L^{2}(PQ) \cong V_{0} \oplus 2V_{1}$, 
since that would imply that $L^{2}(PQ)=L^{2}(P+Q)=L^{2}(QP)=L^{2}(M)$, which would imply that 
$[M:P]=\displaystyle \frac{dim_{N}L^{2}(M)}{dim_{N}L^{2}(P)} < 2$. That leaves the three 
possibilities above. If $dim(N' \cap P_{3}) < 14$ then the argument is essentially the same, 
except there is no $V_{2}$, so only one possibility remains.
\end{proof}

\subsection{Cocommuting quadrilaterals with no extra structure.}
{\large NOTATION: from now on the supertransitivity hypotheses will guarantee
that $[M:P]=[M:Q]$. We introduce the following notational conventions:
$$[M:P]=\beta, \quad [P:N]=\alpha, \quad [M:N]=\gamma=1/\tau$$ which we will use without further 
mention.}

\begin{lemma}
If $N \subset P$ and $N \subset Q$ are $2$-transitive, then 
$e_{P}e_{Q}e_{P}=e_{N} + \lambda (e_{P}-e_{N})$, 
where $\lambda=\displaystyle \frac{tr(e_{\bar{P}\bar{Q}})^{-1}-1}{[P:N]-1}$. 
\end{lemma}

\begin{proof}
That $e_{P}e_{Q}e_{P}=e_{N} + \lambda (e_{P}-e_{N})$ for some $\lambda$ follows from the 
fact that $e_{P}e_{Q}e_{P}$ is an $N-N$ intertwiner of $L^{2}(P) \cong V_{0} \oplus V_{1}$ which 
is the identity on $ L^{2}(N)$. To compute $\lambda$
, note that $tr(e_{P}e_{Q}e_{P})=tr(e_{N})+\lambda tr(e_{P}-e_{N})=\displaystyle \frac{1}{\gamma}
 + \lambda \displaystyle \frac{\alpha-1}{\gamma}$. Solving for $\lambda$ and using
 $tr(e_{P}e_{Q}e_{P})=\displaystyle \frac{1}{\gamma tr(e_{\bar{P}\bar{Q}})}$ (by \ref{trfor2}) 
completes the proof.
\end{proof}

\begin{corollary}
$tr((e_{P}e_{Q}e_{P})^{2})=\displaystyle \frac{1+\lambda^{2}([P:N]-1)}{[M:N]}$.
\end{corollary}

\begin{lemma}
If the quadrilateral cocommutes and $e_{\bar{P}\bar{Q}}e_{\bar{Q}\bar{P}}=
e_{\bar{Q}\bar{P}}e_{\bar{P}\bar{Q}}$ then $dim_{M}L^{2}(\bar{P}\bar{Q} + \bar{Q}\bar{P})=$  \\ 
$[M:P]^{2}(2-\displaystyle \frac{[M:P]}{[P:N]}(1+(\displaystyle 
\frac{[P:N]-[M:P]}{[M:N]-[M:P]})^{2}([P:N]-1)))$. 
\end{lemma}

\begin{proof}
Since $e_{\bar{P}\bar{Q}}e_{\bar{Q}\bar{P}}=e_{\bar{Q}\bar{P}}e_{\bar{P}\bar{Q}}$,
 $dim_{M}L^{2}(\bar{P}\bar{Q} + \bar{Q}\bar{P})=
\gamma(2tr(e_{\bar{P}\bar{Q}})-tr(e_{\bar{P}\bar{Q}}e_{\bar{Q}\bar{P}}))$.

Since the quadrilateral cocommutes, $$tr(e_{\bar{P}\bar{Q}})=\displaystyle
 \frac{dim_{M}L^{2}(\bar{P}\bar{Q})}{\gamma}=\frac{dim_{M}L^{2}(\bar{P})dim_{M}L^{2}(\bar{Q})}
{\gamma}=\displaystyle \frac{\alpha^2}{\gamma}=\frac{\beta}{\alpha}.$$
By (the dual version of) \ref{trepqeqp}  , $tr(e_{\bar{P}\bar{Q}}e_{\bar{Q}\bar{P}})=
(\delta tr(e_{\bar{P}\bar{Q}}))^{2} tr((e_{P}e_{Q}e_{P})^{2})=
tr(e_{\bar{P}\bar{Q}})^{2}(1+\lambda^{2}(\alpha-1))=(\displaystyle 
\frac{\beta}{\alpha})^{2}(1+\lambda^{2}(\alpha-1))$.
Also, since $tr(e_{\bar{P}\bar{Q}})=\displaystyle \frac{\beta}{\alpha}$, 
we have $\lambda=\displaystyle \frac{\alpha-\beta}{\gamma-\beta}$.
Putting all this together gives the result.
\end{proof}

\begin{corollary} \label{crazytraceformula}
In the special case that $[M:P]=[P:N]-1$, the formula becomes
 $dim_{M}L^{2}(\bar{P}\bar{Q} + \bar{Q}\bar{P})=[M:P]^{2}+[M:P]-1$
\end{corollary}

\begin{theorem} \label{cocommthm}
If the quadrilateral cocommutes but does not commute, and $N\subseteq P$ and $N\subseteq Q$ 
are 4-supertransitive  then $N$ is the fixed point algebra of an outer $S_{3}$ action on $M$.
\end{theorem}

\begin{proof}
Since the quadrilateral does not commute, $L^{2}(P) \cong L^{2}(Q)$ as $N-N$-bimodules, by
 \ref{pqeqqp}  . Since the quadrilateral cocommutes, $L^{2}(M)=L^{2}(PQ)$, and since
 $N' \cap P_{3}\leq 14$, by \ref{pqdecomp} the isomorphism type of $L^{2}(M)$ is one 
of: $V_{0} \oplus 2V_{1} \oplus V_{2}$, $V_{0} \oplus 3V_{1} \oplus V_{2}$, or $V_{0} \oplus 3V_{1}$.
 For each of these cases we can explicitly compute $\beta$ as a function of $\alpha$ using the 
formula $\beta=\displaystyle \frac{\gamma}{\alpha}=\displaystyle \frac{dim_{N} L^{2}(M)}{\alpha} $ 
and the dimension formulas of \ref{dimfor}  .

Case 1: $L^{2}(M) \cong V_{0} \oplus 3V_{1} \oplus V_{2}$

In this case, $[\bar{P}:M]=\beta=\displaystyle \frac{dim_{N} L^{2}(M)}{\alpha}=\displaystyle 
\frac{1 + 3(\alpha-1) + \alpha^{2}-3\alpha+1}{\alpha}=\alpha-\displaystyle \frac{1}{\alpha}$. Since
 the quadrilateral cocommutes, by \ref{dimprod}  $dim_{M}(L^{2}(\bar{P} \bar{Q}))=
([P:N]-\displaystyle \frac{1}{\alpha})^{2}$. But then the dimension of its orthogonal 
complement (in $L^{2}(M_{1})$) is $dim_{M}L^{2}(M_{1})-dim_{M}L^{2}(\bar{P} \bar {Q})=
[P:N]^{2}-1-(\alpha-\displaystyle \frac{1}{\alpha})^{2}=1-\displaystyle \frac{1}{\alpha^{2}}<1$,
which is impossible by \ref{pplemma}  .

Case 2: $L^{2}(M) \cong V_{0} \oplus 3V_{1}$

In this case, $\beta=\displaystyle \frac{1+3(\alpha-1)}{\alpha}=3-\displaystyle \frac{2}{\alpha}$, 
which necessarily equals $4cos^{2}\frac{\pi}{5}$. 
(The only other admissible index  value less than three is two, but that would 
imply that the total index is four and then the 
quadrilateral would commute.) Then we have the identity $\beta^{2}=3 \beta -1$, and $\alpha=2 \beta$.
Since $L^{2}(M) \cong V_{0} \oplus 3V_{1}$, any intermediate subfactor must have index 
equal to $\displaystyle \frac{1+3(\alpha-1)}{1+k(\alpha-1)}$ for $k=1$ or $k=2$. So to 
eliminate this case it suffices to find a proper subfactor of $M$ with an integer valued index, 
for which it suffices to find an $M-M$-submodule of $L^{2}(M_{1})$ whose dimension over $M$ is $1$.
 
$L^{2}(\bar{P}+\bar{Q})$ has $M$-dimension $2dim_{M}L^{2}(\bar{P})-dim_{M}L^{2}(M)=2 \beta -1$. Its 
 orthogonal complement in $L^{2}(\bar{P}\bar{Q})$, which we shall call $T$, has $M$-dimension 
 $dim_{M}L^{2}(\bar{P}\bar{Q})-dim_M L^{2}(\bar{P}+\bar{Q})=\beta ^{2} -  (2 \beta -1)=\beta$.
 Since $\beta < 3$, if $T$ is reducible, one of its irreducible components must have $M$-dimension 
$1$, and we 
 are finished. Similarly, if $T'$, the orthogonal complement of
 $L^{2}(\bar{P}+\bar{Q})$ in $L^{2}(\bar{Q}\bar{P})$, is reducible then we get a submodule 
of $M$-dimension $1$.
 
 If $T$ and $T'$ are both irreducible, then $L^{2}(\bar{P}\bar{Q}) 
\cap L^{2}(\bar{Q}\bar{P})=L^{2}(\bar{P}+\bar{Q})$. Then if $S$ is
 the orthogonal complement of $L^{2}(\bar{P}\bar{Q}+\bar{Q}\bar{P})$ in $L^{2}(M_{1})$,
 we have  $dim_{M}S=dim_M L^{2}(M_{1}) - (2dim_M L^{2}(\bar{P}\bar{Q})-dim_{M}L^{2}(\bar{P}+\bar{Q}))
=\\2\beta^{2}-(2\beta^{2}-(2\beta-1))=2\beta-1$. Since $dim(M' \cap M_{2})=dim(N' \cap 
 M_{1})=10$, $S$ must break into 
 $3$ components, one of which must have $M$-dimension $1$.

Case 3: $L^{2}(M) \cong V_{0} \oplus 2V_{1} \oplus V_{2}$

In this case $\beta=\displaystyle \frac{1+2(\alpha-1)+\alpha^{2}-3\alpha+1}{\alpha}=\alpha-1$. 
Note that $dim(N' \cap M_{1})=6$, and therefore also $dim(M' \cap M_{2})=6$. Because
 $L^{2}(M) \subset L^{2}(\bar{P}) \subset L^{2}(\bar{P} + \bar{Q}) \subset L^{2}(\bar{P} \bar{Q})
 \subset L^{2}(M_{1}) $ is a strictly increasing chain of $M-M$ bimodules ($\bar{P} \bar{Q}$ cannot
 be all of $M_{1}$ because the quadrilateral does not commute), $M' \cap M_{2}$ must be Abelian. 
If we let $x=\beta$ (so that $\alpha=x+1$), then $\gamma=x^{2}+x$, and by \ref{crazytraceformula}  
 we have that
$dim_{M} L^{2}(\bar{P}\bar{Q} + \bar{Q}\bar{P})=x^{2}+x-1$, and so the dimension of its orthogonal
 complement in $L^{2}(M_{1})$ is $1$.

It is then easy to see that the dimensions of the six distinct irreducible submodules
 of $L^{2}(M_{1})$ are $1, x-1, x-1, x^{2}-2x-1, x^{2}-2x-1, 1$. But then summing we
 find that $2x^{2}-2x-2=dim_{M}L^{2}(M_{1})=x^{2}+x$, which implies that $x=2$. 
So $[\bar{P}:M]=[\bar{Q}:M]=2$, and $[M_{1}:M]=6$.

So by Goldman's theorem, \cite{G} $M_{1}$ is the crossed product of $M$ by $S_{3}$, or, equivalently,
$N$ is the fixed point subalgebra of an outer $S_{3}$ action on $M$. 

\end{proof}

\section{Restrictions on the principal graph}
If the quadrilateral has no extra structure then we obtain severe restrictions on the principal graph. Specifically, for a noncommuting, noncocommuting quadrilateral with no extra structure the principal graph is completely determined.

\subsection{Structural restrictions}

\begin{lemma} \label{2summands}
If the quadrilateral neither commutes nor cocommutes, and all the elementary subfactors are $6$-supertransitive, then $N' \cap M_1$ and $M' \cap M_2$ both have more than two simple summands.
\end{lemma}

\begin{proof}
First suppose that $N' \cap M_1$ and $M' \cap M_2$ both have exactly two simple summands. 
Then  $L^{2}(M)=V_{0} \oplus kV_{1}$ for some integer $k$.  So we have $\beta=\displaystyle \frac{\gamma}{\alpha}=\displaystyle \frac{dim_{N}(V_{0} \oplus kV_{1})}{\alpha}=\displaystyle \frac{1+k(\alpha-1)}{\alpha}=k-\displaystyle \frac{k-1}{\alpha}<k$. By \ref{pplemma}, $k \leq dim_{N} V = \alpha-1<\alpha$, and so $[M:P]<\alpha$. But we can perform the same calculation in the dual quadrilateral to find that $\alpha=[M_{1}:\bar{P}]<[\bar{P}:M]=[M:P]$, which is a contradiction.

Now suppose that only $M' \cap M_{2}$ has exactly two simple summands, and write $L^{2}(M_{1}) \cong U_{0} \oplus lU_{1}$. Note that because of the 6-supertransitivity hypothesis, the first few tensor powers of $U_{1}$ decompose according to the fusion rules of \ref{fusionrules}. By \ref{pqdecomp}, $L^{2}(\bar{P}\bar{Q}) \cong U_{0} \oplus 3U_{1}$, and since the quadrilateral does not commute, by \ref{pqeqqp} $L^{2}(\bar{P}\bar{Q}) \neq L^{2} (\bar{Q} \bar{P})$, so $l$ must be at least $4$. By \ref{meqpqp} and \ref{pqpmulmap}, $L^{2}(M_{1})$ is a quotient of $L^{2}(\bar{P}\bar{Q}) \otimes_{M} L^{2}(\bar{P}) \cong (U_{0} \oplus 3U_{1})  \otimes_{M} (U_{0} \oplus U_{1}) \cong 4U_{0} \oplus 7U_{1} \oplus 3U_{2}$, where the last isomorphism comes from the fusion rule $U_{1} \otimes_{M} U_{1} \cong U_{0} \oplus U_{1} \oplus U_{2}$  (If $\alpha<3$ then $U_{2}=0$). So we find that $4\leq l \leq 7$.

Similarly, $L^{2}(M)$ is a quotient of $L^{2}(PQ) \otimes_{N} L^{2} (P)$, which in all cases of \ref{pqdecomp} is a quotient of $(V_{0} \oplus 3V_{1} \oplus V_{2}) \otimes_{N} (V_{0} \oplus V_{1}) \cong V_{0} \oplus 8V_{1} \oplus 5V_{2} \oplus V_{3}$. Thus we may write $L^{2}(M) \cong V_{0} \oplus aV_{1} \oplus bV_{2} \oplus cV_{3}$, where $a$, $b$, and $c$ are integers such that $2 \leq a \leq 8$, $0 \leq b \leq 5$, and $0 \leq c \leq 1$, and $b$ and $c$ are not both $0$. 

But because we have $dim(N' \cap M_{1})=dim(M' \cap M_{2})$, we necessarily have 
$a^{2}+b^{2}+c^{2}=l^{2}$. A quick examination reveals that the only possibility is
 that $l=5$, $c=0$, and $\{a,b\}=\{3,4\}$. But if $l=5$ then $\alpha=\displaystyle 
\frac{\gamma}{\beta}=\displaystyle \frac{[M_{1}:M]}{[\bar{P}:M]}=5-\displaystyle 
\frac{4}{\beta} < 5$, which implies that $a \leq dim_{N} V_{1} < 4$ (by \ref{pplemma}), 
so we may assume that $a=3$ and $b=4$.
Then $\beta=\displaystyle \frac{dim_N V_0 \oplus 3V_1 \oplus 4V_2}{\alpha}=\displaystyle
 \frac{1+3(\alpha-1)+4(\alpha^{2}-3\alpha+1)}{\alpha}=4\alpha^{2}-9\alpha+2$, and since 
$\alpha \geq 3$, we must have $\beta \geq 4$, and then also $\alpha =5- \displaystyle 
\frac{4}{\beta} \geq 4$, so the generic fusion rules of \ref{fusionrules} apply.

Then as an $N-N$ bimodule, 
$L^{2}(M_{1}) \cong L^{2}(M) \otimes_{N} L^{2}(M) \cong \\(V_{0} \oplus 3V_{1} 
\oplus 4V_{2}) \otimes_{N} (V_{0} \oplus 3V_{1} \oplus 4V_{2}) \cong 10V_{0} 
\oplus 39V_{1} \oplus 41V_{2} \oplus 12V_{3} \oplus 16(V_{2} \otimes_{N} V_{2})$,
 where the last isomorphism comes from the fusion rules $V_{1} \otimes_{N} V_{1} 
\cong V_{0} \oplus V_{1} \oplus V_{2} $ and $V_{1} \otimes_N V_{2} \cong V_{1} 
\oplus V_{2} \oplus V_{3}$. Since the $N-N$ 
winer space of $L^{2}(M_{1})$ is $N' \cap M_{3}$, this implies that $dim(N' \cap M_{3}) \geq 10^{2} + 39^{2} + 41^{2} +12^{2}=3446$.

On the other hand, as an $M-M$ bimodule,\\ $L^{2}(M_{2}) \cong L^{2}(M_{1}) \otimes _{M} 
L^{2}(M_{1}) \cong L^{2}(M) \oplus 5U_{1} \otimes L^{2}(M) \oplus 5U_{1} \cong \\26U_{0} 
\oplus 35U_{1} \oplus 25U_{2}$, so $dim(M' \cap M_{4})=26^{2}+35^{2} + 25^{2}=2526$. 
But this contradicts the fact that $dim(N' \cap M_{3})=dim(M' \cap M_{4})$. 
\end{proof}

\begin{lemma}
If the quadrilateral neither commutes nor cocommutes and all the elementary subfactors are 6-supertransitive then $[N:P]$ and  $[M:P]$ are both less than four.
\end{lemma}

\begin{proof}
Suppose on the contrary that the hypotheses are satisfied and that $\alpha \geq 4$. (There is no loss of generality here since if only $\beta \geq 4$ we may consider the dual quadrilateral instead.) Then by \ref{2summands}   $N' \cap M_{1}$ has at least three simple summands. Because the quadrilateral is not cocommuting, by \ref{pqeqqp}   $L^{2}(PQ) \neq L^{2}(QP)$, but they must have the same dimension since by \ref{treqpeqltrepq}   $tr(e_{PQ})=tr(e_{QP})$. We consider three cases, corresponding to the three cases of \ref{pqdecomp}  :

Case 1: $L^{2}(PQ) \cong V_{0} \oplus 3V_{1}$. Then also $L^{2}(QP) \cong V_{0} \oplus 3V_{1}$. Note that these two bimodules intersect in $L^{2}(P+Q) \cong V_{0} \oplus 2V_{1}$, so $L^{2}(PQ + QP) \cong V_{0} \oplus 4V_{1}$. Since $N' \cap M_{1}$ has a third summand, $L^{2}(M)$ must also contain an irreducible submodule whose dimension is at least as great as that of $V_{2}$, by \ref{fusionrules} , so we find that $\gamma=dim_{N}L^{2}(M) \geq dim_{N} V_{0} \oplus 4V_{1} \oplus V_{2} \geq 1 + 4(\alpha-1) + (\alpha^{2}-3\alpha+1)=\alpha^{2}+\alpha-2$, and so $\beta=\displaystyle \frac{\gamma}{\alpha}=\alpha+1-\frac{2}{\alpha}>\alpha$.

Case 2: $L^{2}(PQ) \cong V_{0} \oplus 2V_{1} \oplus V_{2}$. Then $L^{2}(PQ + QP) \cong V_{0} \oplus 2V_{1} \oplus 2V_{2}$. So $\gamma
 \geq dim_{N} V_{0} \oplus 2V_{1} \oplus 2V_{2} =2\alpha^{2}-4\alpha+1$, and again we find that $\beta=2\alpha-4+\displaystyle \frac{1}{\alpha}>\alpha $. (Because $\alpha \geq 4$).

Case 3: $L^{2}(PQ) \cong V_{0} \oplus 3V_{1} \oplus V_{2}$. Then $L^{2}(PQ + QP)$ contains either at least four copies of $V_{1}$ or at least two copies of $V_{2}$ and again we find that $\beta>\alpha$.

But since $\beta > \alpha \geq 4$, we can perform these same calculations in the dual quadrilateral to deduce that $\alpha > \beta $, which is absurd.
\end{proof}

\begin{lemma} \label{indiceseq}
If the quadrilateral neither commutes nor cocommutes and all the elementary 
subfactors are 6-supertransitive, then $[P:N]=[M:P]$.
\end{lemma}

\begin{proof}
By the previous lemma we may assume that $\alpha$ and $\beta$ are both less than four. Because $\alpha<4$, $dim_{N} V_{1} < 3$, so by \ref{pplemma}  $L^{2}(M)$ contains at most, and therefore exactly, two copies of $V_{1}$, and so $L^{2}(PQ) \cong V_{0} \oplus 2V_{1} \oplus V_{2}$. Now $L^{2}(M)$ is a quotient of $L^{2}(PQ) \otimes_{N} L^{2}(P) \cong (V_{0} \oplus 2V_{1} \oplus V_{2}) \otimes_{N} (V_{0} \oplus V_{1}) \cong 3V_{0} \oplus 6V_{1} \oplus 4V_{2} \oplus V_{3}$, so it contains at most four copies of $V_{2}$ and at most one copy of $V_{3}$ (and nothing higher). Also, since $L^{2}(QP)$ is isomorphic, but not equal, to $L^{2}(PQ)$, $L^{2}(M)$ contains at least two copies of $V_{2}$.

So we may write $L^{2}(M) \cong V_{0} \oplus 2V_{1} \oplus bV_{2} \oplus cV_{3}$, with $2 \leq b \leq 4$ and $0 \leq c \leq 1$. Similarly, we may write $L^{2}(M_{1}) \cong U_{0} \oplus 2U_{1} \oplus b'U_{2} \oplus c'U_{3}$, with $2 \leq b' \leq 4$ and $0 \leq c' \leq 1$. Since $1^{2} + 2^{2} + b^{2} + c^{2} =dim(N' \cap M_{1})=dim(M' \cap M_{2})=1^{1} + 2^{2} + b'^{2} + c'^{2}$ and $c$ and $c'$ are each either $0$ or $1$, we must have $b=b'$ and $c=c'$.

Define the function $$f_{b,c}(x)=[1 + 2(x-1) + b(x^{2}-3x+1)
 + c(x^{3}-5x^{2}+6x-1)]/x$$ $$=cx^{2}+(b-5c)x+(2-3b+6c)+\displaystyle \frac{(b-c-1)}{x}.$$ Then $f_{b,c}(\alpha)=\beta$ and $f_{b,c}(\beta)=\alpha$. Define $g_{b,c}(x)=f_{b,c}(x)-x$. Then $g_{b,c}'(x)$ is either $b - 1 - \displaystyle \frac{b-1}{x^{2}}$, or $2x+b-6-\displaystyle \frac{b-2}{x^{2}}$, depending upon whether $c$ is $0$ or $1$. In either case, $g'(x)$ is positive when $x \geq 2$ and so $g(x)$ is then an increasing function.

Now if $\alpha > \beta$, then $g_{b,c}(\beta)=f_{b,c}(\beta)-\beta=\alpha-\beta > 0$, and since $\alpha > \beta$ and $g_{b,c}(x)$ is increasing, $g_{b,c}(\alpha)>0$ as well, so we also have $\beta > \alpha$, which is a contradiction. Similarly we find that $\beta>\alpha$ is impossible. Therfore we must have $\beta=\alpha$.
\end{proof}

\subsection{The principal graph}

\begin{lemma} \label{a11}
There does not exist a noncommuting quadrilateral of subfactors with 
$L^{2}(M) \cong V_{0} \oplus 2V_{1} \oplus 2V_{2}$ and with the principal graph of the elementary subfactors equal to $A_{11}$.
\end{lemma}

\begin{proof}
Suppose such a quadrilateral exists. \\
Then $L^{2} (M_{1}) \cong L^{2}(M) \otimes_{N} L^{2}(M) \cong 9V_{0} \oplus 20V_{1}
 \oplus 20V_{2} \oplus 12V_{3} \oplus 4V_{4}$, and $L^{2} (M_{2}) \cong L^{2}(M_{1}) 
\otimes_{N} L^{2}(M_1) \cong 89V_{0} \oplus 222V_{1} \oplus 254V_{2} \oplus 196V_{3} \oplus 108V_{4} \oplus 32V_{5}$, by the $A_{11}$ fusion rules. (\ref{fusionrules} with $n=12$ gives $V_i \otimes_{N} V_j =\oplus_{|i-j|}^{5-|5-(i+j)|}V_k$.) 

Recalling the principle that each level of the Bratteli diagram for the tower of relative commutants is obtained by reflecting the previous level and adding some "new stuff", with the rule that the "new stuff" connects only to the "old new stuff" 
(see \cite{GHJ}), it is easy to deduce that the Bratteli diagram must include the graph in \ref{dia11brat} .

\begin{diagram} \label{dia11brat}
\vpic{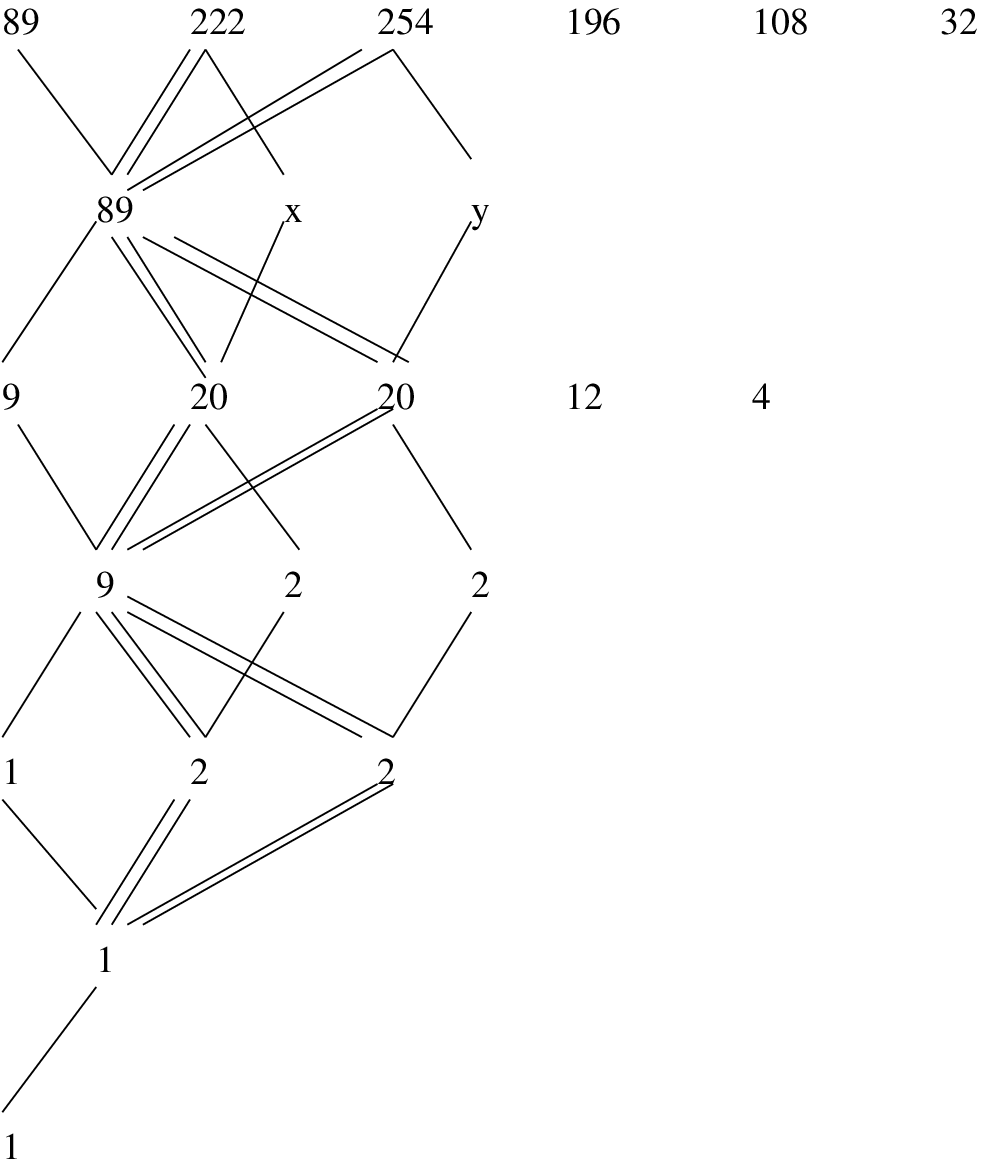} {4 in}    
\end{diagram}

Let $m$ and $n$ be the number of bonds which connect the two "2"s in the fourth row with "12" in the fifth row, respectively. Then we must have $2m+2n=12$, or $m+n=6$. By the reflection principle, there must also be $m$ and $n$ bonds connecting "12" with "x" and "y" respectively, as well as "x" and "y" with "196". This implies that $x \geq 20 + 12m $, $y \geq 20 + 12n$, and $196 \geq m (20+12m) + n (20 +12n)=20(m+n)+12(m^2 + n^2)$ which is absurd since $m+n=6$.

\end{proof}

\begin{lemma}
If the quadrilateral neither commutes nor cocommutes, and the elementary inclusions 
are 6-supertransitive, then $[P:N]=[M:P]=2+ \sqrt{2}$ and 
$L^{2}(M) \cong V_{0} \oplus 2V_{1} \oplus 2V_{2} \oplus V_3 $. 
\end{lemma} 
\begin{proof}
As in the proof of \ref{indiceseq}  , there are six possible isomorphism types 
for $L^{2}(M) \cong V_{0} \oplus 2V_{1} \oplus bV_{2} \oplus cV_{3}$, 
corresponding to $b=2,3,4$ and $c=0,1$. We will eliminate them all except
$b=2,c=1$.

Let $x=\alpha$. From the proof, and the conclusion, of \ref{indiceseq} we have
$$cx^{3}+(b-5c-1)x^2+(2-3b+6c)x+(b-c-1)=0.$$ Let us consider the cases one at a
time:\\
\underline{$c=0,b=2$}\\
Then $x=2+\sqrt 3$ and the only principal graphs possible for $N\subseteq P$ are $A_{11}$ and
$E_6$. But $E_6$ is not $4$-supertansitive and $A_{11}$ was eliminated in \ref{a11}.\\
\underline{$c=0,b=3$}\\
Then $2x^2-7x+2=0$, neither of whose roots is an allowed index value.\\
\underline{$c=0,b=4$}\\
Then $3x^2-10x+3=0$ so $\alpha=3$ which implies $\dim_N(V_2)=1$ which is impossible by
\ref{pplemma}.\\  
\underline{$c=1,b=3$}\\
Then $x^3-3x-x+1=0$ or $x(x^2-3x+1)=2x-1$ which implies $\dim_N(V_2)<2$. Again by 
\ref{pplemma} this is impossible.\\
\underline{$c=1,b=4$}\\
Then $x^3-2x^2-4x+2=0$. The largest root of this equation is between $3$ and $4\cos^2\pi/7$
so it is not a possible index value.

Finally, in the case $c=1,b=2$, $x(x^2-4x+2) =0$ so $\alpha=2+\sqrt 2$ (which is  $4\cos ^2\pi/8)$). 
\end{proof}
\begin{corollary}\label{anglePQ}
With the hypotheses of the previous lemma, $tr(e_{PQ}=\frac{1}{\sqrt 2}$,$tr(e_Pe_Q)=\frac{1}{4+3\sqrt 2}$
 and the angle between $P$ and $Q$ is $\cos^{-1}(\sqrt 2 -1)$.
\end{corollary}
\begin{proof}
By 2-transitivity we know that $e_Pe_Qe_P=e_N + t(e_P-e_N)$ for some number $t$ which
is the square of the cosine of the angle. Moreover by \ref{pqdecomp} we know that 
$\dim_N(L^2(PQ))=1+3(1+\sqrt 2)$. Taking the trace, using \ref{mulfor} and solving
for $t$ we are done.
\end{proof}
\begin{theorem}\label{mainchance}
Let $N\subset P, Q\subset M$ be a noncommuting noncocommuting quadrilateral
with all elementary inclusions 6-transitive. Then $[M:P]=[P:N]=[Q:N]=2+\sqrt 2$ 
and the principal and dual principal graphs for $N\subset M$ are both 
\vpic{princgraph} {0.8 in}  .
\end{theorem}
\begin{proof}
Reduction to this one case is a consequence of the previous results. We need
only compute the principal graph.
Since there is no subfactor with principal graph $D_5$, all the elementary 
subfactors must have principal graph $A_7$. Thus there are only the 4 possible isomorphism
types $V_0,V_1,V_2$ and $V_3$ for the $N-N$ bimodules in $L^2(M),  L^2(M_1),...$
, i.e. the Bratteli diagram for the tower of relative commutants $N'\cap M_K$ has 
at most 4 simple summands for $k$ odd. Since there are 4 simple summands in 
$N'\cap M_1=End_{N-N}L^2(M)$, the subfactor $N\subset M$ is of depth 3.
Moreover if we let $V_{a}=V_{0} \oplus V_{3}$ and 
$V_{b}=V_{1} \oplus V_{2}$, then $L^{2}(M) \cong V_{a} \oplus 2V_{b}$, and
 the fusion rules are very simple: $V_{a} \otimes V_{a}=2V_{a}$, $V_{a} \otimes V_{b}=2V_{b}$,
 and $V_{b} \otimes V_{b}=2V_{A} \oplus 4 V_{b}$. 
So $L^{2}(M_{1}) \cong L^{2}(M) \otimes L^{2}(M) \cong 10 V_{a} \oplus 24 V_{b} \cong 10 V_{0} 
\oplus 24 V_{1} \oplus 24V_{2} \oplus 10 V_{3}$, and there is only one way to fill 
in the $N'\cap M_2$ level of the Bratteli diagram for the tower of relative commutants which 
will thus begin as in Fig \ref{brattprinc}. 

\begin{diagram}\label{brattprinc}
\vpic{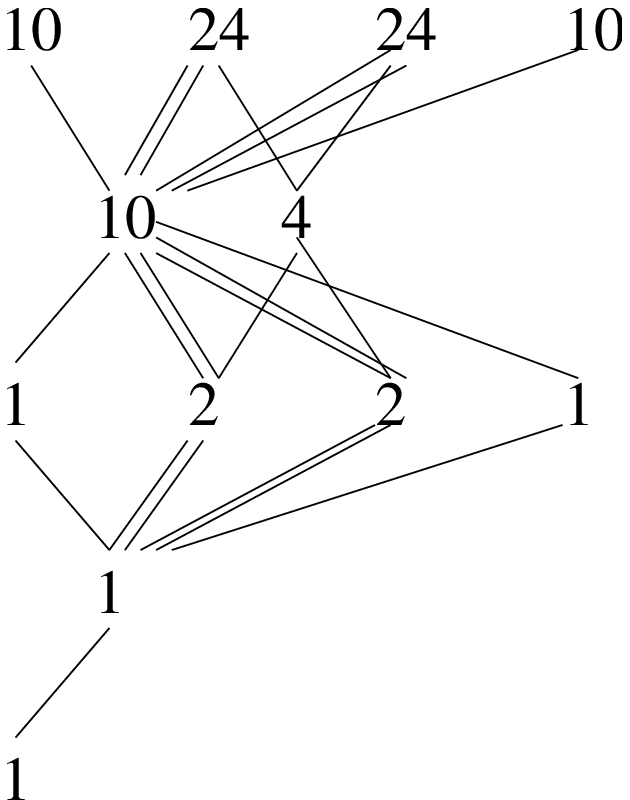} {1.6 in}  
\end{diagram}
By depth 3 we are done. 

The dual principal graph has to be the same as the principal graph since
$M\subset M_1$ satisfies the same hypotheses as $N\subset M$.
\end{proof}

\section {The $6+4\sqrt 2$ example.} \label{the example}

\subsection{Material from ``Coxeter graphs and towers of algebras'.}
We give a general construction for pairs of intermediate
subfactors which seems to be of some interest. Recall two
constructions of subfactors from \cite{GHJ}:

 Let $\Gamma$ be a Coxeter-Dynkin diagram of type
type $A, D$ or $E$ with Coxeter number $k$, with $\Gamma =\Gamma_0 \sqcup
\Gamma_1$ a particular bipartite structure. Construct a pair
$A_0\subset A_1$ of finite dimensional C$^*$-algebras the underlying graph of whose
Bratelli
diagram is $\Gamma$. Thus the minimal central projections in $A_i$ are 
indexed by $\Gamma_i$ for $i=0,1$.
Using the Markov trace $tr$ on $A_1$ iterate the basic construction
to obtain the tower $A_{i+1}=\langle A_i, e_i \rangle$, $e_i$ being the 
orthogonal projection onto $A_{i-1}$. There is a unitary braid group representation
inside the tower obtained by sending the usual generators $\sigma_i$ of
the braid group (see \cite{J4}) to the elements $g_i=(t+1)e_i-1$ with
$t=e^{2\pi i /k}$.\\

First construction-commuting squares.\\

If we attempt to obtain a commuting square from the tower by conjugating
$A_1$ inside $A_2$ by a linear combination of $e_1$ and $1$, we find that
there are precisely two choices up to scalars: $g_1$ and $g_1^{-1}$.
Then the following is a commuting square:

$$\begin{array}{ccc}
B_1=g_1A_1g_1^* &\subset &A_2 \cr
\cup& &\cup \cr
B_0=A_0&\subset &A_1

\end{array}$$

We may then define $B_i$ to be the C$^*$-algebra generated by $B_{i-1}$
and $e_i$ to obtain II$_1$ factors $B_\infty \subseteq A_\infty$ with
index $4\cos^2 \pi/k$. This construction is known to give all subfactors
of index less than 4 of the hyperfinite II$_1$ factor. The Dynkin
diagram $\Gamma$ is the principal graph of the subfactor in the cases
$A_n$, $D_{2n}$, $E_6$ and $E_8$ but not otherwise. For $D_{2n+1}$ 
the principal graph is $A_{4n-1}$.  See \cite{EK7}.
\vskip 10pt

Second construction-GHJ subfactors.\\

The $e_i$'s in the II$_1$ factor $A_\infty$ above generate a II$_1$ factor $TL$
and by a lemma of Skau (see \cite{GHJ})  $TL'\cap A_\infty=A_0$. Thus one may obtain
irreducible subfactors $N\subseteq M$ by choosing a minimal projection 
$p$ in $A_0$, i.e. a vertex of $\Gamma$ in $\Gamma_0$, and setting
$N=pTL$ and $M=pA_\infty p$. These subfactors are known as "GHJ" subfactors
as they first appeared in \cite{GHJ}. We will call the subfactor $TL\subseteq A_\infty$
the "full GHJ subfactor". The indices of the GHJ subfactors are all
finite and were calculated 
in \cite{GHJ} (but note the error there: for $D_n$ using the two univalent vertices connected
to the trivalent one-it should be divided by 2).

\begin{remark}
The cut-down Temperley-Lieb projections $pe_1, pe_2, ...$ satisfy the same relations in the cut-down algebra $pA_{\infty}p$ that the projections $e_1, e_2,...$ do in $A_\infty$. Therefore when discussing $pA_{\infty}p$ we will denote the cut-down Temperley-Lieb projections simply by $e_i$.
\end{remark}

 Using Skau's lemma
Okamoto in \cite{Ok} calculated the principal graphs for the GHJ subfactors as
follows: if $TL_n$ is the C$^*$ algebra generated by $e_1, e_2,...,e_{n-1}$ then
the inclusions:
$$\begin{array}{ccc}
pTL_{n+1} &\subset &pA_{n+1}p \cr
\cup& &\cup \cr
pTL_n&\subset &pA_np

\end{array}$$
are commuting squares for which the Bratteli diagram of the unital
inclusion $pTL_n \subseteq pA_{n}p$ 
may be calculated 
explicitly inductively using one simple rule which follows from the basic
construction.\\

Rule:  If $q$ is a minimal projection in $pTL_n$ and $r$ is a minimal 
projection in $pA_np$ then $e_{n+1}q$ and $e_{n+1}r$
are minimal projections in $pTL_{n+2}$ and $pA_{n+2}p$ respectively, and
the number of edges connecting $q$ to $r$
is equal to the number connecting $e_{n+1}q$ to $e_{n+1}r$. \\

Thus one obtains two Bratteli diagrams depending on the parity of $n$.
For sufficiently large $n$ the inculsion matrices for these Bratteli
diagrams do not change and  the
principal graph for the GHJ subfactor is the underlying bipartite graph of
the stable Bratteli diagram for the 
inclusion $pTL_n \subseteq pA_{n}p$, with distinguished vertex $*$ being the 
$*$ vertex in the Temperley-Lieb type $A$ graph. This specifies
the parity of $n$ needed. Note that the dual principal graph is \underline{not}
in general the inclusion graph with the other parity!\\

\begin{example}\label{d5ghj1}{\rm We take $\Gamma$ to be the Coxeter graph $D_5$ with the minimal projection
$p$ being that corresponding to the trivalent vertex. The two vertical Bratteli diagrams
are those for $pA_\infty p$ and $pTL$, and the inclusions $pTL_n \subset pA_{n}p$ are given
by approximately horizontal dashed lines, except the one which is the GHJ subfactor principal
graph which is made up of the heavy lines at the top of the figure.
We have suppressed the dashed lines for $pTL_5\subset pA_{5}p$ to avoid confusion and because
this inclusion graph is not the principal graph. The figure has been constructed from
the bottom up using the basic construction and the above Rule.}

\vpic{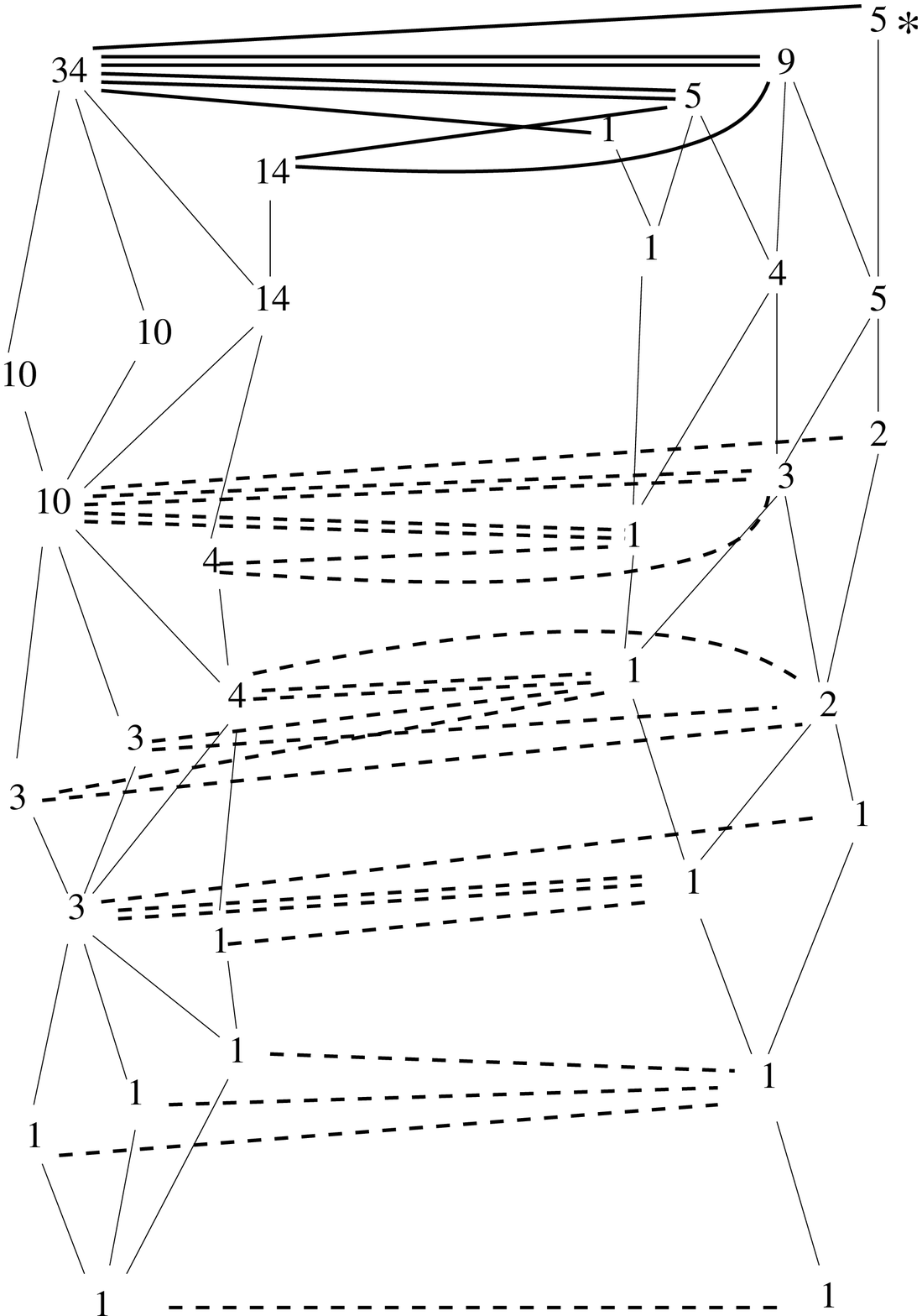} {4in}
\vskip 5pt
{\rm Making the principal graph more visible we obtain:}
\vskip 5pt
\vpic{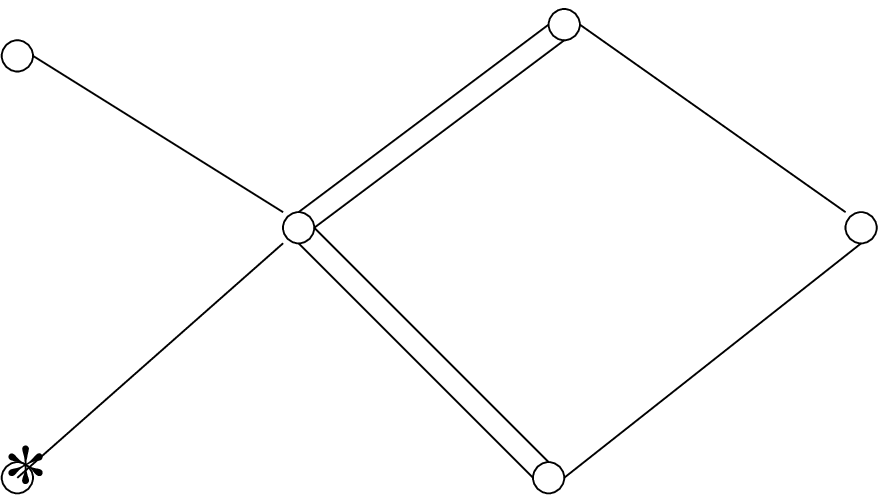} {2in}
\end{example}
\subsection{GHJ Subfactor Pairs.}
Looking again at the commuting square construction from the original Coxeter-Dynkin
diagram we see that we may in fact construct \underline{two} subfactors
of $A_\infty$ by conjugating initially by $g$ and $g^{-1}$!
This construction works in great generality and gives a pair of subfactors 
whenever a subfactor is constructed using the endomorphism method of \cite{GHJ},\cite{J30}.
In fact there is a way to obtain the quadrilateral with no extra structure
by a simpler method, with simpler angle calculation and using only the
real numbers. It seems to be a bit less general than the method using the braid group
so we present it second.

\begin{definition}The \underline{full GHJ subfactor pair} is the pair $\mathfrak P$ and
 $\mathfrak Q$ of subfactors
of the (hyperfinite) II$_1$ factor $A_\infty$ defined as the von Neumann algebras
generated by the $P_n$ and $Q_n$ in the following towers:
\begin{diagram}\label{twotowers}

$$\begin{array}{ccccc}
\cup   &        &   \cup &          & \cup    \cr
P_{n+1}&\subset &A_{n+1} &\supset & Q_{n+1} \cr
\cup   &        & \cup   &          & \cup    \cr
P_n    &\subset & A_n    &\supset & Q_n    \cr
\cup   &        &  \cup  &          &    \cup \cr
 \end{array}$$
 \end{diagram}
Where $A_n$ is as above, $P_1=Q_1=A_0$, $P_2=g_1 A_1 g_1^*$, $Q_2=g_1^* A_1 g_1$ and
$P_{n+1}=\{P_n,e_{n}\}'', Q_{n+1}=\{Q_n,e_{n}\}''$.
\end{definition}

Note that in \ref{twotowers}, all squares involving just  $A$'s and $P$'s or just
$A$'s and $Q$'s are commuting but squares involving $P$'s and $Q$'s may not be.

\begin{definition} Let $TL2$ be the subfactor of $A_\infty$ generated by
all the $e_i$ with $i\geq 2$.
\end{definition}
 
\begin{proposition} $[A_\infty:{\mathfrak P}\cap {\mathfrak Q}]<\infty$.
\end{proposition}
\begin{proof} By construction $e_i \in {\mathfrak P}\cap {\mathfrak Q}$ for all $i\geq 2$. Moreover
$TL2$ is of index $4\cos^2 \pi/k$ in the full GHJ subfactor $TL$ which is in turn of finite index in
$A$ by \cite{GHJ}.
\end{proof} 

Note that $A_0$ is in $TL_2'\cap A_\infty$ and $A_0\subseteq {\mathfrak P}\cap {\mathfrak Q}$.
 We suspect that
${\mathfrak P}\cap {\mathfrak Q}$
 is the von Neumann algebra $TL2\otimes A_0$ generated by $TL2$ and $A_0$. 
We hope to answer this question in a future systematic study of the GHJ subfactor pairs.

Our interest in this paper has been in pairs of subfactors $P,Q\subseteq M$ with
$(P\cap Q)'\cap M =\mathbb C id$. 

\begin{definition} Let $p$ be a projection in $A_0$ that is minimal in $A_1$.
Then the \underline{GHJ subfactor pair} corresponding to $p$ is the pair of
subfactors $$P=p{\mathfrak P}p,Q= p{\mathfrak Q}p \quad\subseteq \quad M=pA_\infty p$$.
\end{definition}

\begin{proposition} If $P,Q\subseteq M$ is a GHJ subfactor pair then
$(P\cap Q)'\cap M =\mathbb C id$.
\end{proposition}
\begin{proof}
By Skau's lemma we know  that the commutant of $TL2$ in $M$ is $A_1$. 
\end{proof} 

A projection in $A_0$ which is minimal in $A_1$  is the same thing as
a univalent vertex in $\Gamma_0$. Note that the subfactor $TL2\subseteq A_\infty$
is then the full GHJ subfactor for the other bipartite structure on $\Gamma$, 
and the subfactor $pTL2 \subseteq p A_\infty p$ is the GHJ subfactor obtained
by choosing the unique neighbour of the original univalent vertex. (This is because
the inclusion $A_1 \subseteq A_2$  can be used as the initial inclusion to 
construct the full GHJ subfactor for the other bipartite structure and $p$ 
is a minimal projection in $A_1$ since we started with a univalent vertex.)

There are not too many choices of univalent vertex, especially up to symmetry. We enumerate them below,
the chosen univalent vertex being indicated with a * :\\

$A_n$ \vpic{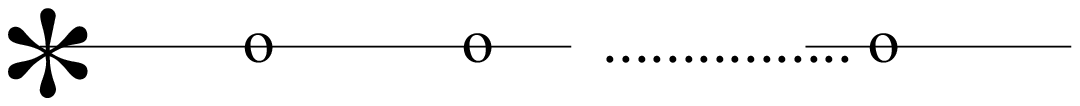} {1.0in} \\

$D_{n,1}$ \vpic{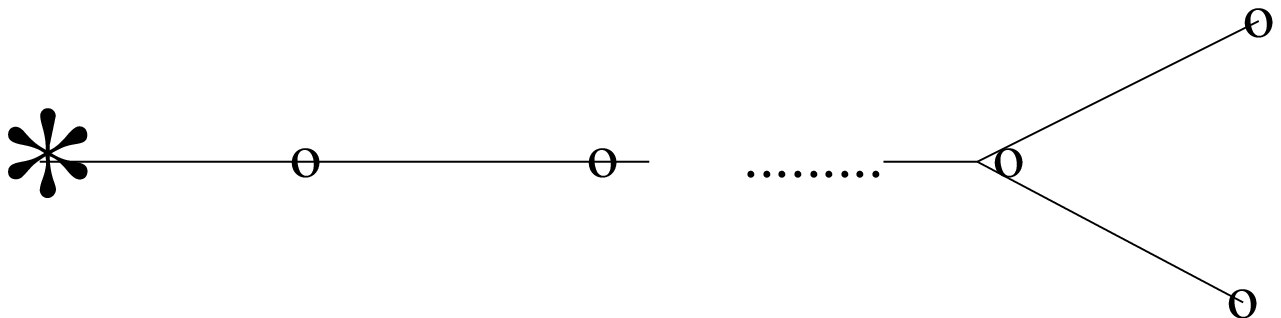} {1.0in} \quad $D_{n,2}$ \vpic{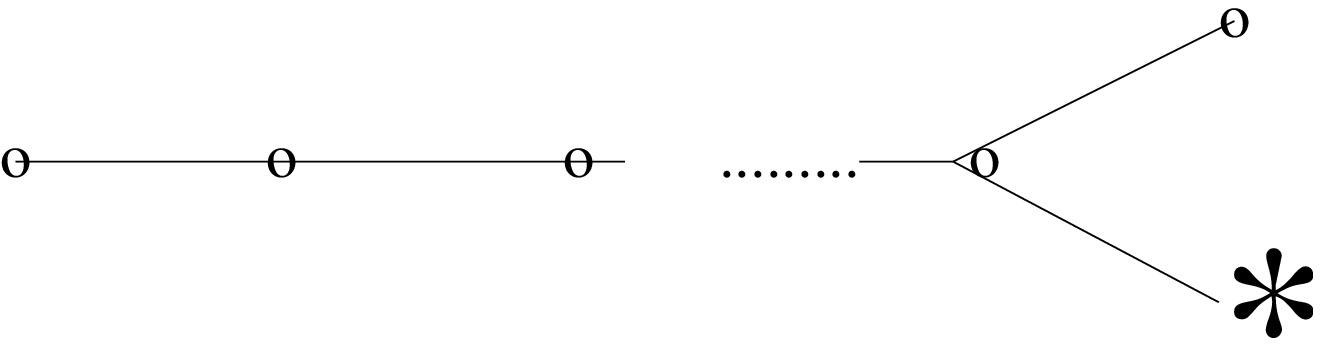} {1.0in} \\

$E_{6,1}$ \vpic{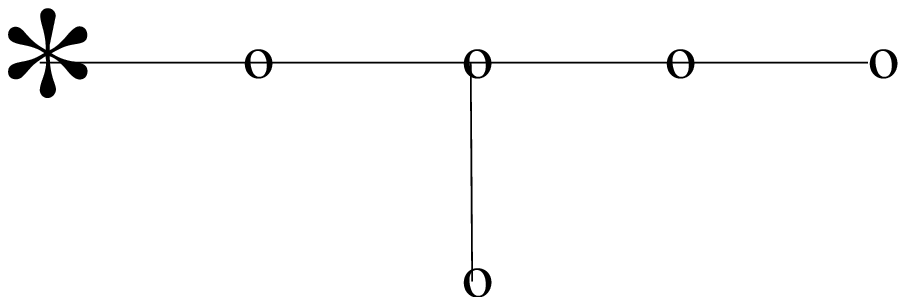} {1.0in} \quad $E_{6,2}$ \vpic{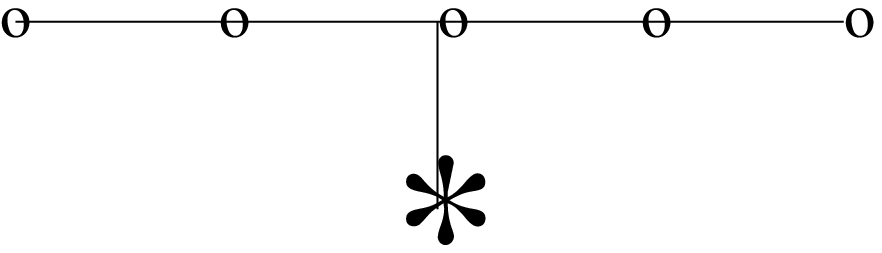} {1.0in} \\

$E_{7,1}$ \vpic{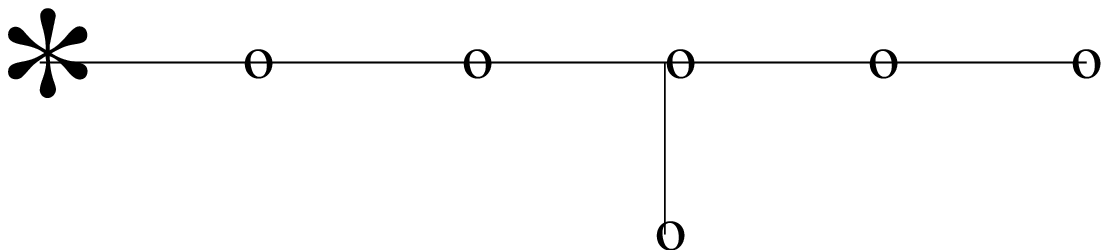} {1.0in}
 \quad $E_{7,2}$ \vpic{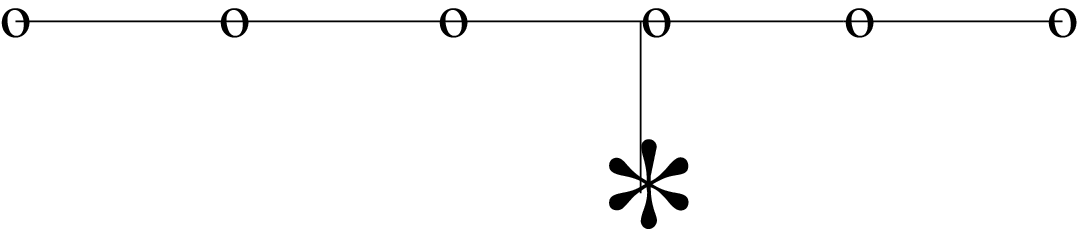} {1.0in} \quad $E_{7,3}$ \vpic{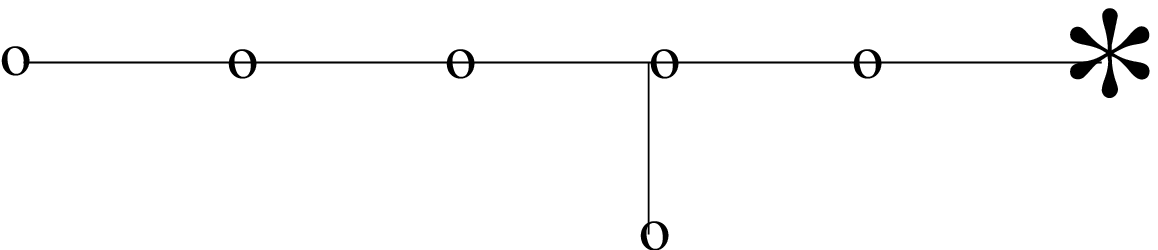} {1.0in} \\

$E_{8,1}$ \vpic{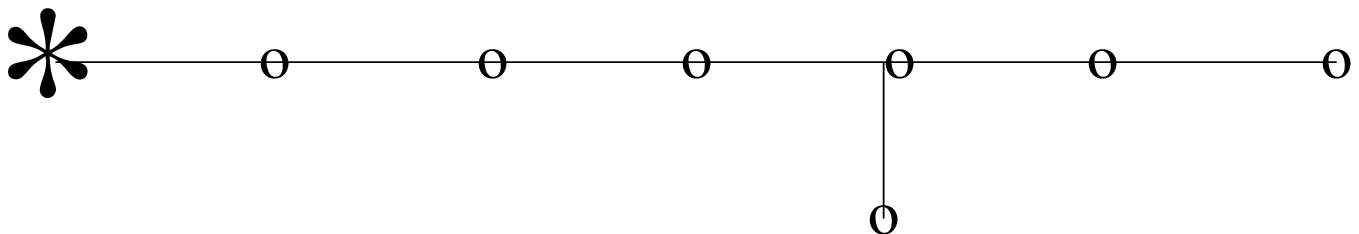} {1.0in}
 \quad $E_{8,2}$ \vpic{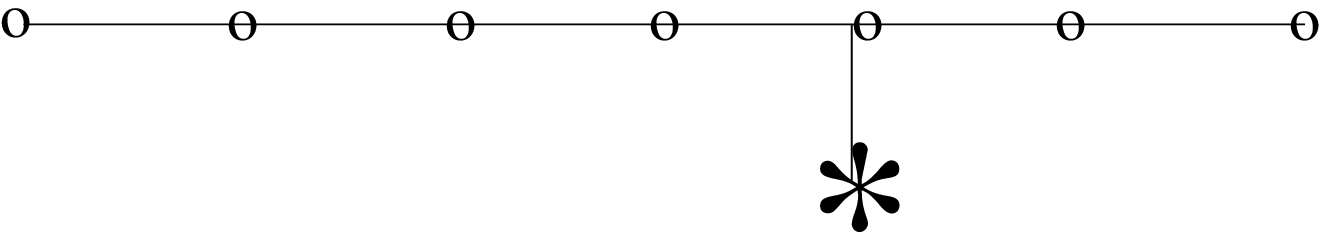} {1.0in} \quad $E_{8,3}$ \vpic{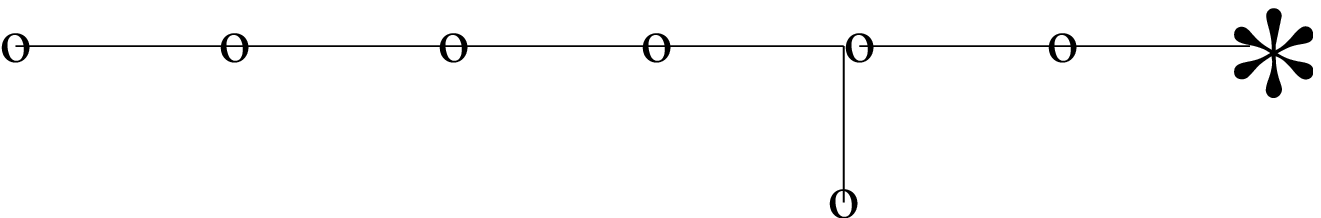} {1.0in} \\

\begin{proposition} The subfactor $pTL2 \subseteq M$ in the
case $D_{5,2}$ has index $(2+\sqrt 2)^2$
and principal graph \vpic{principalgraph} {1.5in} .
\end{proposition}
\begin{proof}
This is just the calculation done in example \ref{d5ghj1}.
\end{proof}

At this stage it looks very likely that the $D_{5,2}$ pair realises the 
case \ref{mainchance} of a no extra structure quadrilateral. In order to be sure 
of this we need to know that $P$ and $Q$ in this case do not commute.
To do this we shall compute the angle between them. At this stage we do
not even know if $P$ and $Q$ are distinct.

\subsection{Angle Computation}

 Our strategy for calculating the angle between $P$ and $Q$ will work whenever
the subfactors $TL2\subseteq P$ and $TL2\subseteq Q$ are 2-transitive so in this
subsection we only assume that of the Coxeter graph with chosen univalent vertex.

\begin{definition} Let $\Gamma$ be a pointed Coxeter graph of type $D_n$ for $n>4$ or
 $E$ on the
list above. Then $d=d(\Gamma)$ will denote the distance from * to the trivalent
vertex.
\end{definition}
Thus $d(E_{6,1})=2$ and $d(D_{5,2})=1$.

\begin{theorem}\label{anglecomp}
Suppose $\Gamma$ be a pointed Coxeter graph of type $D_n$ for $n>4$ or
 $E$, with Coxeter number $\ell$, and that the GHJ subfactor with the starred vertex is 2-transitive. 
Then the angle between the two intermediate subfactors is
 $$\{0,\pi/2, \cos^{-1}(|\frac{\cos \hspace {2pt}(2d+3)\pi/\ell}{\cos \pi/\ell}| )\}$$
\end{theorem}

\begin{proof}
The idea is as follows: by 2-transitivity $E_PE_QE_P$ is a multiple of the indentity
on the orthogonal complement of $TL2$ in $P$ so it suffices to find an element $x$ of
this orthogonal complement and calculate $||E_Q(x)||_2$. 
We will find our element $x$ in $pP_{d+2}p$ which is the smallest $pP_kp$ which is
strictly bigger than $pTL2_k$. It will be convenient to pull back the calculations
to $pA_np$ so in the next lemma we give the unitaries which conjugate $A_n$ to
$P_{n+1}$ and $Q_{n+1}$. These unitaries may be deduced from \cite{GHJ} but we give
a proof here for the convenience of the reader.

\begin{lemma}\label{shift}Let $v_n=g_1g_2...g_n$ and $w=g_1^{-1}g_2^{-1}...g_n^{-1}$.
Then \\(a) $P_{n+1}=v_nA_nv_n^*$ and $Q_{n+1}=w_nA_nw_n^*$. \\
\qquad(b) $TL2_n= v_n TL_n v_n^* = w_n TL_n w_n^*$  \end{lemma}

\begin{proof}
Braid group relations give $v_ng_iv_n^*=g_{i+1}$ and $w_ng_i^{-1}w_n^*=g_{i+1}^{-1}$
hence $v_ne_iv_n^*=e_{i+1}$ and $w_ne_iw_n^*=e_{i+1}$ for $1\leq i\leq n-1$.
This proves the assertion (b) about the Temperley-Lieb algebras. Since
$[e_i,A_1]=0$ for $i\geq 2$ we get $v_nA_1v_n^*=g_1A_1g_1^*=P_1$ and $w_nA_1w_n^*=Q_1$.   
By the definition of $P_n$ and $Q_n$ we are done.
\end{proof} 

As in figure \ref{d5ghj1} the Bratelli diagram for $pA_\infty p$ is
given by taking the full Bratelli diagram for $A_\infty$ and 
considering only edges emanating from the starred vertex. Thus by the
definition of $d(\Gamma)$ there is an element $y$ of $pA_{d+1}p$ which
is orthogonal to $TL_d$, and is unique up to a scalar multiple.
We may assume $||y||_2=1$ and $y=y^*$. Define $x\in P_{d+2}$ by $x=v_{d+1}yv_{d+1}^*$.\
By \ref{shift} we know that $x$ is orthogonal to $e_2,e_3,...e_{d+1}$.
Moreover since $tr(x)=0$ (since $x\perp 1$), $E_{P_{d+1}}(x)=0$ so
$e_{d+2}xe_{d+2}=0$ and taking the trace, $x\perp e_{d+2}$. By the usual
properties of the Markov trace in a tower, $x\perp e_n$ for $n>d+2$.
Thus $x\perp TL2$. 

Since the inclusions of $pQ_n p$ in $pA_n p$ are commuting squares we may 
calculate $E_Q(x)$ by $E_{pQ_{d+2}p}(x)$ (inside $pA_{d+2}p$ ).
%%$$ (Ad \hspace{3pt} w_{d+1})E_{pA_{d+1}p}(Ad \hspace{3pt}  w_{d+1}^*)(x).$$
But this element of $pQ_{d+2}p$ is orthogonal to $TL2$ so is a multiple of
$w_{d+1}yw_{d+1}^*$. So the cosine of the angle between $P$ and $Q$ is the
absolute value of the inner product
$$tr(xw_{d+1}yw_{d+1}^*)=tr(v_{d+1}yv_{d+1}^*w_{d+1}yw_{d+1}^*).$$

The algebras $pA_np$ are all included in the planar algebra for the bipartite
graph $\Gamma$ as defined in \cite{J20} so we may use the diagrams therefrom.
In particular the inner product we need to calculate is given by the 
following partition function (up to a power of $\delta =2\cos \pi/\ell$).\\
\begin{diagram}\label{innerproduct}
       \vpic {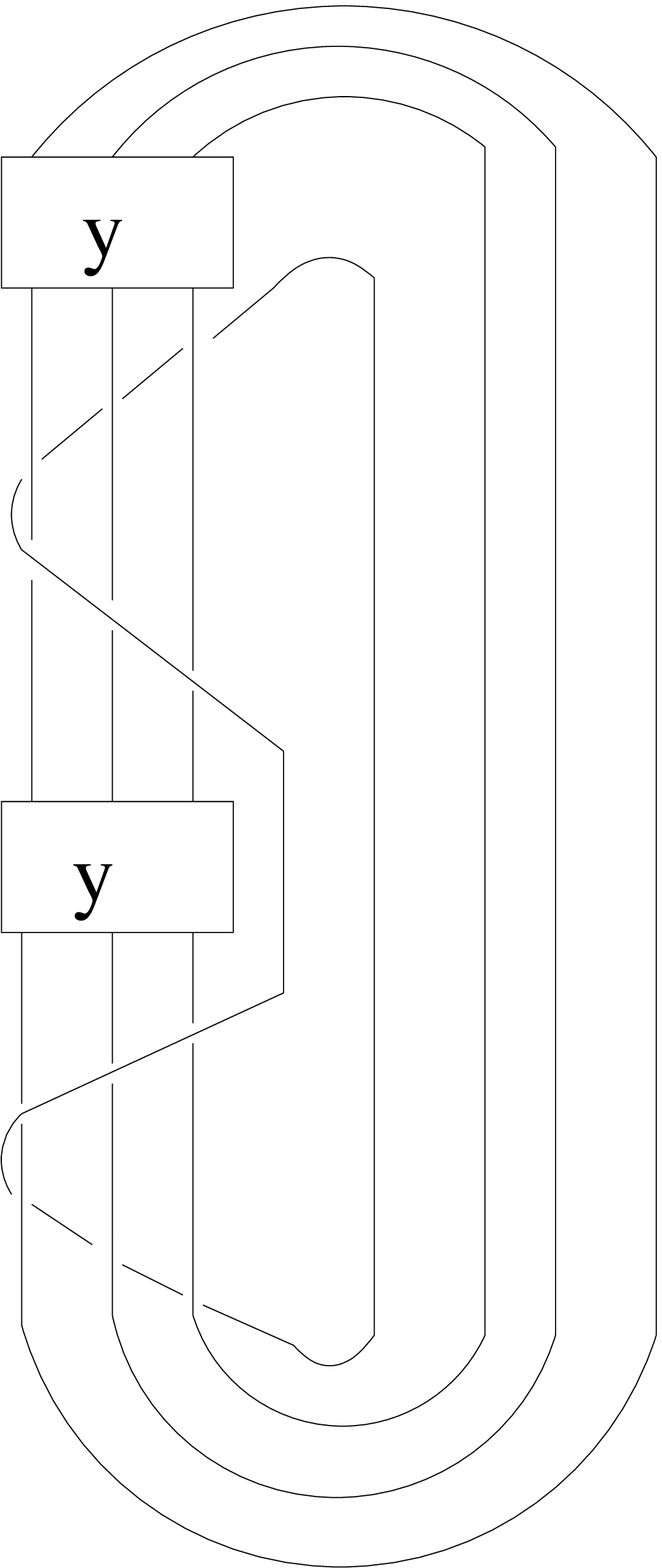} {1.5in} 
\end{diagram}  

The crossings in the picture are the braid elements $g_i$ with some convention
as to which is positive and which is negative, read from bottom to top. We have
illustrated with $d=2$ for concreteness.
They may be evaluated using the Kauffman picture:\\

\vpic{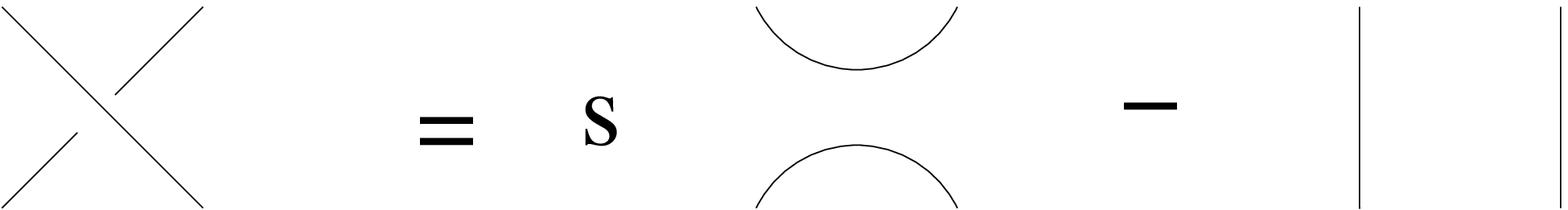} {2in} 

where $s=e^{\pi i/\ell}$.

The orthogonality of $y$ to $TL$ is equivalent to the fact that, if any
tangle contains a $y$ box with two neighbouring boundary points connected
by a planar curve(in which case we say the box is "capped off"),
the answer is zero. Thus one may evaluate \ref{innerproduct} as follows.

Using the Kauffman relation in \ref{innerproduct} inside the dotted circle
below one obtains\\

\begin{diagram}\label{tanglereduction}
\vpic{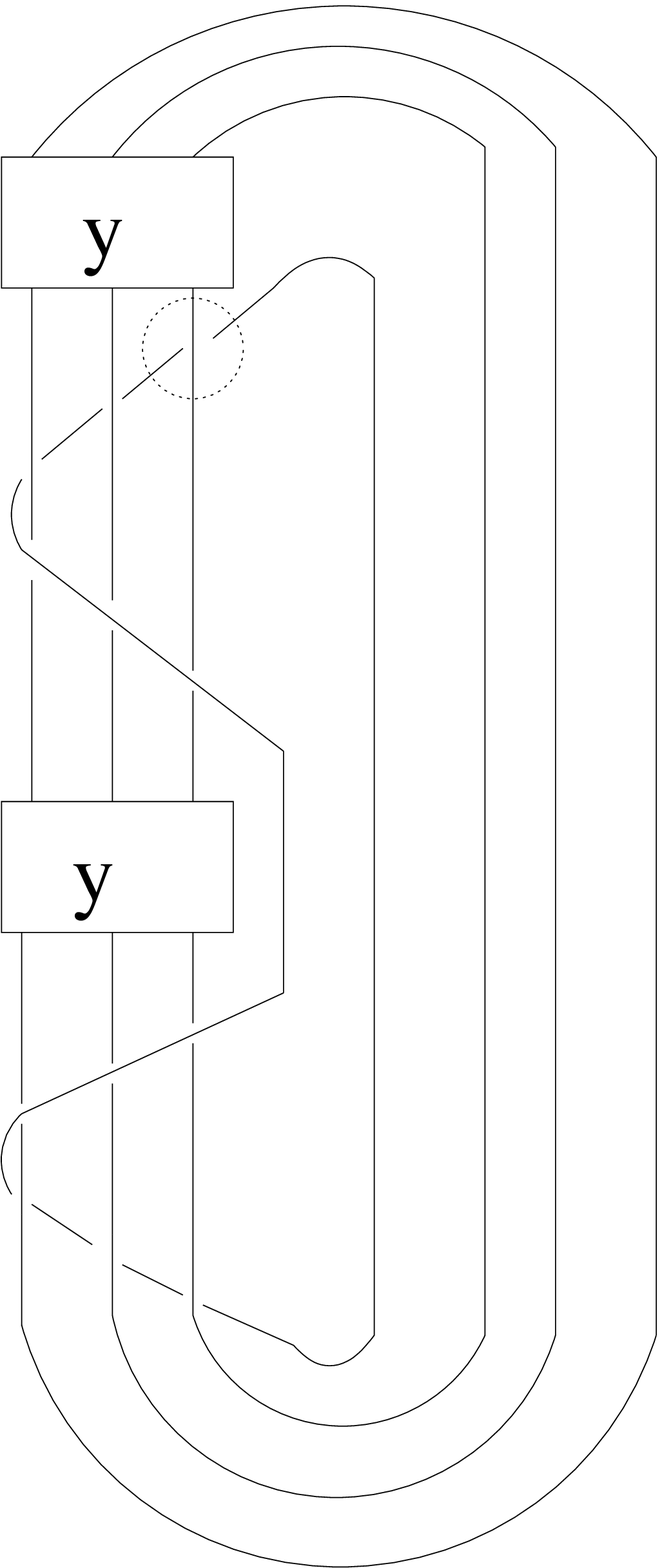} {1in} $= s$ \vpic{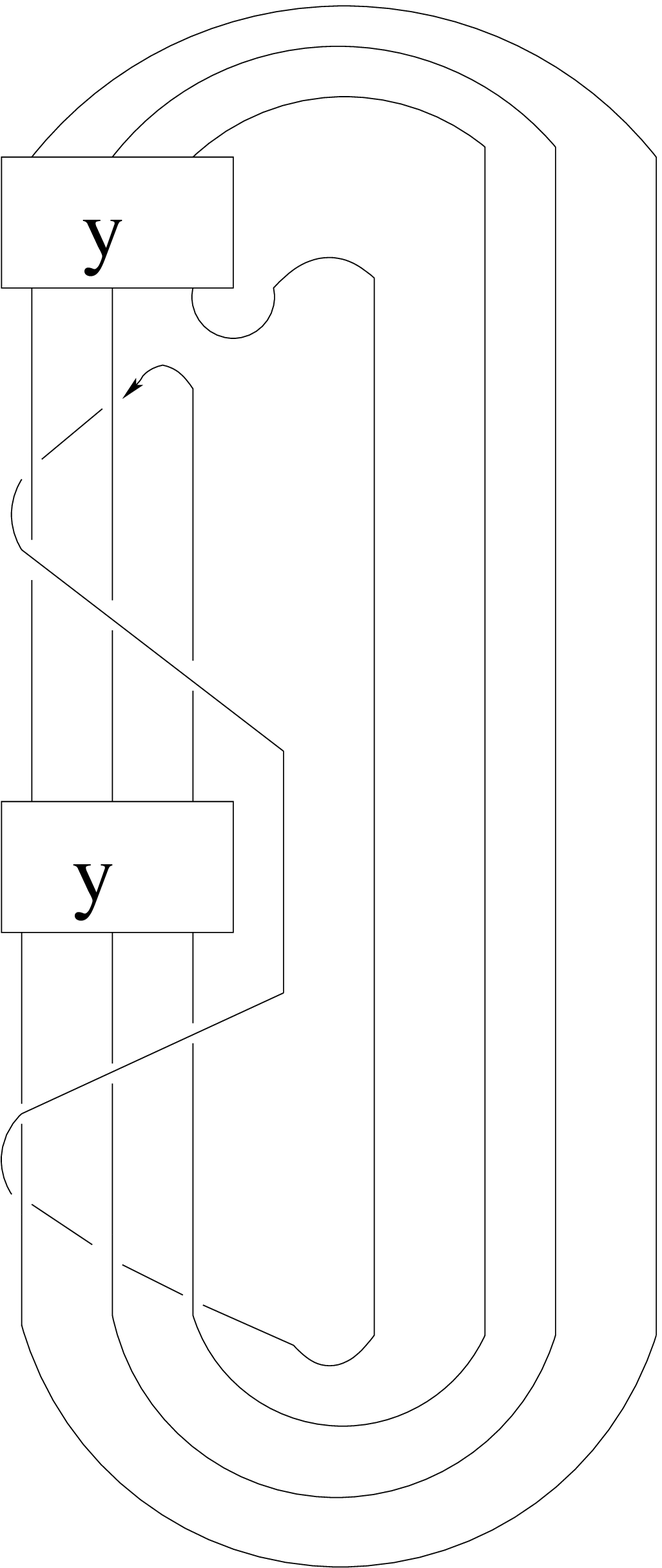} {1in} $-$ \vpic{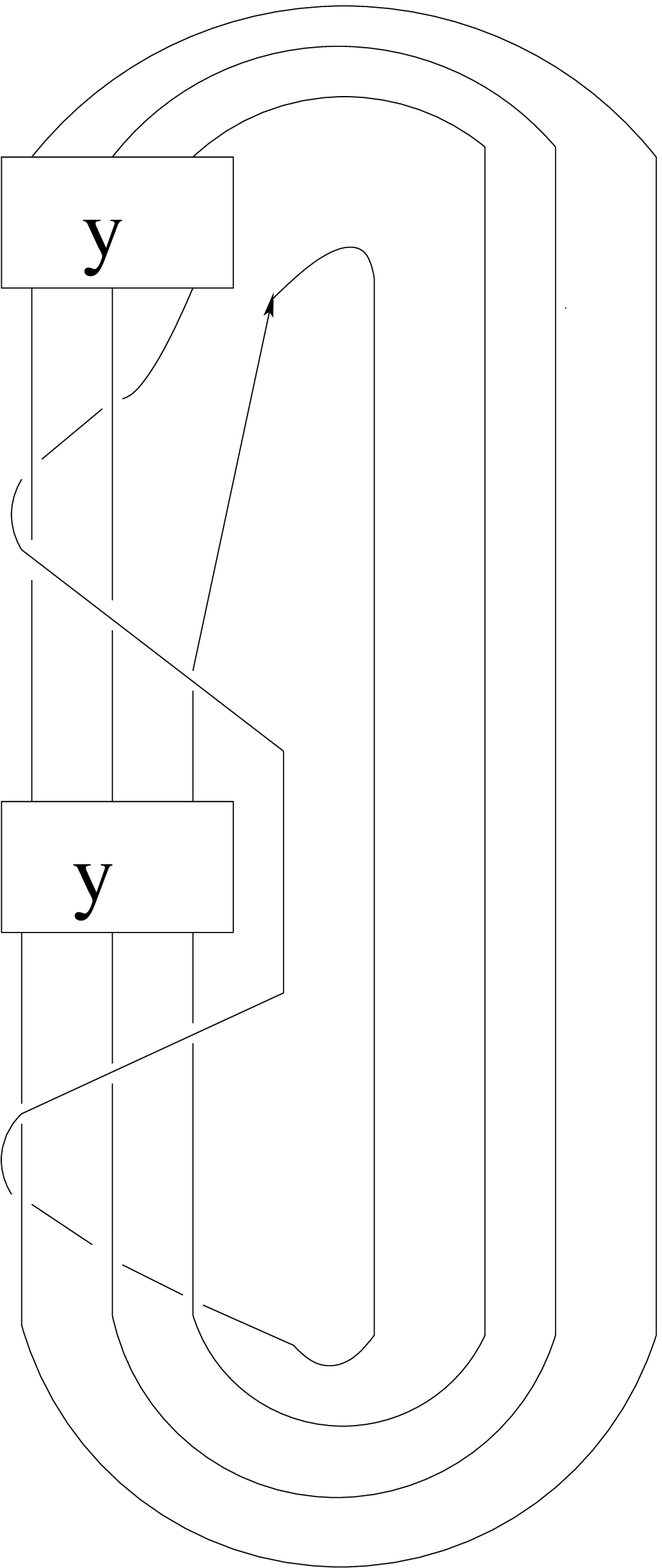} {1in} .
\end{diagram}
Consider the first diagram on the right hand side of the equation \ref{tanglereduction}. Following the 
curve in the direction indicated by the arrow, observe that one choice of the two
possibilities in applying the Kauffman relation at each crossing always results in 
one of the $y$ boxes being capped off. The first $d$ such crossings thus contribute
a factor of $s$ each. Then one meets the situation \vpic{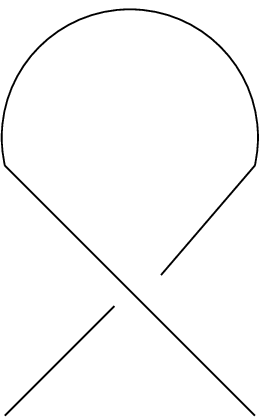} {0.3in} which is
easily seen to be the same as $s^2$ times \vpic{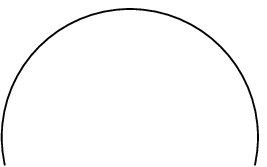} {0.2in} . One then meets
$d$ more crossings each of which which contributes $s$. After this (the crossings below the bottom $y$
box in \ref{innerproduct}) the only terms in the Kauffman relation contributing just
give the sign $-1$. Since there are an even number of them we deduce that 
the diagram of the first term on the right hand side of \ref{tanglereduction} is $s^{2d+2}$ times a
tangle which is $tr(y^2)$ up to a power of $\delta$. A similar analysis of the 
diagram of the second term gives $-s^{-(2d+3)}$ times  $tr(y^2)$. A little thought
concerning the powers of $\delta$ gives the final result that  

$$tr(v_{d+1}yv_{d+1}^*w_{d+1}yw_{d+1}^*)=\frac {s^{2d+3}+s^{-2d-3}}{s+s^{-1}}$$

This ends the proof of theorem \ref{anglecomp} 
\end{proof}
\begin{corollary} For the GHJ subfactor pair given by $D_{5,2}$, there is no extra 
structure, the
angle between $P$ and $Q$ is $\cos^{-1} (\sqrt 2 -1)$, 
and $P\cap Q = TL_2$. 
\end{corollary}
\begin{proof}
We have $[M:P]=4\cos^2\pi/8$ from the $D_5$ commuting square. Also $pTL2\subseteq P$ has
the same index from a GHJ calculation, or from the one already done for $D_5$. So there cannot
be subfactors between $pTL2$ and $P$ or $Q$, and $pTL2\subseteq P$ is 2-transitive.
So we can apply the previous theorem to get the angle. The only possible principal 
graph with index $4\cos^2\pi/8$ is $A_7$ so there is 
no extra structure.

\end{proof}

\subsection{A simpler quadrilateral with no extra structure.}\label{simpler}

Note that the definition of the GHJ pair will require the use of 
certain roots of unity. But at least in the $D_{n,2}$ case it is possible
to find another pair $\tilde P$ and $\tilde Q$ between $pTL2$ and $M$,
which is defined over $\mathbb R$! We will see that both $\tilde P$ and $\tilde Q$
form commuting cocommuting squares with both $P$ and $Q$. 
One of these two intermediate subfactors is quite canonical and exists whenever
$P\cap Q=TL_2$.
\begin{definition} Let $\Gamma$ etc. be as above. Let $\tilde P$ be the GHJ subfactor
for $p$, i.e. the subfactor generated by $pTL2$ and $pe_1$.
\end{definition}

\begin{proposition} The quadrilaterals $N\subset \tilde P, P \subset M$ and  
$N\subset \tilde P, Q \subset M$ are commuting squares.
\end{proposition}
\begin{proof} Reducing by $p$ is irrelevant so we can do the computation in the
full GHJ factor.
As in the proof of \ref{anglecomp} it suffices to find a non-zero element of $\tilde P$
orthogonal to $TL2$ and show that its projection onto $P$ is zero. Let $x = e_1 - \tau id$
where $\tau=(4\cos^2 \pi/ell)^{-1}$. Then
since the $P_n$'s form commuting squares with the $A_n$'s and $e_1\in A_2$ we need 
only project onto $P_2=Ad g_1(A_1)$. But $E_{P_2}=Ad g_1 E_{A_1} Ad g_1^{-1} $ 
and $Ad g_1 (x)=x$. But $E_{A_1}(x)=0$ is just the Markov property for the trace
on $A_2$. The same argument applies to $\tilde Q$.
\end{proof}

\begin{lemma}\label{effe}Let $\Gamma$ be $D_{n,2}$ for $n\geq 5$ or $E_{6,2}$. Then there
is a projection $f$ in $pA_2$ with the following properties.\\
(a)$tr(f)= \tau$.\\
(b)$fpe_1=0$.\\
(c)$pe_2 fpe_2 =\tau e_2$ and $fpe_2f=\tau f$.

 \end{lemma}
\begin{proof} From the Bratteli diagram for $pA_2$, it has three minimal
projections, which are central. One is clearly $pe_1$ and one of the other 
two has the same trace by symmetry. Let $f$ be that other one. Then (a) and
(b) are obvious. The first part of (c) follows from $\dim(pA_1)=1$ and the second
part follows since, from the Bratteli diagram, $f$ is a minimal projection in $pA_3$.
\end{proof}
\begin{definition}Let $\Gamma$ be $D_{n,2}$ for $n\geq 5$ or $E_{6,2}$. Let
$\tilde Q$ be the von Neumann algebra generated by $pTL_2$ and the $f$ of
\ref{effe}.
\end{definition}
\begin{theorem}
Let $\Gamma$ be $D_{n,2}$ for $n\geq 5$ or $E_{6,2}$. Then $\tilde Q$ is a 
II$_1$ factor with $[\tilde Q : pTL_2]=4\cos^2\pi/\ell$, and the angle between
$\tilde P$ and $\tilde Q$ is $\displaystyle \cos^{-1}(\frac{\tau}{1-\tau})$.
\end{theorem}
\begin{proof}
Lemma \ref{effe} and the properties of the basic construction show that $f$
has exactly the same commutation relations and trace properties 
with $pe_i$ for $i\geq 2$ as does $pe_1$. Thus by \cite{J3} $\tilde Q$ is a II$_1$ 
factor with the given index. Moreover the subfactor $pTL_2 \subset \tilde Q$ is
2-transitive so we can speak of \underline{the} angle between $\tilde P$ and
$\tilde Q$.

The angle calculation is not hard. As in \ref{anglecomp} it suffices to compute
the length of the projection onto $\tilde P_1$ of 
a unit vector in $\tilde Q$ orthogonal to $pTL_2$. By \ref{effe}, the 
element $x=f-\tau id$ is orthogonal to the two-dimensional algebra 
$pTL_2$ and $tr(x^*x)=\tau(1-\tau)$.
Since the $pTL_n$'s form commuting squares with the $pA_n$'s, $E_{\tilde P}(x)$
is just the projection $E(x)$ of $x$ onto $pTL_2$. By the bimodule property
of $E$, $E(x)pe_1=-\tau pe_1$ so $E(x)=\tau pe_1 +\lambda (p-pe_1)$.
Using $tr(x)=0$ we find $\displaystyle \lambda= -\frac{\tau^2}{1-\tau}$.
So $$||E(x)||^2=\tau^3 +(\frac{\tau^2}{1-\tau})^2(1-\tau)= \frac{\tau^3}{1-\tau}.$$
And finally $\displaystyle \frac{||E(x)||^2}{||x||_2^2}=\frac{\tau^2}{(1-\tau)^2}$.
\end{proof}
Observe that for $\tau^{-1}=4\cos^2\pi/\ell$, $\displaystyle \frac{\tau}{1-\tau}=\sqrt 2 -1$
so the angle between $\tilde P$ and $\tilde Q$ is indeed the same as that between $P$ 
and $Q$, and the quadrilateral formed by $\tilde P$ and $\tilde Q$ has no extra structure
for the same reasons as the one formed by $P$ and $Q$. As a last detail observe
that the quadrilaterals $N\subset \tilde Q, P \subset M$ and  
$N\subset \tilde Q, Q \subset M$ are commuting squares. We leave the argument to
the reader.
\section {Uniqueness.} \label{uniqueness section}
Outer actions of finite groups are extremely well understood so we need say nothing
more in the case $[M:N]=6$. Uniqueness up to conjugacy in the hyperfinite case follows
from \cite{J2}. 

So from now on we assume that $[M:N]=6+4\sqrt 2$ and that there are two 
intermediate subfactors $P$ and $Q$ which neither commute nor cocommute.
We will eventually show that all the constants in a planar algebra presentation of
the standard invariant of $N\subseteq M$ are determined by this data.

From the structure of the principal graph we see that there is exactly one projection in 
$N'\cap M_1$ different from $e_1$ but with the same trace as $e_1$.
 By \cite{PP2} this means that there is  a self adjoint unitary in the normaliser of $M$ in $M_1$
(and in the normaliser of $M_1$ in $M_2$). We record some useful diagrammatic facts 
about normalisers below. It is convenient to work with the normaliser of $M_1$ in $M_2$
but any subfactor is dual so the result can be modified for the normaliser of $M$.

\subsection{Diagrammatic relations for the normaliser.}
If $N\subseteq M$ is an irreducible finite index subfactor we will consider an element
$u$ in the normaliser of $M_1$ inside $M'\cap M_2$, that is to say a unitary in $M'\cap M_2$ 
with $uM_1 u^* = M_1$. First observe that such a unitary defines an automorphism 
$\alpha$ of $M_1$ by $\alpha (x) = uxu^*$. 
\begin{proposition} $\alpha(x)=x$ for all $x\in M$.
\end{proposition}
\begin{proof} Follows immediately from $u\in M'$.
  
\end{proof}
The automorphism $\alpha$ in turn 
defines a unitary on $L^2(M_1)$ which is in $M'\cap M_2$ and by irreducibility
differs from $u$ by a scalar. Thus we may alter $u$ so that $u=\alpha$ as
maps on $L^2(M_1)$.  The element $u$ is in $N'\cap M_2$ so in the planar algebra
picture it may be represented by a diagram: \vpic{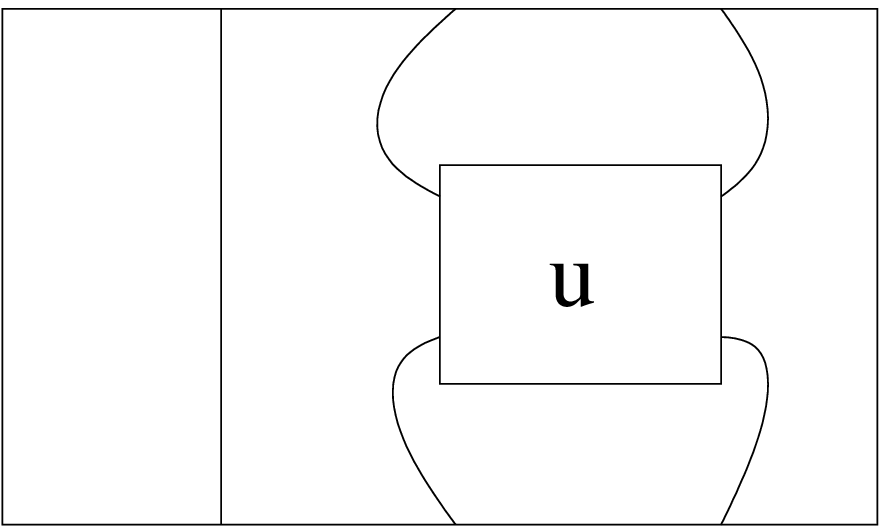} {1 in}  and the
relation $uxu^*=\alpha(x)$ for $x\in N'\cap M_1$ is the equality
\begin{diagram} \label{adu}\vpic{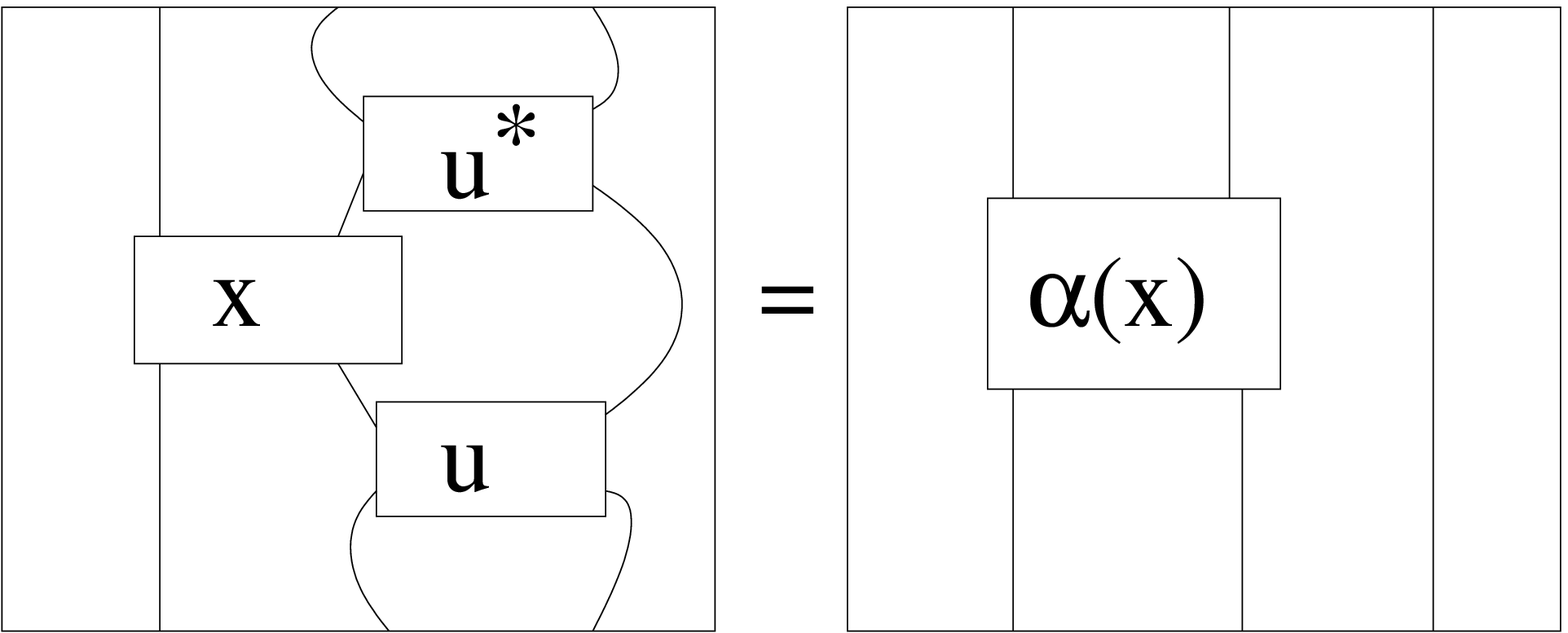} {2 in}
.
\end{diagram}
We will make considerable use of the following result:
\begin{lemma}\label{reduceu} If $u=u^*$ is in the normaliser as above then\\

\vpic{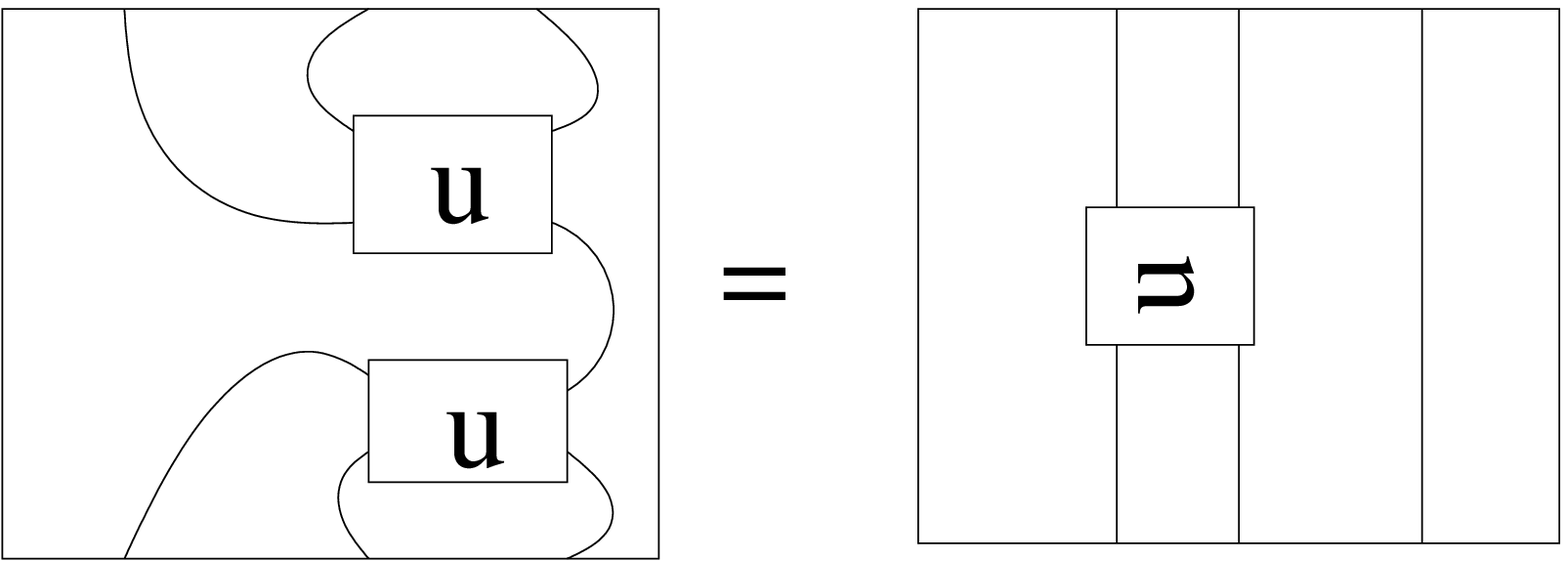} {2 in}  
\end{lemma}  
\begin{proof}We first establish that for any $u$ in the normaliser with 

$u=\alpha$ as above, and $x\in N'\cap M_1$,
\begin {diagram}\label {actionnormaliser}
\vpic{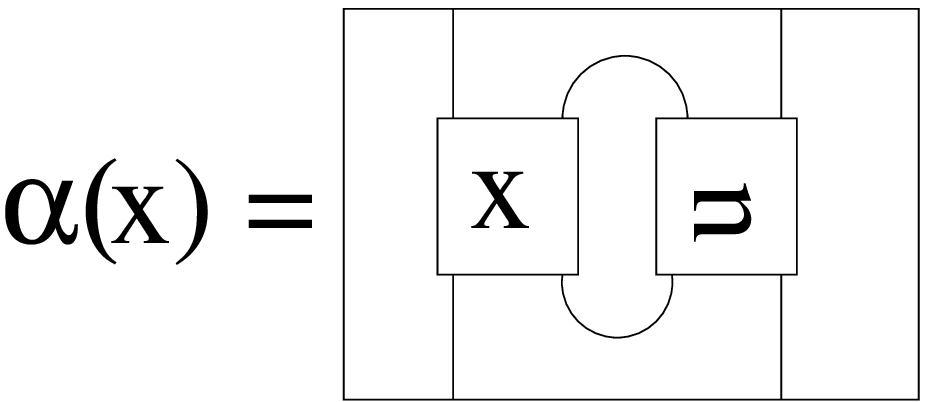} {1.7 in}  .\\
\end{diagram}
For this observe that if $a=xe_2y $ for $x, y\in M_1$ and $b\in M_1\subseteq L^2(M_1)$,
$E_{M_1}(abe_2)=\delta^{-2}xE_M(yb)= a(b)$. Since linear combinations of elements
of the form $xe_2y$ span $M_2$ we have 
$$E_{M_1}(abe_2)=\delta^{-2}a(b)$$ for all $a\in M_2$ and $b\in M_1$.
Drawing this relation diagramatically for $a=u$ and $b=x$ in $N'\cap M_1$ we obtain the diagram
for $\alpha (x)$. Finally apply \ref{adu} with $x=e_1$, and the above diagram 
to obtain the lemma.
\end{proof}
\begin{corollary}\label{ucoproj} With notation as above, $u$ is a coprojection.
\end{corollary}
\begin{proof} Use the property that $\alpha$ is a $*$-automorphism in
the previous lemma.
\end{proof}
\subsection{The structure of $N'\cap M_1$.}

We need to adopt some conventions for the position of certain operators in $N'\cap M_1$.
Since the angle between $P$ and $Q$ consists of one value (different from $0,\pi/2$),
we know that $e_P$ and $e_Q$ generate a $2\times 2$ matrix algebra modulo $e_N$.
We also know from the dual principal graph that there is an intermediate subfactor
$S$ with $[S:N]=2$. If $e_S$ is the projection onto $S$ then the trace of 
$e_S$ is $\frac{2}{6+4\sqrt 2}$ and it is $e_N$ plus a minimal projection
in $N'\cap M_1$. This means that $e_S$ must be orthogonal to both of the 
$2\times 2$ matrix algebras in $N'\cap M_1$ since the traces of minimal projections
therein do not match.

\begin{definition} We write $N'\cap M_1= e_N \mathbb C \oplus A \oplus B \oplus (e_S-e_N) \mathbb C$
where $A$ and $B$ are $2\times 2$ matrix algebras with $e_P A \neq 0$.
\end{definition}
This definition specifies $A$ uniquely since $tr(e_P) = (2+\sqrt 2)^{-1}, tr(E_N)=(2+\sqrt 2)^{-2}$
and the trace of a minimal projection in $A$ is $\frac{1+\sqrt 2}{(2+\sqrt 2)^2}$. Thus
$e_P B=0$.

\subsection{Relations between elements in $N'\cap M_1$}

From \ref{mainchance} we know that the principal and dual principal graphs are the same
and that there is a single projection of trace equal to that of $e_N$ in all
the (second) relative commutants. This means by \cite{PP2} that each for each inclusion
$M_i\subset M_{i+1}$ there is an intermediate inclusion $R_i$ 
with $[R_i:M_i]=2$. By duality there are thus $S_i$ with $M_i\subset S_i\subset M_{i+1}$
so that $S_i\subset M_{i+1}\subset R_{i+1}$ is a fixed point/crossed product pair
for an outer action of $\mathbb Z/2\mathbb Z$. In particular there are unitaries 
$u_i$ satisfying the conditions of the previous section at every step in the towwer.
 So let  $\alpha$ be the period two automorphism of $M$
(which is the identity on $N$) defining an element $u$ of $N'\cap M_1$.
Then $\frac{u+1}{2}$ is the projection onto an intermediate subfactor of
index 2 for $N\subset M$ which we shall call $R$. Thus\\
\vskip 5pt
$[M:R]=2$ or $tr(e_R)=\frac{1}{2}$\quad, and
$u=2e_R - 1$.\\

\begin{lemma}\label{position of R}The subfactors $P$ and $R$ cocommute but do not commute,
$e_P e_R e_P = e_{N} + (1-\frac{1}{\sqrt{2}}) (e_P-e_N) $ and $e_R B\neq 0$.
\end{lemma}
\begin{proof}
Since $L^{2}(M_{1}) \cong U_{0} \oplus 2U_{1} \oplus 2U_{2} \oplus U_{3}$ 
as $M-M$ bimodules, where 
$L^{2}(\bar{P}) \cong U_{0} \oplus U_{1}$ and $L^{2}(\bar{R}) \cong U_{0} \oplus U_{3}$,
  the dual subfactors $\bar{P}$ and $\bar{R}$ commute. However, $[\bar{P}:M][\bar{R}:M]<[M_1 :M]$ so
  by \ref{commcocomm}  $\bar{P}$ and $\bar{R}$ do not cocommute. Thus $P$ and $R$ cocommute but do not commute.
Then $L^{2}(R)$ must be of the form $V_{0} \oplus V_{1} \oplus V_{2}$, so $e_{R} B\neq 0$.
 Since $N \subset P$ is 2-supertransitive, by 5.3.1 
we have $\displaystyle e_P e_R e_P =e_N + \frac{tr(e_{\bar{P}\bar{R}})^{-1}-1}{[P:N]-1}
 (e_P-e_N)$. Since the dual quadrilateral commutes, by 4.1.1 we 
have $\displaystyle tr(e_{\bar{P}\bar{R}})=\frac{[\bar{P}:M][\bar{R}:M]}{[M:N]}=\frac{2}{2+\sqrt{2}}$. 
Combining these equations gives the result.
\end{proof}

We want to investigate the algebraic and diagrammatic relations between
$e_P$, $e_Q$ and $u$.
First a simple but crucial computation:
\begin{lemma}\label{orthogonality of u} $tr(ue_P)=tr(ue_q)=0$.
\end{lemma}
\begin{proof} Since $P$ and $R$ cocommute, by \ref{mulfor} $tr(e_Pe_R)=tr(e_P)tr(e_R)=1/2 tr(e_P)$,
and $u=2e_R-1$.
\end{proof}
We will use on several occasions the following result which is no doubt extremely 
well known. We include a proof for the convenience of the reader.

\begin{lemma}\label{four subspaces}
 Let $P,Q,R,S$ be distinct projections onto four one-dimensional subspaces of $\mathbb C^2$
all making the same angle with respect to one another. Then that angle is $\cos^{-1} \frac{1}{\sqrt 3}$.
\end{lemma}
\begin{proof} If we choose a basis so that $P= \left( \begin{array}{cc} 
1&0 \cr
0&0 \end{array} \right)$ then any other projection at $\cos^{-1}(\sqrt a)$ to $P$ is
of the form  $P= \left( \begin{array}{cc} 
a&\omega \sqrt{a(1-a)} \cr
\omega^{-1}\sqrt{a(1-a)} &1-a \end{array} \right)$ where $|\omega |=1$.  Equating $a$ to 
the traces of $QR$, $RS$ and $QS$ we see that $\omega$ must be a proper cube root of unity 
and that $3a^2-4a+1=0$.
\end{proof}

\begin{corollary}\label{switchPQQP} $ue_P u=e_Q$ and $ue_{PQ}u=e_{QP}$.
\end{corollary}
\begin{proof}
These are equivalent to $\alpha(P)=Q$. By \ref{position of R} $ue_Pu\neq P$.
If $\alpha(P)$ were
not equal to $Q$ then $P,Q,\alpha(P)$ and $\alpha(Q)$ 
are four distinct intermediate subfactors. But $ue_Pu=e_{\alpha (P)}$
and $ue_Qu=e_{\alpha(Q)}$ so the 
$N-N$ bimodules defined by these four intermediate subfactors  are all isomorphic to $L^2(P)$
and none of them commutes with any other. By \ref{position of R} which guarantees that
$\alpha(P)$ and $P$ do not commute,  the angles between all four subfactors are
the same and, by \ref{anglePQ}, equal to $\cos^{-1}(\sqrt 2 - 1)$. By \ref{four subspaces} this is
impossible. 
\end{proof}
\begin{corollary}\label{uep} $ue_P=e_N + \frac{1}{1-\sqrt 2}(e_Qe_P-e_N)$
 and \\$ue_Q=e_N + \frac{1}{1-\sqrt 2}(e_Pe_Q-e_N)$.
\end{corollary}
\begin{proof}$u(e_P-e_N)$ and $e_Q(e_P-e_N)$ are in $A$ and both multiples of a partial isometry
with intial domain $e_P-e_N$ and final domain $e_Q-e_N$. They are thus proportional. Taking
the trace we get the result using \ref{orthogonality of u} and \ref{anglePQ}.
\end{proof} 
This yields a different derivation of the angle between $P$ and $Q$. We see that modulo
the ideal $\mathbb C e_N$ we have $ue_P=\frac{1}{1-\sqrt 2} e_Qe_P$ so that
mod this ideal $e_P=e_Puue_P=(\frac{1}{1-\sqrt 2})^2 e_Pe_Qe_P$ which determines
the constant in the angle formula $e_Pe_Qe_P -e_N = constant(e_Q-e_N)$.

\begin{corollary}\label{idA} The identity $1_A$ of the $2\times 2$ matrix algebra $A\subseteq N'\cap M_1$
is $\frac{\sqrt 2 +1}{2}(e_P+e_Q)+1/2(ue_P+ue_Q)-(2+\sqrt2)e_N$ 
\end{corollary}
\begin{proof} From \ref{anglePQ},
 $(e_P-e_N)(e_Q-e_N)(e_P-e_N)= (\sqrt 2 -1)^2(e_P-e_N)$ so 
$1_A =\frac{\sqrt 2 +1}{2} (e_P-e_Q)^2$. \ref{uep} gives $e_Pe_Q=\sqrt 2 e_N + (1-\sqrt 2) ue_Q$
hence the result.
\end{proof}
\begin{lemma}$tr(ue_{PQ})=0$
\end{lemma}
\begin{proof} Since $u=2e_R-1$, $tr(ue_{PQ})=2tr(e_Re_{PQ})-1/\sqrt 2$ by \ref{anglePQ}.
But $tr(e_Re_{PQ})$ is given by $\displaystyle \frac{1}{\delta^3 tr(e_Pe_Q)}$ times the following diagram:
\vpic{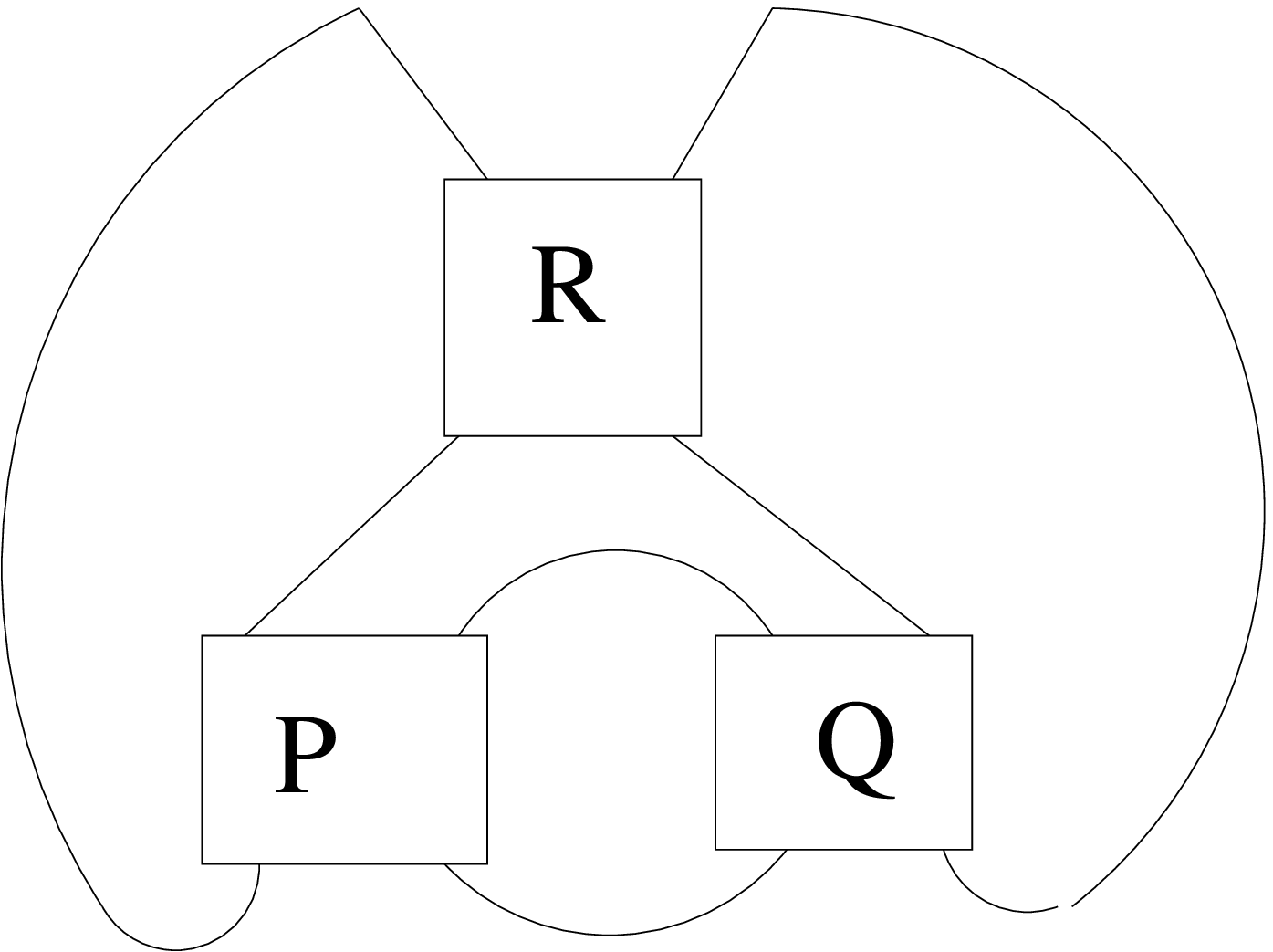} {1.0 in}
. This is essentially the cotrace of $e_R\circ e_P\circ e_Q$, and we know that 
$e_R\circ e_P$ is $(2+\sqrt 2)tr(e_R)tr(e_P)id$ by \ref{landauproj}
since $P$ and $R$ cocommute. Using this in the figure we obtain 
$tr(e_Re_{PQ})=\displaystyle \frac{1}{\delta^3 tr(e_Pe_Q)}(2+\sqrt 2)tr(e_R)tr(e_P)\delta^2 tr(e_Q)
 = \frac{1}{2\sqrt 2}$.
\end{proof}
\begin{lemma}$tr(e_{PQ}e_{QP})=\frac{5\sqrt 2 - 6}{2}$\end{lemma}
\begin{proof} As in \ref{trepqeqp} we recognise $tr(e_{PQ}e_{QP})$ as being $\frac{1}{2[M:N]}$ times  the 
cotrace of $e_P\circ e_Q\circ e_P\circ e_Q$. But since $[M:P]=[P:N]$, $e_P$ and
$e_Q$ are coprojections and the angles between them as coprojections are the same as the angles
between them as projections. So $tr(e_{PQ}e_{QP})=\frac{1}{2}tr((e_Pe_Qe_P)^2)$.
However from \ref{anglePQ} $e_Pe_Qe_P=e_N+\frac{\sqrt 2 -1}{\sqrt 2 + 1}(e_P-e_N)$.
Squaring and taking the trace gives the answer.

\end{proof}
\begin{corollary}\label{uePQ} $ue_{PQ}=e_N+u1_A -(\sqrt 2 +1)( e_{QP}e_{PQ}- (e_N +  1_A))$.
\end{corollary}
\begin{proof}
As in \ref{uep}, $u(e_{PQ}-e_N-1_A)$ and $e_{QP}e_{PQ}-e_N-1_A$ are both in $B$ (certainly
$e_{PQ}>e_Q$ and the trace of $e_{PQ}$ is the trace of $e_N$ plus 3 times the trace
of a minimal projection in $A$ so that $e_{PQ}e_S=0$) 
and are multiples of a the same partial isometry. Taking the trace using the
last two lemmas we get  $u(e_{PQ}-e_N-1_A)=\frac{3+2\sqrt 2}{\sqrt 2 -1}(e_{QP}e_{PQ}-e_N-1_A)$ 
 and the result
follows.

\end{proof}
\begin{corollary}\label{ePQQPPQ}
$e_{PQ}e_{QP}e_{PQ}-1_A-e_N=(\sqrt 2 -1)^2(e_{PQ}-1_A-e_N)$
\end{corollary}
\begin{proof} Modulo the ideal spanned by $e_N$ and $A$,
 $ue_{PQ}=-(\sqrt 2 +1)e_{QP}e_{PQ}$. So mod this ideal
 $e_{PQ}uue_{PQ}= (\sqrt 2 +1)^2e_{PQ}e_{QP}e_{PQ}$.
 The left and right hand sides are proportional and this determines
 the constant.
\end{proof}
Taking the trace of this equality provides a useful check on our
calculations. It is curious that $e_{PQ}$ and $e_{QP}$ make
the same angles as $e_P$ and $e_Q$.

\subsection{A basis and its structure constants.}

\begin{definition}Let {\gothic {C}} $=\{e_N,1\} \cup  ${\gothic{A}} $\cup ${\gothic{B}}
where {\gothic{A}}$=\{e_P,e_Q,ue_P,ue_Q\}$ and {\gothic{B}}$=\{e_{PQ},e_{QP},ue_{PQ},ue_{QP}\}$.
\end{definition}
\begin{theorem}\label{basis1} {\gothic C} is a basis for $N'\cap M_1$ and all multiplication and 
comultiplication structure constants for this basis are determined.
\end{theorem}
\begin{proof} That {\gothic C} is a basis follows easily from the previous results-
$\{e_N\} \cup ${\gothic A} is a basis for $\mathbb C e_N \oplus A$ by \ref{uep} and
$2\times 2$-matrix calculations. Similarly {\gothic B} forms a basis for 
$B$ modulo $\mathbb C e_N \oplus A$ by \ref{ePQQPPQ}. The identity spans $N'\cap M_1$
modulo $\mathbb C c_P \oplus A \oplus B$.

With the results so far, it is easy to see that all the structure constants for multiplication are
determined: Multiplication of any basis element by $e_N$ produces $e_N$. Multiplication
within {\gothic A} is determined by \ref{uep} and \ref{anglePQ}. Similarly multiplication
within {\gothic B} is determined by \ref{uePQ}, \ref{ePQQPPQ} and the explicit form
of $1_A$ in \ref{idA}. This leaves only multiplication between {\gothic A} and {\gothic B}.
But $e_{PQ}e_P=e_P$ (and other versions with $P$ and $Q$ interchanged) takes care of this.
Note also that {\gothic C}={\gothic C}$^*$ so that the $*$-algebra structure of 
$N'\cap M$ is explicitly determined on the basis {\gothic C}. 

We now turn to comultiplication. The $*$ structure for comultiplication
is rotation by $\pi$ and insertion of $*$'s of elements. Inspection shows that
the basis {\gothic C} is stable under this operation since $u=u^*$ is a projection for
comultiplication by \ref{ucoproj}.
The subsets {\gothic A} and {\gothic B} no
longer correspond to the algebraic structure but it will be convenient to organise
the calculation according to them. Determination of all the structure constants 
will just be a long sequence of cases, the most difficult of which will be 
diagrammatic and make frequent use of \ref{reduceu}. Note that the shading of
the picture will be 
the opposite of that in \ref{reduceu} since $u$ is in $M_1$ and not in $M_2$.
Occasionally the diagrammatic reductions will produce the element $u$ itself.
It is easy to express $u$ as a linear combination of basis elements since
$u(1-e_N-1_A-1_B)=1-e_N-1_A-1_B$ and $u$ times any element of {\gothic A}$\cup${\gothic B}
is another element of {\gothic A}$\cup${\gothic B}.

We will also use the exchange relation for biprojections from \cite{Bs3}:\\ \vpic{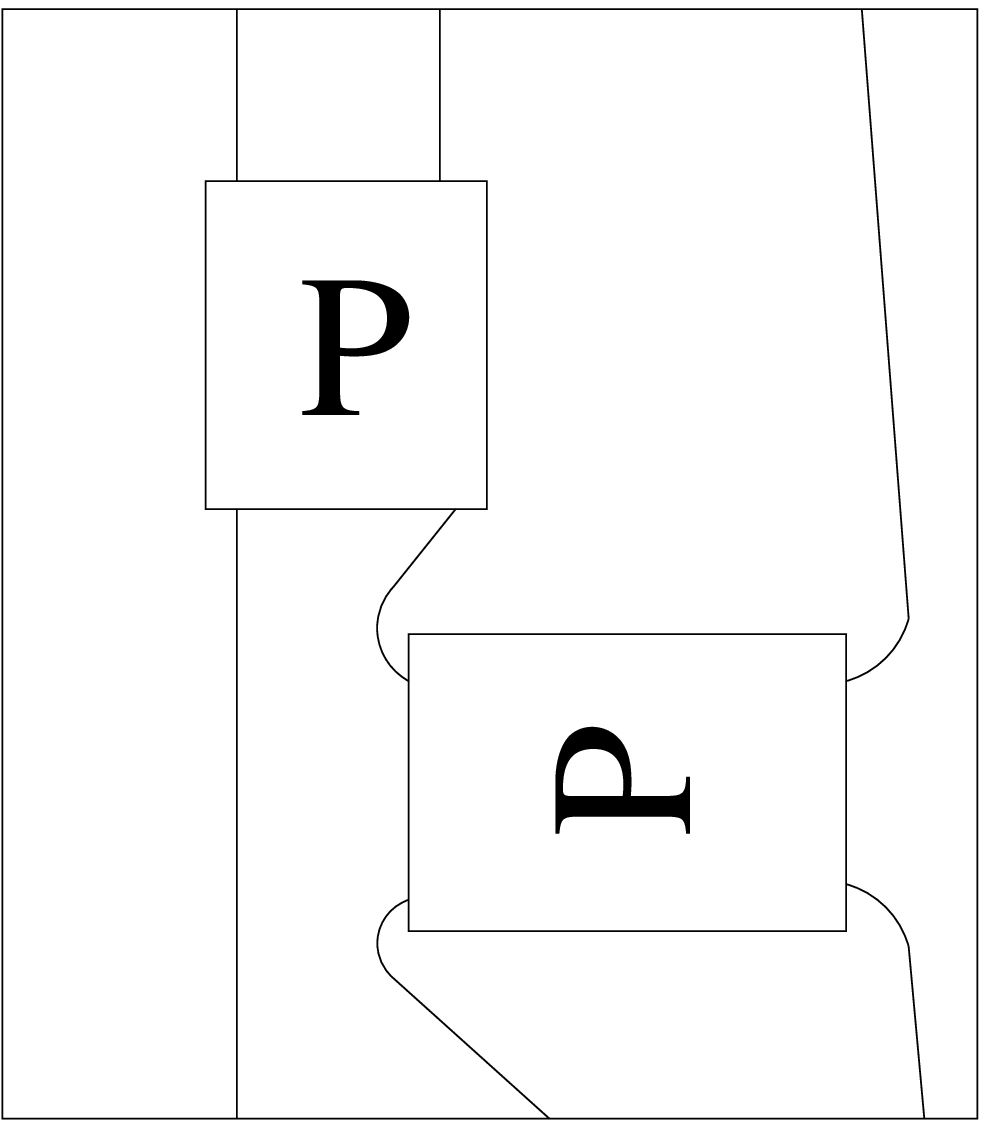} {0.8 in}  $=$
\vpic{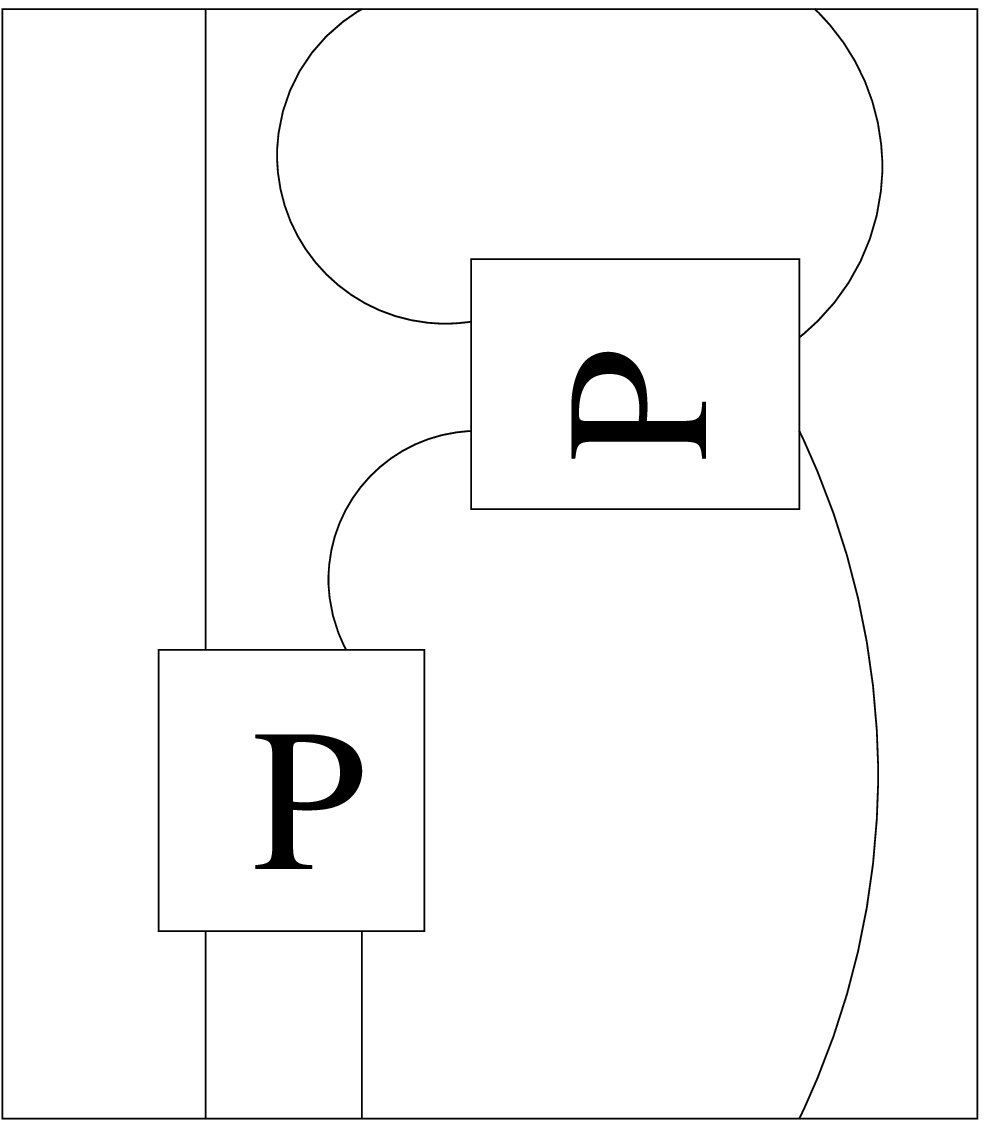} {0.8 in}  .

We have no need for the exact values of the structure constants, we only need to
know that they could be calculated explicitly. Thus we introduce the notation
$x\approx y$ to mean that the elements $x$ and $y$ of $N'\cap M_1$ are equal
up to multiplication by a constant that could be calculated explicitly.

Thus for instance $e_N\approx \tilde 1$ when $\tilde 1$ is the identity for comultiplication.
So all structure constants for comultiplication by $e_N$ are determined.
Comultiplication by $1$ is easy by the formula $x\circ 1 \approx tr(x)1$ for
$x\in N'\cap M_1$ and the only
trace that requires any work at all is that of $ue_{PQ}$ which is 
determined from \ref{uePQ} and \ref{ePQQPPQ}.
\vskip 5pt
Case 1. Comultiplication within {\gothic A}.
We may replace $ue_P$ by $e_Qe_P$ which is $\approx$ the projection onto
$L^2(PQ)$ for comultiplication. It is thus greater than $e_P$ and $e_Q$
so $e_P\circ (e_Pe_Q)\approx e_P$. The first case where any work is required is $(ue_P)\circ (ue_Q)$ and
up to simple modifications of the argument this handles all comultiplications 
within {\gothic A}. The labelled tangle defining $(ue_P)\circ (ue_Q)$ is: 
\hbox{ \qquad \qquad  \qquad \qquad                                            \vpic{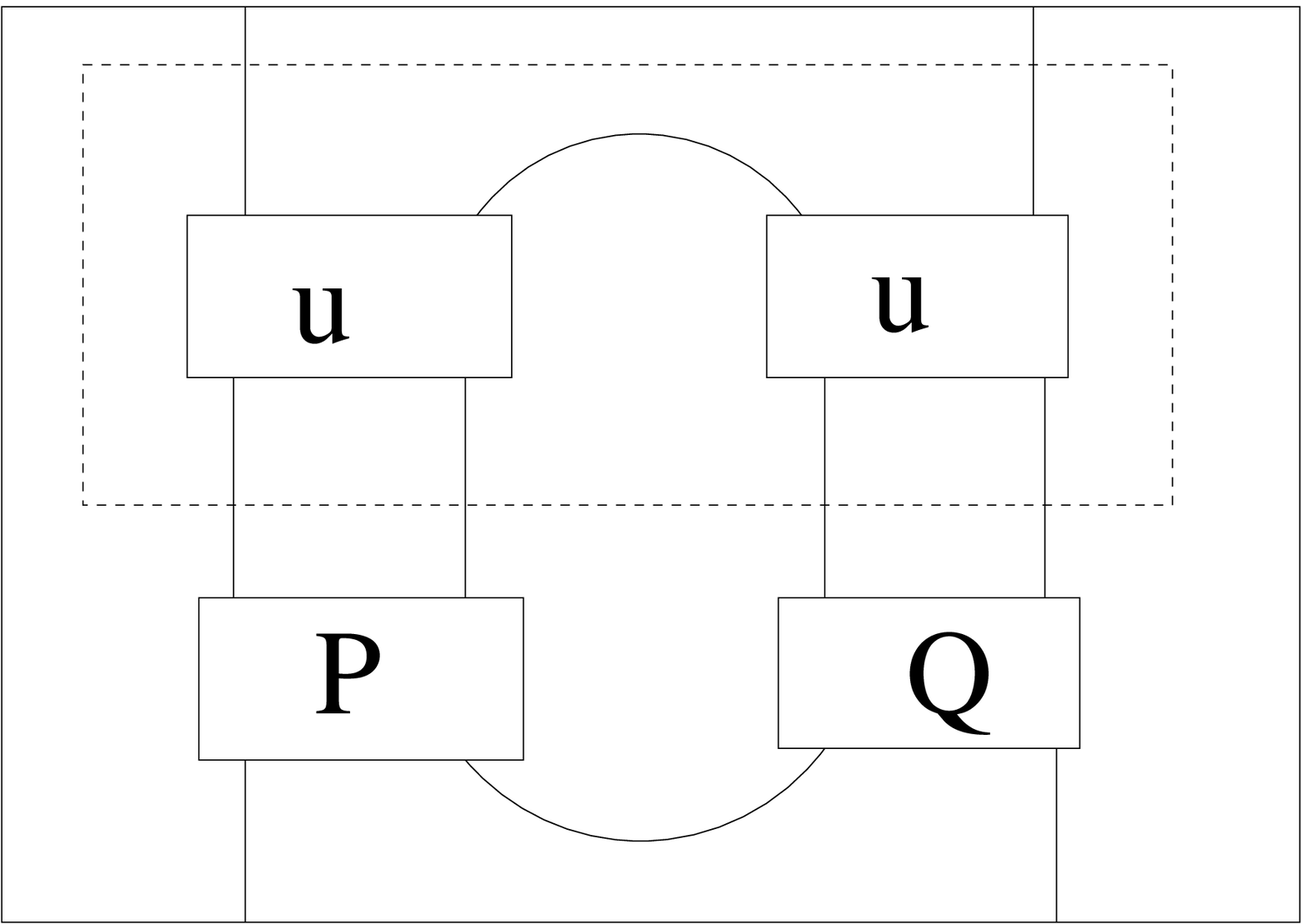} {1.5in}   }.
\vskip 5pt
Applying \ref{reduceu} to the region inside the dotted rectangle we obtain:
\vskip 5pt
\hbox{ \qquad \qquad  \qquad \qquad                                            \vpic{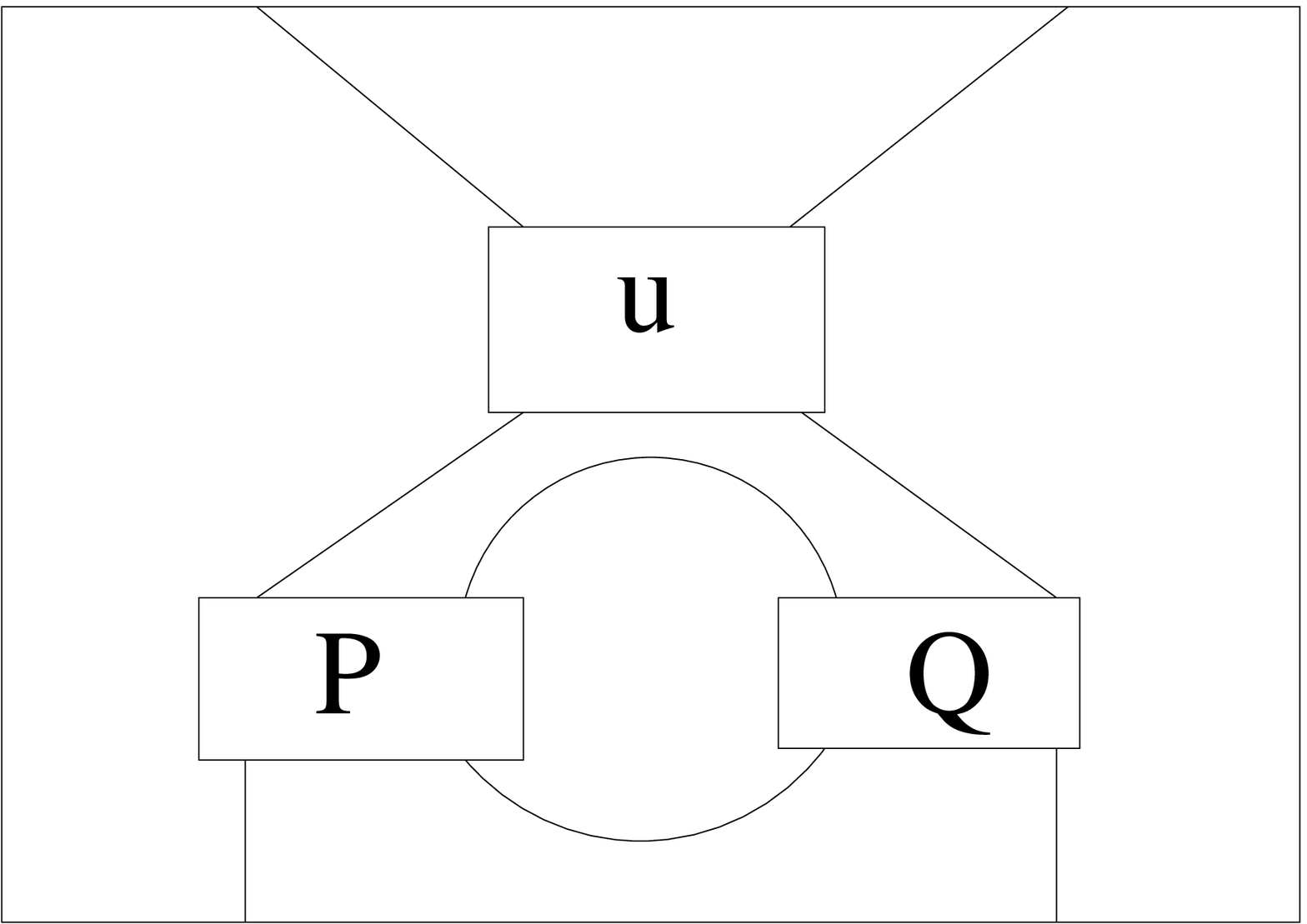} {1.5in}   .}

But this is $\approx e_{PQ}u$ which is a basis element.
\vskip 5pt
Case 2. Comultiplication within {\gothic B}.

Comultiplying $e_{PQ}$ with itself or with $e_{QP}$ is easy since under comultiplication
$e_P$ and $e_Q$ generate a $2\times 2$ matrix algebra mod $1$ and $e_P \circ e_Q \approx e_{PQ}$.
Comultiplying $e_{PQ}$ or $e_{QP}$ with $ue_{PQ}$ or $ue_{PQ}$ can, after applying 
\ref{switchPQQP} if necessary, a labelled tangle like:
\hbox{ \qquad \qquad                                            \vpic{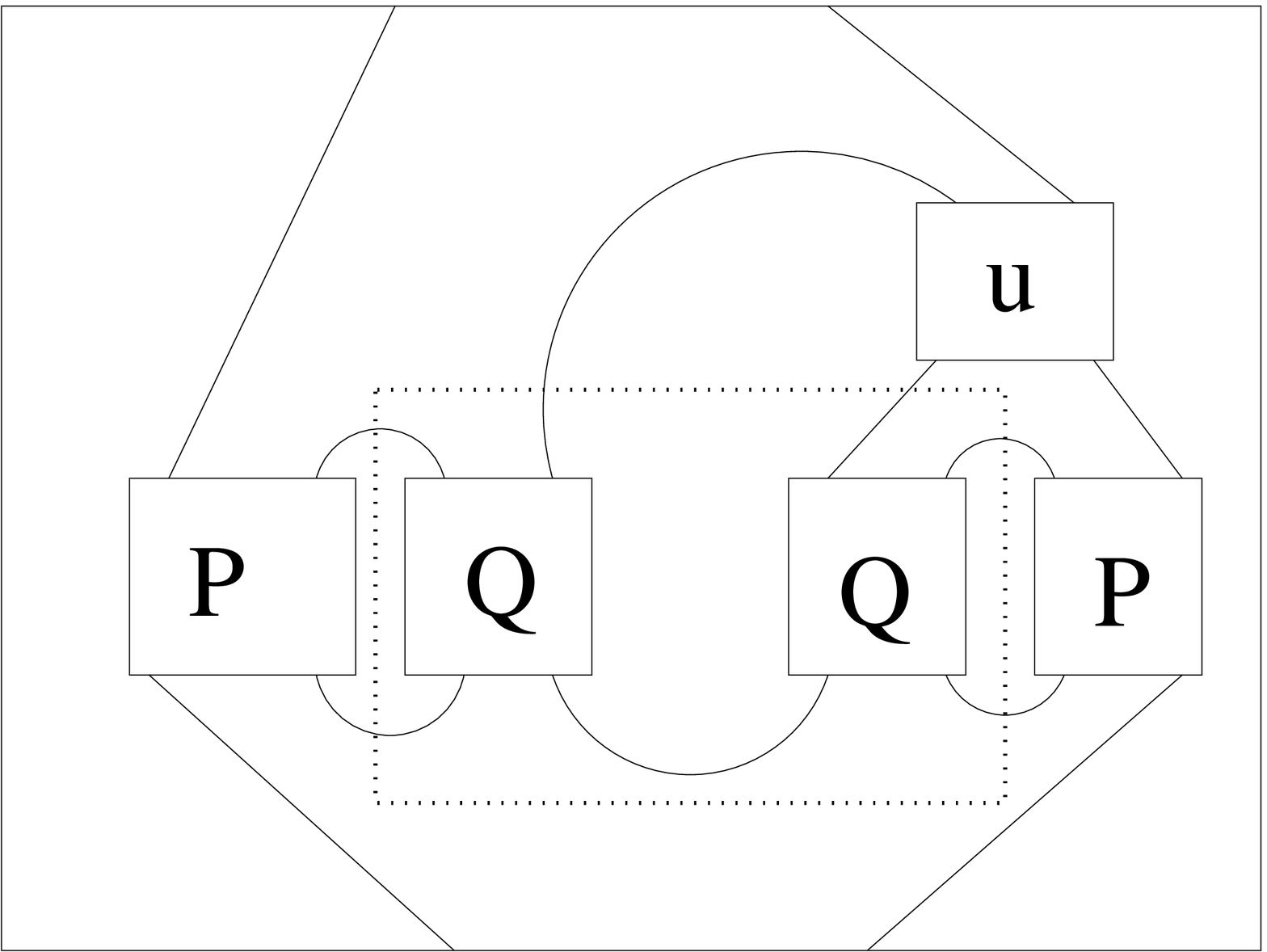} {1.5in}   }.
\\
The point of using \ref{switchPQQP} is to ensure that in the dotted rectangle we see
either two $P$'s or two $Q$'s. The $u$ may thus end up below the $P$'s and $Q$'s but
that does not affect the rest of the argument. In the dotted rectangle we may thus
apply the exchange relation for $Q$ to obtain, after a little isotopy:
\hbox{ \qquad         \vpic{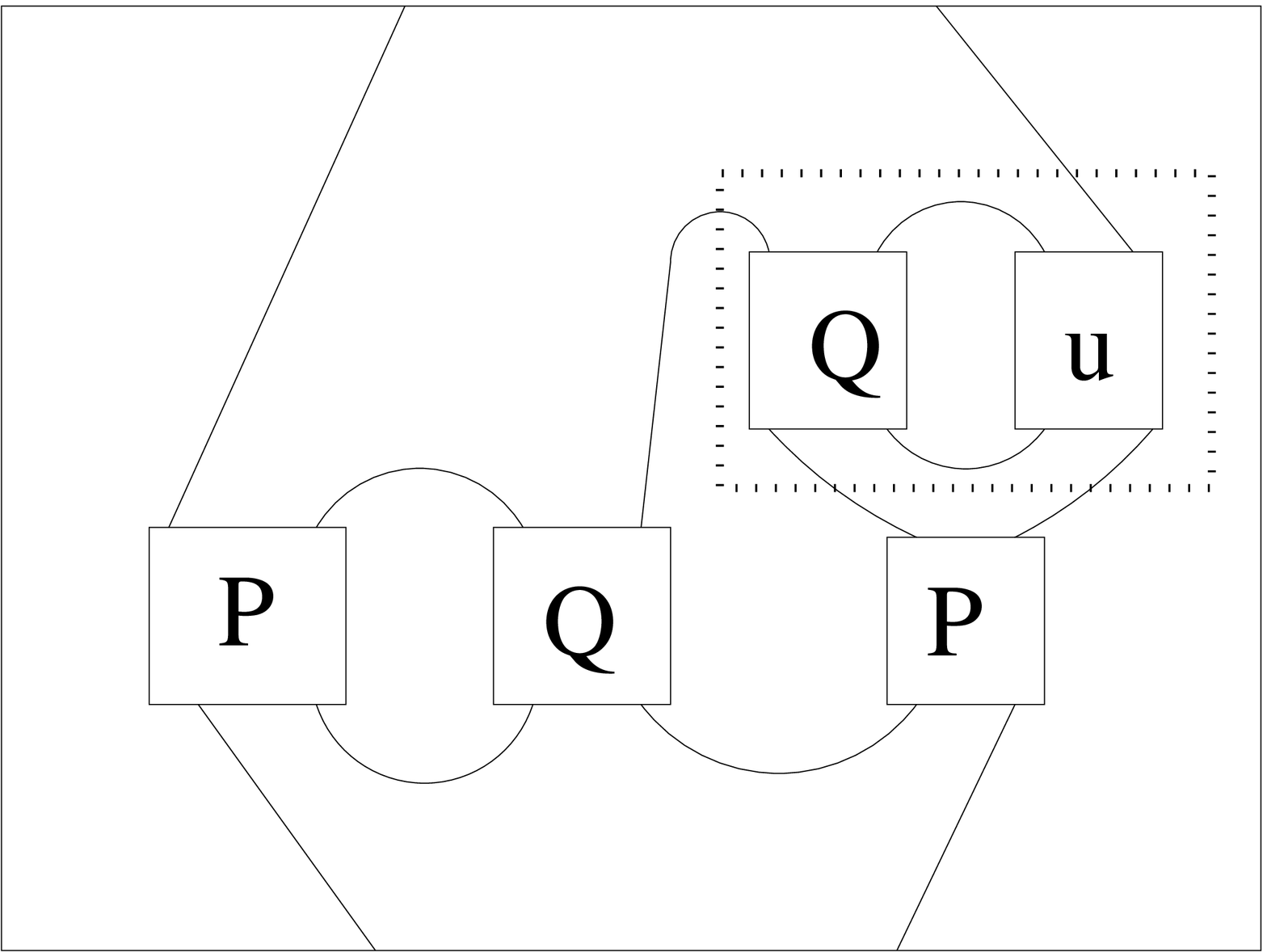} {1.5in}   }. Notice that inside the dotted rectangle
we see the comultiplication of $e_Q$ and $u$. Replacing $u$ by $2e_R -1$ gives 2 tangles,
the one with the identity being $\approx e_P\circ e_Q \circ e_P$. The tangle with $e_R$ can be handled easily since $e_Q \circ e_R =1$
which also yields $e_P\circ e_Q \circ e_P$.\\
Finally we need to be able to comultiply $ue_{PQ}$ with itself and $ue_{QP}$. This
goes very much like comultiplying $ue_P$ and $ue_Q$ except that after applying
\ref{reduceu} we find  a 
coproduct of more than two terms on $e_P$ and $e_Q$. These words may be reduced
to $e_P$, $e_Q$, $e_{PQ}$ or $e_{QP}$ modulo $e_N$. The term with $e_N$ will produce
a $u$ by itself but as observed above we know how to write $u$ as an explicit
linear combination of basis elements.

Case 3. Comultiplication between {\gothic A} and {\gothic B}. Terms without
$u$ like $e_P\circ e_{PQ}$ are simple.
The most difficult case is of the form $e_P \circ ue_{QP}$ but as above we
may rearrange it so that there are two like terms in the dashed rectangle below:\\

\hbox{ \qquad \qquad       \qquad \qquad       \vpic{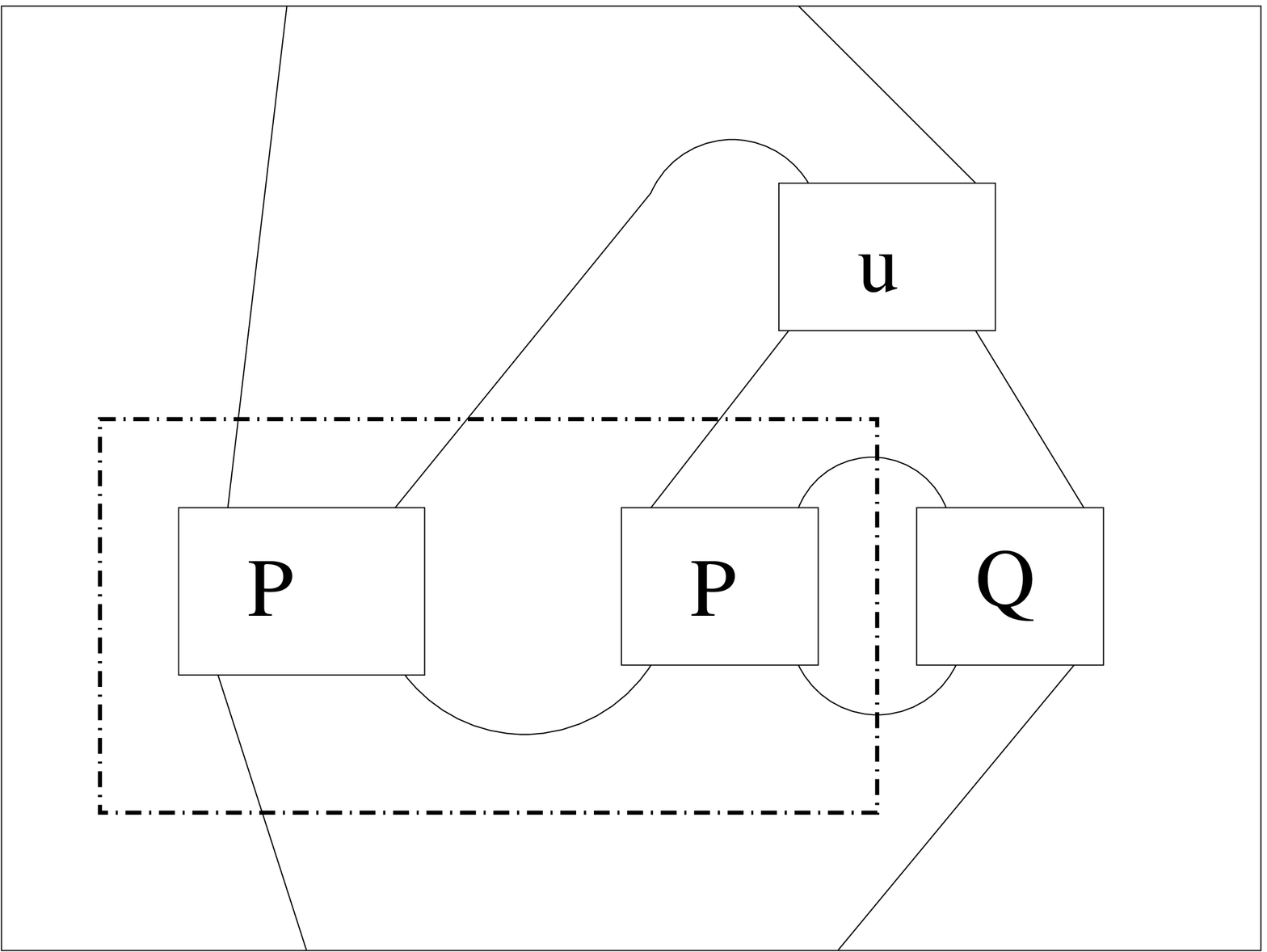} {1.5in}  .   }

Applying the exchange relation as before we obtain:

\hbox{ \qquad \qquad       \qquad \qquad       \vpic{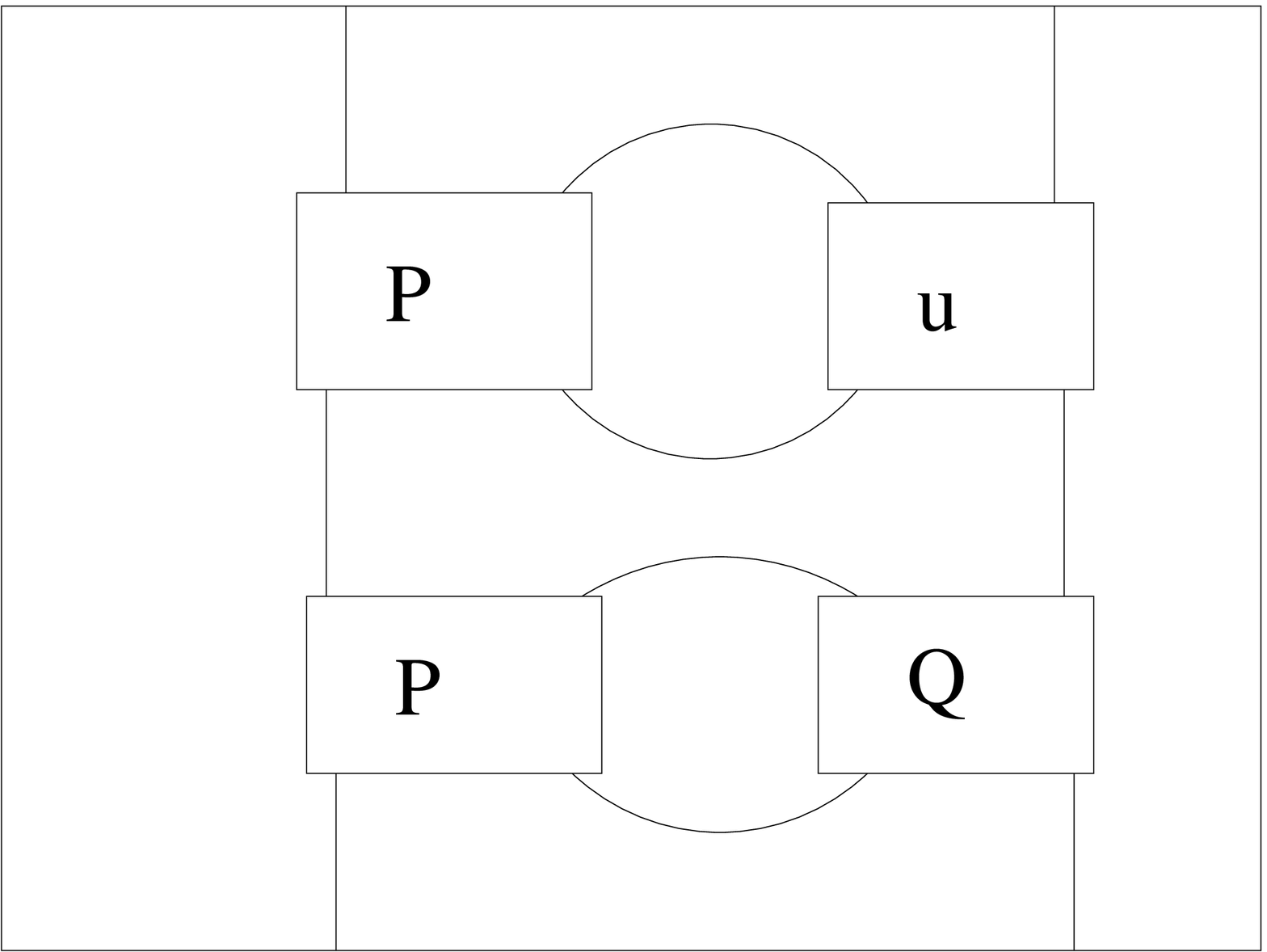} {1.5in} .  }

Note the comultiplication of $u$ and $e_P$ which can be reduced to an explicit
linear combination of basis elements using $u=2e_R-1$ and $e_R\circ e_P \approx 1$.

The coproduct of $ue_P$ with $e_{PQ}$ works similarly except that applying
the exchange relation immediately produces an explicit multiple of a basis 
element. Finally terms like $ue_P \circ ue_{PQ}$ can be reduced to explicit
linear combinations of basis elements using \ref{reduceu} and comulitplication
of words on $e_P$ and $e_Q$. Once again $u$ terms may be produced.
\end{proof}

\begin{lemma}
Let $v\in M'\cap M_2$ be the self-adjoint unitary in the normaliser of
$M_1$ guaranteed by the form of the dual principal graph in \ref{mainchance}.
Then $vAv=B$.
\end{lemma}
\begin{proof} By \ref{actionnormaliser} we have $ve_Pv=$ \vpic{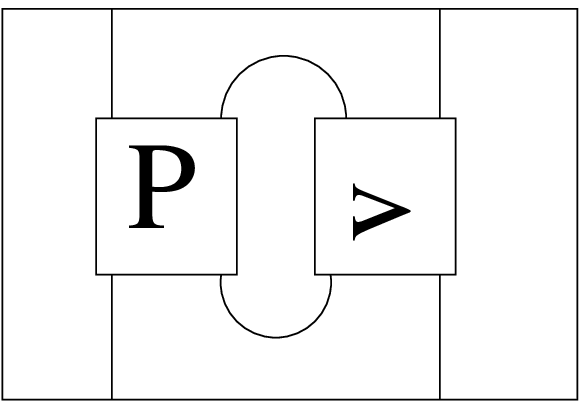} {1.0in} . So
$e_Pve_Pv=$ \vpic{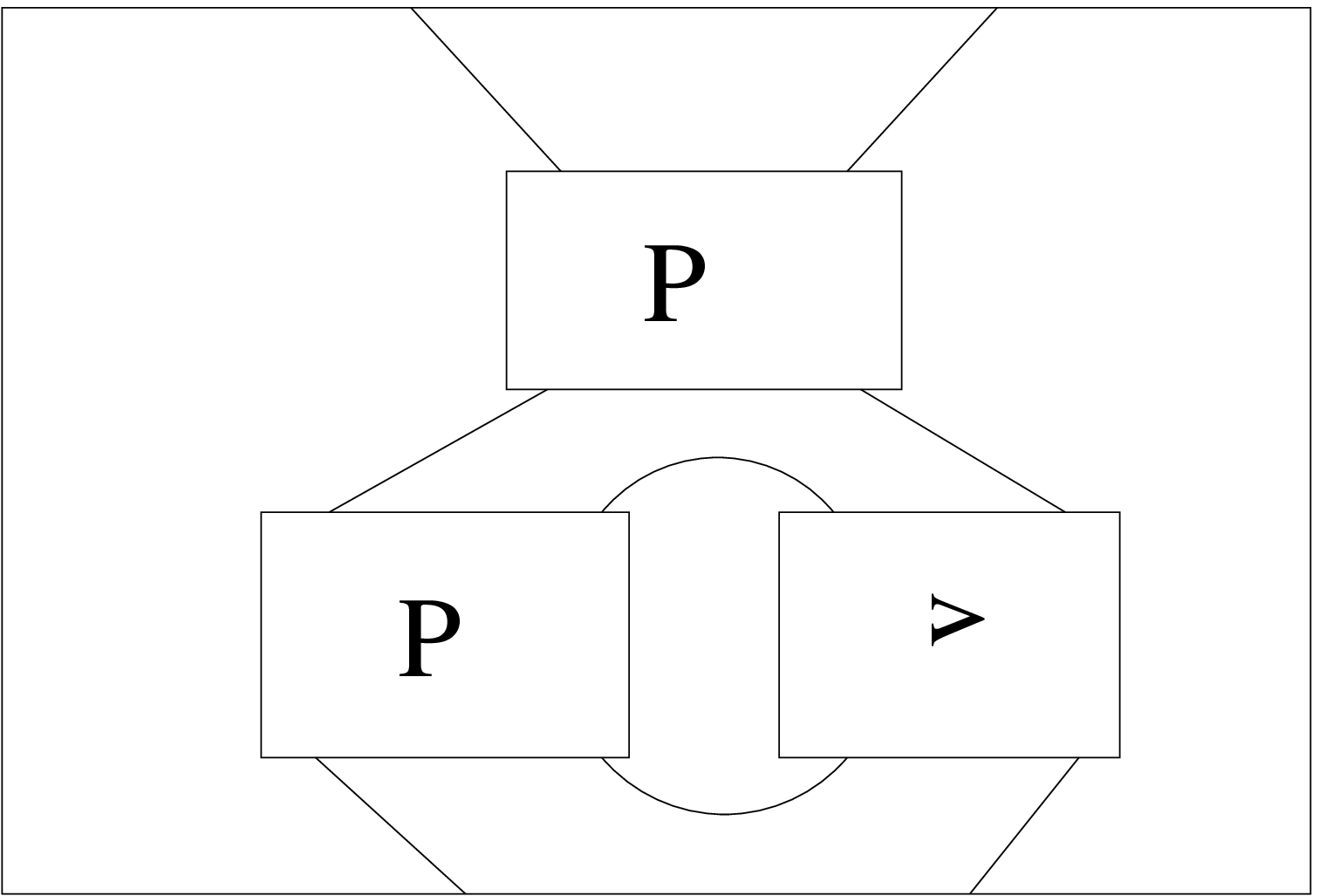} {1.0in} . Applying the exchange relation to this we obtain
\vpic{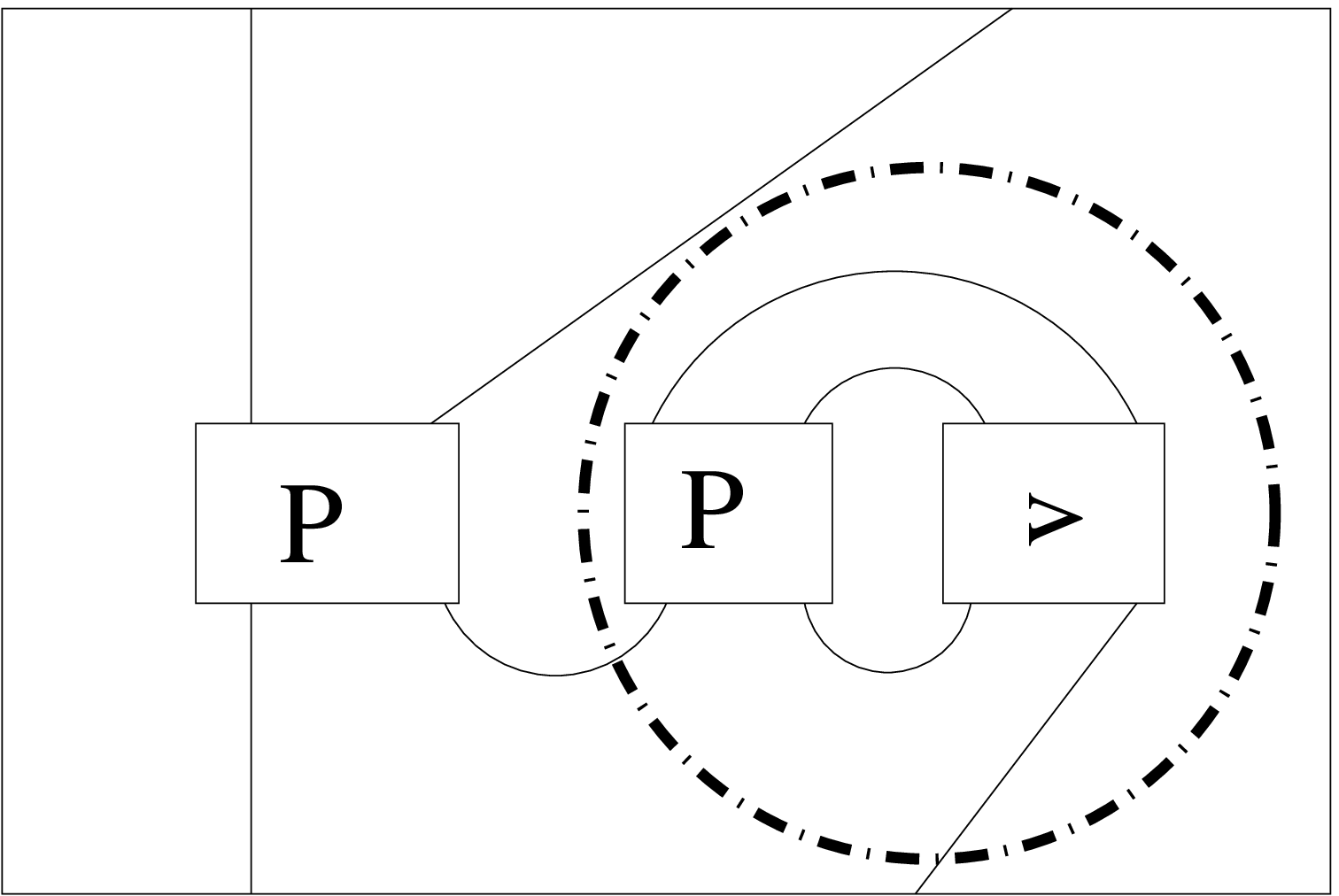} {1in} . Inside the dotted circle we recognise a multiple of the trace
in $M_2$ of the product in  $M'\cap M_2$ of the projection $e_{\overline P}$ defined by $e_P$ and $v$. 
But $v$ bears the same relation to this coprojection as $u$ does to $e_P$ so by
\ref{orthogonality of u} we obtain zero. Thus $e_PVe_PV=0$. We may apply 
\ref{switchPQQP} to $e_{\overline P}$ and $v$ to deduce in the same way that 
$e_QVe_PV=0$. This is enough to conclude that $vAv=B$ from the structure of $N'\cap M_1$
which is normalised by $v$.
\end{proof}

\begin{corollary} \label{3boxspace}If $e_M$ is the projection onto $L^2(M)$ in the basic construction
of $M_2$ from $M_1$ then
$\mathfrak D =${\gothic C}$e_M${\gothic C} $\cup$ {\gothic A} $\cup$
 $v${\gothic A} $\cup$ {\gothic B} $\cup v${\gothic B}
is a basis for $N'\cap M_2$. 
\end{corollary}
\begin{proof}
From the principal graph, $N'\cap M_2$ is the direct sum of the ideal $\mathfrak I$
generated 
by $e_M$, which is isomorphic
to a basic construction coming from the pair $N'\cap M \subseteq N'\cap M_1$,
and a $4\times 4$ matrix algebra. Since $N\subseteq M$ is irreducible the 
map $x\otimes y \mapsto xe_My$ is a vector space isomorphism from 
$N'\cap M_1\otimes N'\cap M_1$ to $\mathfrak I$. Thus 
{\gothic C}$e_M${\gothic C}  is a basis for $\mathfrak I$.

Since $v$ is in the normaliser of $M_1$, it is orthogonal to $M_1$ by
irreducibility and $N'\cap M_2$ contains a copy  of the crossed product
of $N'\cap M_1$ by the period 2 autormorphism given by $Ad \hspace{3pt}v$. 
By the previous lemma the algebra generated by $A$, $B$, and $v$ is a
$4\times 4$ matrix algebra-call it $\mathfrak E$. It is spanned modulo 
$\mathfrak I$ by 
{\gothic A} $\cup$
 $v${\gothic A} $\cup$ {\gothic B} $\cup v${\gothic B}
 since $A$ and $B$ are spanned modulo
$e_N$ by {\gothic A} and {\gothic B} respectively (see the proof of \ref{basis1}).
Since a matrix algebra is simple, to check that $\mathfrak E$ spans 
$N'\cap M_2$ mod $\mathfrak I$ we need only show that it is not contained in
$\mathfrak I$.  But from the principal graph we see that $A$ itself is non-zero
mod $\mathfrak I$.

\end{proof}

\subsection{The uniqueness proof and some corollaries.}

We can now give the main argument for the uniqueness of a subfactor of index $(2+\sqrt 2)^2$
with noncommuting intermediate subfactors. It relies on the "exchange relation"
developed by Landau in \cite{La2}. We begin with a planar algebra result from which
our uniqueness will follow. 

NOTE: We will assume that all planar algebras $P$ satisfy $\dim P_1 =1$.

\begin{definition}Let $P= P_n$ be a planar algebra and $\mathfrak R$ a self-adjoint subset
of $P_2$. Let $\mathfrak Y$ be the set of planar $3$-tangles
labelled with elements of $\mathfrak R$, with at most one internal
disc. We say that $\mathfrak R$ satisfies an {\rm exchange relation} if
there are constants $b_{Q,R,Y}$, $c_{Q,R,S,T}$ and $d_{Q,R,S,T}$ such that\\

\noindent \vpic{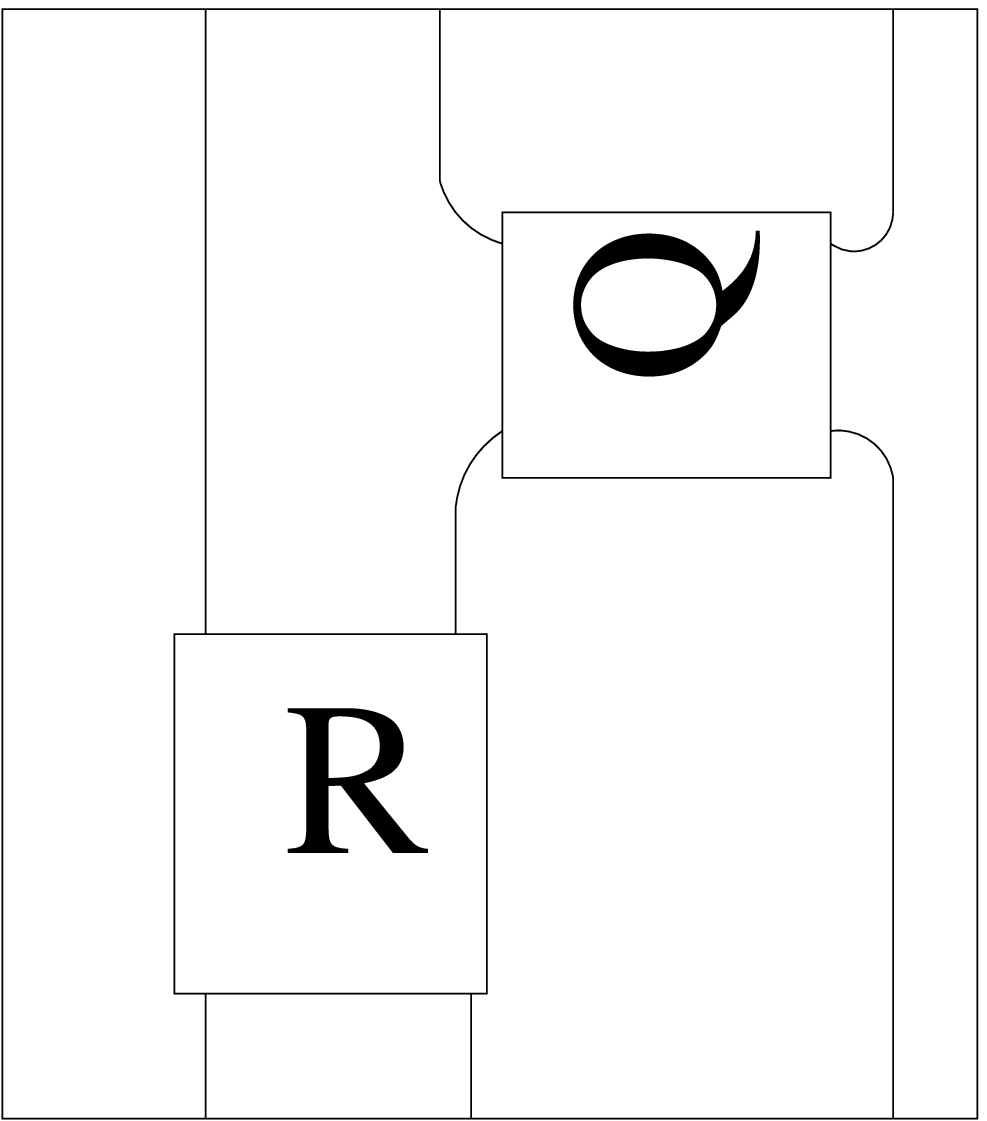} {0.6in} $\displaystyle =\sum_{S,T \in \mathfrak R} c_{Q,R,S,T}$ 
\vpic{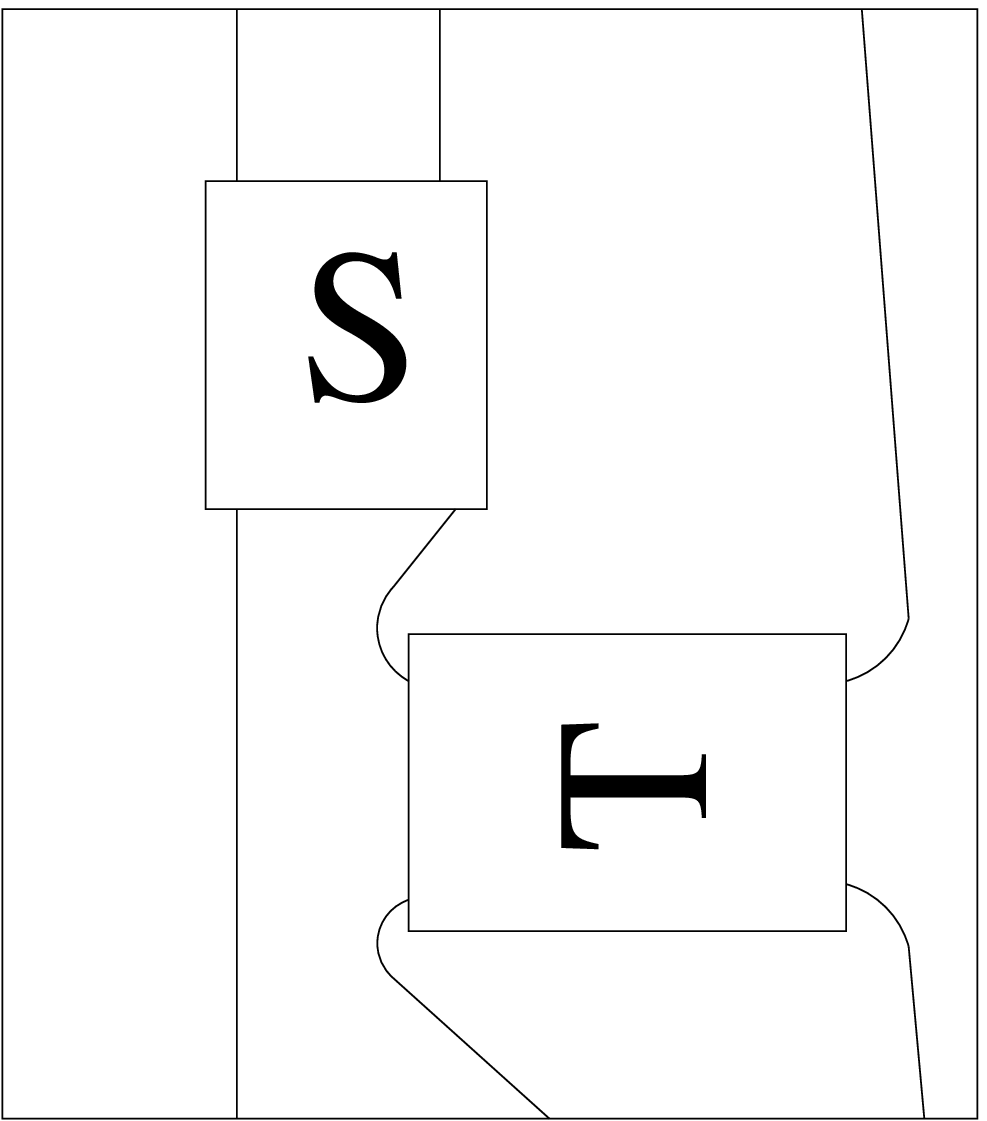} {0.6in} $\displaystyle +\sum_{S,T \in \mathfrak R} d_{Q,R,S,T}$
\vpic{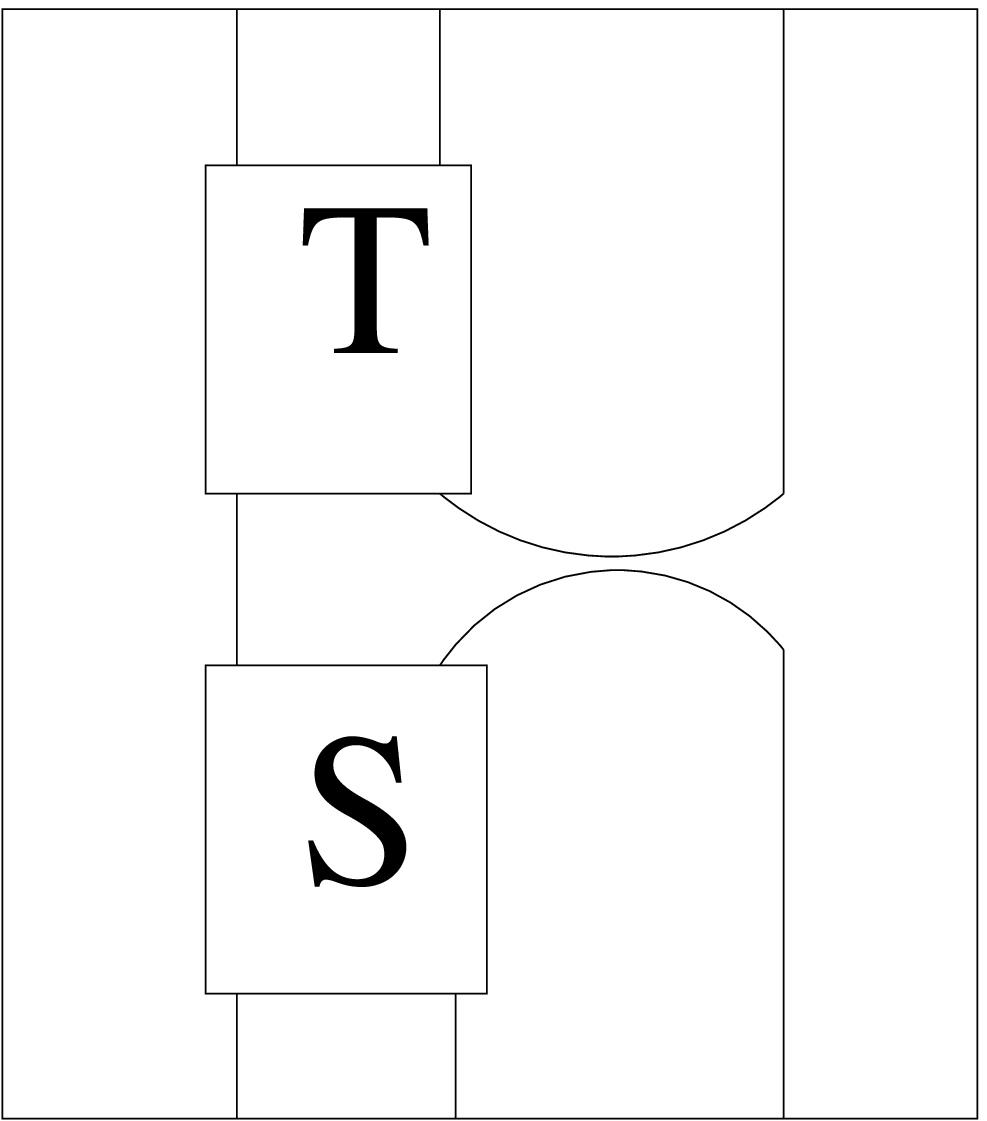} {0.6in} $\displaystyle +\sum_{Y\in \mathfrak Y} b_{Q,R,Y} Y.$

The constants will be called the {\rm exchange constants} for $\mathfrak R$.
\end{definition}

\begin{theorem}\label{zeph}(Landau,\cite{La2})
A subfactor planar algebra $P$ generated by $\mathfrak R =\mathfrak R ^*\subseteq P_2$ is determined
up to isomorphism by the exchange constants for $\mathfrak R$ and the traces and
cotraces of elements in $\mathfrak R$. 
\end{theorem}

The idea of the proof is that one may calculate the partition function of
any labelled tangle in $P_0$ by applying the exchange relation. The strategy
is to take any face and reduce it to a bigon, which is either a multiplication
or comultiplication of elements in $\mathfrak R$. But multiplication and
comultiplication are also determined by the exchange relation by suitably
capping off the pictures in the above definition. As soon as the planar 
algebras in question are non-degenerate in the sense that they are determined
by the partition functions of labelled planar tangles in $P_0$,
the theorem will hold. The isomorphism between two planar
algebras with the same subset $\mathfrak R$ is defined by extending
the identity map from $\mathfrak R$ to itself to all labelled tangles on
$\mathfrak R$. Then any relation for one planar algebra is necessarily
a relation for the other by nondegeneratess of the partition function as 
a bilinear/sesquilinear form on the $P_n$. 
This strategy for proving uniqueness was already used
for a proof of the uniqueness of the $E_6$ and $E_8$ subfactors in \cite{J21}.

\begin{lemma}\label{determined} Let $P$ be a subfactor planar algebra with $\mathfrak R$
a self-adjoint subset of $P_2$ which satisfies an exchange relation. Then the exchange
constants for $\mathfrak R$ are determined by the traces and cotraces of
elements of $\mathfrak R$ together with the structure constants for
multiplication and comultiplication of elements of $\mathfrak R$.
\end{lemma}

\begin{proof}
Using positive definiteness 
of the inner product given by the trace on $P_3$, it suffices to 
prove that the partition function of any planar diagram with at most
4 internal discs, all labelled with elements of $\mathfrak R$,  is determined
by the given structure constants.

For this, we may suppose the labelled diagrams are connected and by our
hypothesis on $\dim P_1$, we may suppose that no $2-box$ is connected to 
itself. If there
are 4 internal discs one must be connected to another with a multiplication
or a comultiplication. This reduces us to the case of 3 internal boxes where
it is even clearer. To see these assertions it is helpful to view the
labelled tangles as the generic planar projections of links in $\mathbb R^3$ which 
are obtained by shrinking the internal 2-boxes to points.
\end{proof}

Putting the previous results together we have:
\begin{theorem}\label{mainuniqueness}
 Let $N_1\subseteq M_1$ and $N_2\subseteq M_2$ be two irreducible
 II$_1$ subfactors
of index $(2+\sqrt 2)^2$ with pairs $P_1$, $Q_1$ and $P_2$, $Q_2$ of non-commuting
intermediate subfactors of index $2+\sqrt 2$. Then there is a unique isomorphism 
from the planar algebra for $N_1\subseteq M_1$ to the planar algebra of
 $N_2\subseteq M_2$ which extends
the map sending $e_{P_1}$ and $e_{Q_1}$ to $e_{P_2}$ and $e_{Q_2}$ respectively.
\end{theorem}

\begin{proof}  The only allowed principal graph for the elementary 
subfactors is $A_7$ so there is no extra structure and we know the 
principal graph and dual principal graph. The normalising unitaries $u_i , i=1,2$
 can be written
as an explicit  linear combination of $e_N, 1$ and products and coproducts of
$e_{P_i}$ and $e_{Q_i}$. Then form the sets {\gothic A}$_i$ and {\gothic B}$_i$,
$i=1,2$ in the obvious way. The planar algebra for $N_i\subseteq M_i$ is generated
by {\gothic A}$_i$ and {\gothic B}$_i$ by \ref{3boxspace}. By \ref{basis1} and
\ref{determined}  we may apply  \ref{zeph} to the sets $\mathfrak R_i =${\gothic A}$_i  \cup${\gothic B}$_i$ to
deduce the result. (The traces and cotraces of the basis elements 
of {\gothic C} were determined in the course of proving \ref{basis1}.) 
\end{proof}
 \begin{corollary}
Given a quadrilateral $N\subset P, Q \subset M$ with $[M:N]=6+4\sqrt 2$ and such that $P$ and
$Q$ do not commute, there are further subfactors $\tilde P$ and $\tilde Q$ with $[M:\tilde P]$
and $[M:\tilde Q]$ equal to $2+\sqrt 2$, which
commute with both $P$ and $Q$   and are at an angle $\cos^{-1}(\sqrt 2 -1)$ to each other.

\end{corollary}
\begin{proof} This is the case for the example so by uniqueness it is always true.
\end{proof}
It is obvious that the projections onto $\tilde P$ and $\tilde Q$ are in $B$ mod $e_N$.
\begin{corollary} The only subfactors between $N$ and $M$ are $P,Q,\tilde P,\tilde Q, R$ and $S$
so the intermediate subfactor lattice is
\vpic{d5lattice} {1.3in} .
\end{corollary}
\begin{proof}
Let $T$ be a seventh intermediate subfactor. From the principal graph and obvious
index restrictions the possible values of $(6+4\sqrt 2)tr(e_T)$ are $2+\sqrt 2$, $3+2\sqrt 2$
and $2$. The cases $3+2\sqrt 2$ and $2$ correspond to index 2 subfactors and would show up
as extra vertices on either the dual or dual principal graphs, so we must have 
$\displaystyle tr(e_T)= \frac{1}{2+\sqrt 2}$. This forces $e_T-e_N$ to be a minimal projection
in either $A$ or $B$, so by the previous corollary and the 
observation after it we may suppose wolog that 
$e_T-e_N\in A$. If $e_P e_T=e_N$ then by a $2\times 2$ matrix calculation $T$ makes 
a forbidden angle with $Q$. So the angle between all three of $P$, $Q$ and $T$ is $\cos^{-1}(\sqrt 2-1)$.
But by lemma \ref{position of R} applied to $T$, $T$ and $R$ do not commute
so there must be a fourth subfactor $\alpha(T)$ which makes the same angle with all
of $P,Q$, and $T$. By lemma \ref{four subspaces} this is not allowed. This contradicts
the existence of $T$.
\end{proof}
\begin{corollary} If $M$ is hyperfinite there
 is an automorphism of $M$ sending  $P$ to $\tilde P$ and 
$Q$ to  $\tilde Q$.
\end{corollary} 
\begin{proof}This follows from theorem \ref{mainuniqueness} and Popa's
classification theorem \cite{P6} which states that in finite depth one may construct
the subfactor directly as the completion of the inductive limit of the tower
of relative commutants. 
\end{proof}
It is not obvious what the automorphism of the previous corollary looks like
in the GHJ realisation of section \ref{the example}. It will certainly require
the complex numbers to write it down as guaranteed by the next result.
Observe first that the $D_5$-based GHJ example of \ref{simpler} is defined
over the real numbers so the intermediate subfactors exist in the setting
of real II$_1$ factors. That the GHJ pair for $D_{5,2}$ needs the complex
numbers is the next result.
 
\begin{corollary}If $N\subset P,Q \subset M$ is a noncommuting quadrilateral of real 
II$_1$ factors with $[M:N]=6+4\sqrt 2$, then $P$ and $Q$ are the \underline{only} intermediate
subfactors of index $2+\sqrt 2$.
\end{corollary}
\begin{proof} Let $N\subset M$ be the subfactor for the $D_{5,2}$
Coxeter graph. Since this subfactor may be defined over the reals (as the
GHJ subfactor for the trivalent vertex) complex conjugation defines a
conjugate linear *-automorphism $\sigma$ of $N\subset M$ with $\sigma (\tilde P)=\tilde P$
and $\sigma (\tilde Q)=\tilde Q$ but with 
$\sigma(g_i)=g_i^*$ so $\sigma (P)=Q$. Thus $\sigma$ will act on the planar
algebra of $N\subset M$ exchanging $e_P$ and $e_Q$. However the fixed points for
$\sigma$ acting on the planar algebra is again a planar algebra so there is 
a real subfactor $N_{\mathbb R}\subset M_{\mathbb R}$ with 
$[M_{\mathbb R}:N_{\mathbb R}]=6+4\sqrt 2$ having a pair ($\tilde P^\sigma$
and $\tilde Q^\sigma$) of noncommuting
intermediate subfactors of index $2+\sqrt 2$ and no other intermediate subfactors
of the same index since $\sigma(e_P)=e_Q\neq e_P$. Our uniqueness result never
used the complex numbers (all the structure constants were real) so that
no other such real subfactor can have more than two intermediate subfactors
of index $2+\sqrt 2$.
\end{proof}

\5
%\hspace{1in}\vpic{fffr} {2.5in}
\5
%\subsection{The principal and dual principal graphs.}
%In the case of index $(2+\sqrt 2)^2$ we know that $N'\cap M_2$ is
%the direct sum $\mathbb C \oplus M_2(\mathbb C) \oplus M_2(\mathbb C) \oplus \mathbb C$.
%The subfactors dual to $P$ and $Q$ satisfy the same hypotheses as $P$ and $Q$ 
%so $M'\cap M1$ is also  $\mathbb C \oplus M_2(\mathbb C) \oplus M_2(\mathbb C) \oplus \mathbb C$.
\end{document}